\documentclass[11pt,notitlepage,twoside,a4paper]{amsart}
 \usepackage{amsfonts}

\usepackage{amsmath,amssymb,enumerate}

\usepackage{epsfig,fancyhdr,color}%,showkeys,amsmidx
 
\usepackage[table]{xcolor}
 
\usepackage{amssymb}
\usepackage{amsmath,amsthm}
\usepackage{latexsym}
\usepackage{amscd}
\usepackage{psfrag}
\usepackage{graphicx}
\usepackage[latin1]{inputenc}
\usepackage[all]{xy}
\usepackage[mathcal]{eucal}

\definecolor{NoteColor}{rgb}{1,0,0}

\renewcommand{\textsc}{\textcolor{red}}

\newtheorem*{theorem 1}{\rm\bf Proposition 1}
\newtheorem*{theorem 2}{\rm\bf Proposition 2}

\theoremstyle{definition}

\theoremstyle{remark}

\def\interieur#1{\mathord{\mathop{\kern 0pt #1}\limits^\circ}}

\title[Looking backward: From Euler to Riemann]{Looking backward:  From Euler to Riemann}

\author{Athanase Papadopoulos}
\address{A. Papadopoulos, Institut de Recherche Math{\'e}matique Avanc\'ee,
Universit{\'e} de Strasbourg and CNRS,
7 rue Ren\'e Descartes,
 67084 Strasbourg Cedex, France, and 
  Brown University, Mathematics Department, 
 151 Thayer Street
Providence, RI 02912, USA.}

\makeindex

 \date{\today}

\begin{document}

\maketitle

  \hfill{\emph{\small{Il est des hommes auxquels on ne doit pas adresser d'\'eloges,}}}
  
  \hfill{\emph{\small{si l'on ne suppose pas qu'ils ont le go\^ut assez peu d\'elicat}}}
  
  \hfill{\emph{\small{pour go\^uter les louanges qui viennent d'en bas}.}}
  
    \hfill{\small{(Jules Tannery, \cite{Tannery-souvenir} p. 102)}}

\begin{abstract}  

We survey the main ideas in the early history of the subjects on which Riemann worked and that led to some of his most important discoveries. The subjects discussed include the theory of functions of a complex variable, elliptic and Abelian integrals, the hypergeometric series, the zeta function, topology, differential geometry,  integration, and the notion of space. We shall see that among Riemann's predecessors in all these fields, one name occupies a prominent place, this is Leonhard Euler. The final version of this paper will appear in the book \emph{From Riemann to differential geometry and relativity} (L. Ji, A. Papadopoulos and S. Yamada, ed.)  Berlin: Springer, 2017.

\end{abstract}
 \medskip


\medskip


\tableofcontents

 \section{Introduction} \label{s:intro}
 
 More than any other branch of knowledge, mathematics is a science in which every generation builds on the accomplishments of the preceding ones, and where reading the old masters has always been a ferment for new discoveries. Examining the roots of Riemann's ideas takes us into the history of complex analysis, topology, integration, differential geometry and other mathematical fields, not to speak of physics and philosophy, two domains in which Riemann was also the heir of a long tradition of scholarship.
 
Riemann himself was aware of the classical  mathematical literature, and he often quoted his predecessors. For instance, in the last part of his Habilitation lecture, \emph{\"Uber die Hypothesen, welche der Geometrie zu Grunde liegen}\index{Riemann! habilitation lecture}\index{habilitation lecture!Riemann} \cite{Riemann-Ueber} (1854), he writes:\footnote{In all this paper, for Riemann's habilitation, we use Clifford's translation \cite{Riemann-Clifford}.}  
\begin{quote}\small The progress of recent centuries in the knowledge of mechanics depends almost entirely on the exactness of the construction which has become possible through the invention of the infinitesimal calculus, and through the simple principles discovered by Archimedes,\index{Archimedes (c. 287 B.C.--c. 212 B.C.)} Galileo\index{Galilei, Galileo (1564--1642)} and Newton,\index{Newton, Isaac (1643--1727)} and used by modern physics.
\end{quote}

 The references are eloquent: Archimedes,\index{Archimedes (c. 287 B.C.--c. 212 B.C.)}  who developed the first differential calculus, with his computations of length, area and volume, Galileo,\index{Galilei, Galileo (1564--1642)} who introduced the modern notions of motion, velocity and acceleration, and Newton,\index{Descartes, Ren\'e (1596--1650)} who was the first to give a mathematical expression to the forces of nature, describing in particular the motion of bodies in resisting media, and most of all, to whom is attributed a celebrated notion of space, the ``Newtonian space." As a matter of fact, the subject of Riemann's habilitation lecture includes the three domains of Newton's  \emph{Principia}:\index{Newton!Principia} mathematics, physics and philosophy. It is interesting to note also that Archimedes,\index{Archimedes (c. 287 B.C.--c. 212 B.C.)} Galileo\index{Galilei, Galileo (1564--1642)} and Newton\index{Newton, Isaac (1643--1727)} are mentioned as the three founders of mechanics in the introduction (Discours pr\'eliminaire) of Fourier's\index{Fourier, Joseph (1768--1830)} \emph{Th\'eorie analytique de la chaleur}\index{Fourier!Th\'eorie analytique de la chaleur} (\cite{Fourier}, p. i--ii), a work in which the latter lays down the rigorous foundations of the theory of trigonometric series. Fourier's quote and its English translation are given in \S \ref{s:trigo} of the present paper.  In the historical part of his Habilitation dissertation, \emph{\"Uber die Darstellbarkeit einer Function durch eine trigonometrische Reihe} (On the representability of a function by a trigonometric series)\index{trigonometric series} \cite{Riemann-Trigo}, a memoir which precisely concerns trigonometric series, Riemann gives a detailed presentation of the history of the subject, reporting on  results and conjectures by Euler,\index{Euler, Leonhard (1707--1783)} d'Alembert,\index{Alembert@d'Alembert, Jean le Rond (1717--1783)} Lagrange,\index{Lagrange, Joseph-Louis (1736--1813)} Daniel Bernoulli, \index{Bernoulli, Daniel (1700--1782)}  Dirichlet\index{Dirichlet, Johann Peter Gustav Lejeune (1805--1859)},  Fourier and  others. The care with which Riemann analyses the evolution of this field, and the wealth of historical details he gives, is another indication of the fact that he valued to a high degree the history of ideas and was aware of the first developments of the subjects he worked on. In the field of trigonometric series and in others, he was familiar with the important paths and sometimes the wrong tracks that his predecessors took for the solutions of the problems he tackled. Riemann's  sense of history is also manifest in the announcement of his memoir  \emph{Beitr\"age zur Theorie der durch die Gauss'sche Reihe $F(\alpha,\beta,\gamma,x)$ darstellbaren Functionen}  (Contribution to the theory of functions representable by Gauss's series $F(\alpha,\beta,\gamma,x)$),  published in the \emph{G\"ottinger Nachrichten}, No. 1, 1857, in which he explains the origin of the problems considered, mentioning works of Wallis, Euler, Pfaff, Gauss\index{Gauss, Carl Friedrich (1777--1855)} and  Kummer.  There are many other examples.
  
Among Riemann's forerunners in all the fields that we discuss in this paper, one man fills  almost all the background; this is Leonhard Euler.\index{Euler, Leonhard (1707--1783)} 
Riemann was an heir of Euler for what concerns functions of a complex variable, elliptic integrals, the zeta function,\index{zeta function}\index{Riemann!zeta function} the topology of surfaces, the differential geometry of surfaces, the calculus of variations, and several topics in physics. 

Riemann refers to Euler at several places of his work, and Euler\index{Euler, Leonhard (1707--1783)} was himself a diligent reader of the classical literature: Euclid, Pappus, Diophantus, Theodosius, Descartes,\index{Descartes, Ren\'e (1596--1650)} Fermat, Newton,\index{Newton, Isaac (1643--1727)} etc. All these authors are mentioned all along his writings, and many of Euler's works were motivated by questions that grew out of his reading of them.\footnote{Cf. for instance Euler's \emph{Problematis cuiusdam Pappi Alexandrini constructio} (On a problem posed by Pappus of Alexandria), \cite{E543}, 1780.}
Before going into more details, I would like to say a few words about the lives of Euler and Riemann, highlighting analogies, but also differences between them.

 Both Euler\index{Euler, Leonhard (1707--1783)} and Riemann\index{Riemann, Bernhard  (1826--1866)} received their early education at home, from their fathers, who were protestant ministers, and who both were hoping that their sons will become like them, pastors. At the age of fourteen, Euler\index{Euler, Leonhard (1707--1783)} attended a Gymnasium in Basel, while his parents lived in Riehen, a village near the city of Basel.\footnote{Today, Riehen belongs to the Canton of the city of Basel, and it hosts the famous Beyeler foundation.} At about the same age, Riemann was sent to a Gymnasium in Hanover, away from his parents. During their Gymnasium years, both Euler and Riemann lived with their grandmothers.\footnote{In 1842, at the death of his grandmother, Riemann quitted Hanover and attended the Gymnasium at the Johanneum L\"uneburg.} They both enrolled a theological curriculum (at the Universities of Basel and G\"ottingen respectively), before they obtain their fathers' approval to shift to mathematics. 
 
 There are also major differences between the lives of the two men.\index{Euler, Leonhard (1707--1783)}\index{Riemann, Bernhard  (1826--1866)}  Euler's productive period lasted 57 years (from the age of 19, when he wrote his first paper, until his death at  the age of 76). His written production comprises more than 800 memoirs and 50 books. He worked on all domains of mathematics (pure and applied)  and physics (theoretical and practical) that existed at his epoch. He also published on geography, navigation, machine theory, ship building, telescopes, the making of optical instruments, philosophy, theology and music theory. Besides his research books, he wrote elementary schoolbooks, including a well-known book on the art of reckoning \cite{E17}. The publication of his collected works was decided in 1907, the year of his bicentenary, the first volumes appeared in 1911, and the edition is still in progress (two volumes appeared in 2015), filling up to now more than 80 large volumes. Unlike Euler's, Riemann's\index{Riemann, Bernhard  (1826--1866)} life was short. He published his first work at the age of 25 and he died at the age of 39. Thus, his productive period lasted only 15 years. His collected works stand in a single slim volume. Yet, from the points of view of the originality and the impact of their ideas, it would be unfair to affirm that either of them stands before the other. They both had an intimate and permanent relation to mathematics and to science in general. 
 Klein writes in his \emph{Development of mathematics in the 19th century} (\cite{Klein-development}, p. 231 of the English translation): 
 \begin{quote}\small
 After a quiet preparation Riemann came forward like a bright meteor, only to be extinguished soon afterwards. 
 \end{quote}
 
 On Euler, I would like to quote Andr\'e Weil,\index{Weil, Andr\'e (1906--1998)} from his book  on the history of number theory, \emph{Number Theory: An approach through history from Hammurapi to Legendre} \cite{Weil-Hammurapi}. He writes, in the concluding section:

   \begin{quote}\small
 [...] Hardly less striking is the fact that Euler never abandoned a problem after it has once aroused his insatiable curiosity. Other mathematicians, Hilbert for instance, have had their lives neatly divided into periods, each one devoted to a separate topic. Not so Euler. All his life, even after the loss of his eyesight, he seems to have carried in his head the whole of the mathematics of his day, both pure and applied. Once he has taken up a question, not only did he come back to it again and again, little caring if at times he was merely repeating himself, but also he loved to cast his net wider and wider with never failing enthusiasm, always expecting to uncover more and more mysteries, more and more ``herrliche proprietates" lurching just around the next corner. Nor did it matter to him whether he or another made the discovery. ``\emph{Penitus obstupui}", he writes (``I was quite flabberggasted": \emph{Eu}.I-21.1 in E 506$\vert$1777, cf. his last letter to Lagrange, \emph{Eu}.IV A-5.505$\vert$1775) on learning Lagrange's additions to his own work on elliptic integrals; after which he proceeds to improve upon Lagrange's achievement. Even when a problem seemed to have been solved to his own satisfaction (as happened with his first proof of Fermat's theorem $a^p\equiv a$ mod $p$, or in 1749 with sums of two squares) he never rested in his search for better proofs, ``more natural" (\emph{Eu}. I-2.510 in E 262$\vert$1755; cf. \S VI), ``easy" (\emph{Eu}.1-3.504 in E 522$\vert$1772; cf. \S VI), ``direct" (\emph{Eu}.I-2.363 in E 242$\vert$1751; cf. \S VI); and repeatedly he found them.\footnote{In Weil's book, every piece of historical information is accompanied by a precise reference. Works that attain this level of scholarship are very rare.}
\end{quote}

Let us say in conclusion that if we had to mention a single mathematician of the eighteenth century, Euler would probably be the right choice. For the nineteenth century, it would be Riemann. Gauss, who will also be mentioned many times in the present paper, is the main figure astride the two centuries.

Euler's results are contained in his published and posthumous writings, but also in his large correspondence, available in several volumes of his \emph{Opera Omnia}.  We shall mention several times this correspondence in the present paper. 
It may be useful to remind the reader that at the epoch we are considering here, there were very few mathematical journals (essentially the publications of the 
few existing Academies of Sciences). The transmission of open problems and results among mathematicians was done largely through correspondence. 
On this question, let us quote the mathematician Paul Heinrich Fuss,\index{Fuss, Paul Heinrich (1755--1826)} who published the first set of letters of Euler, and who was his great-grandson. He writes in the introduction to his \emph{Correspondence} \cite{Fuss-corresp}, p. xxv:\footnote{Unless otherwise stated, the translations from the French in this paper are mine.}
\begin{quote}\small 
Since sciences ceased to be the exclusive property of a small number of initiates, correspondence between scholars was taken over by the periodical publications. The progress is undeniable. However, this freeness with which ideas and discoveries were communicated in the past, in private and very confidential letters, we do not find it any more in the ripe and printed pieces of work. At that time, the life of a scholar was, in some way, all reflected in that correspondence. We see there the great discoveries being prepared and gradually developed; no link and no transition is missing;  the path which led to these discoveries is followed step by step, and we can draw there some information even in the errors committed by these great geniuses who were the authors. This is sufficient to explain the interest tied to this kind of correspondence.\footnote{Depuis que les sciences ont cess\'e d'\^etre la propri\'et\'e exclusive d'un petit nombre d'initi\'es, ce commerce \'epistolaire des savants a \'et\'e absorb\'e par la presse p\'eriodique. Le progr\`es est incontestable. Cependant, cet abandon avec lequel on se communiquait autrefois ses id\'ees et ses d\'ecouvertes, dans des lettres toutes confidentielles et priv\'ees, on ne le retrouve plus dans les pi\`eces m\^uries et imprim\'ees. Alors, la vie du savant se refl\'etait, pour ainsi dire, tout enti\`ere dans cette correspondance. On y voit les grandes d\'ecouvertes se pr\'eparer et se d\'evelopper graduellement ; pas un cha\^\i non, pas une transition n'y manque ; on suit pas \`a pas la marche qui a conduit \`a ces d\'ecouvertes, et l'on puise de l'instruction jusque dans les erreurs des grands g\'enies qui en furent les auteurs. Cela explique suffisamment l'int\'er\^et qui se rattache \`a ces sortes de correspondances.}
\end{quote}

In the case of Euler, particularly interesting is his correspondence with Christian Goldbach,\index{Goldbach, Christian (1690--1764)}  published recently in two volumes of the \emph{Opera Omnia} \cite{Euler-Goldbach}. It contains valuable information on Euler's motivations and progress in several of the domains that are surveyed in the sections that follow, in particular, topology, the theory of elliptic functions and the zeta function. A few lines of biography on this atypical person are in order. 

Christian Goldbach (1690--1764) was one of the first German scholars whom Euler met at the Saint Petersburg Academy of Sciences when he arrived there in 1827. He was very knowledgeable in mathematics, although he was interested in this field only in an amateurish fashion, encouraging others' works rather than working himself on specific problems. He was also a linguist and thoroughly involved in politics. Goldbach studied law at the University of K\"onigsberg.  In Russia, he became closely related to the Imperial family. In 1732, he was appointed secretary  of the Saint Petersburg Academy of Sciences and in 1737 he became the administrator of that institution. In 1740, he held an important position at the Russian ministry of foreign affairs and became the official cryptographer there. Goldbach\index{Goldbach, Christian (1690--1764)} had a tremendous influence on Euler\index{Euler, Leonhard (1707--1783)}, by being attentive to his progress, by the questions he asked him on number theory, and also by motivating him to read Diophantus and Fermat.  Goldbach,\index{Goldbach, Christian (1690--1764)} who was seventeen years older than Euler, became later on one of his closest friends\index{Euler, Leonhard (1707--1783)} and the godfather of his oldest son, Johann Albrecht, the only one among Euler's thirteen children who became a mathematician.  Paul Heinrich Fuss writes in the introduction to his \emph{Correspondence} \cite{Fuss-corresp}, p. xxii: 
\begin{quote}\small It is more than probable that if this intimate relationship between Euler and this scholar, a relationship that lasted 36 years without interruption, hadn't been there, then the science of numbers would have never attained the degree of perfection which it owes to the immortal discoveries of Euler.\footnote{Il me semble plus que probable que si cette liaison intime entre Euler et ce savant, liaison qui dura 36 ans sans interruption, n'e\^ut pas lieu, la science des nombres n'aurait gu\`ere atteint ce degr\'e de perfection dont elle est redevable aux immortelles d\'ecouvertes d'Euler.} 
\end{quote}
Goldbach kept a regular correspondence with Euler, Nicolas and Daniel Bernoulli, Leibniz\index{Leibniz, Gottfried Wilhelm (1646--1716)} (in particular on music theory) and many other mathematicians.

  After Goldbach and his influence on Euler, we turn to Gauss,\index{Gauss, Carl Friedrich (1777--1855)} who, among the large number of mathematicians with whom Riemann was in contact, was certainly the most influential on him.\footnote{Some historians of mathematics claimed that when Riemann enrolled the University of G\"ottingen, as a doctoral  student of Gauss, the latter was old and in poor health, and that furthermore, he disliked teaching. From this, they deduced that Gauss's influence on Riemann was limited. This is in contradiction with the scope and the variety of the mathematical ideas of Riemann for which he stated, in one way or another, but often explicitly, that he got them under the direct influence of Gauss  or by reading his works. The influence of a mathematician is not measured by the time spent talking with him or reading his works. Gauss died the year after Riemann obtained his habilitation, but his imprint on him was permanent.} 
We shall see in the various sections of the present paper that this influence was crucial for what concerns the fields of  complex analysis, elliptic integrals, topology, differential geometry -- the same list as for Euler's influence on Riemann --  and also for what concerns his ideas on space. There are other topics in mathematics and physics which were central in the work of Riemann and where he used ideas he learned from Gauss: the Dirichlet principle, magnetism, etc.; they are addressed in several other chapters of the present book.

The first contact between Gauss and Riemann took place probably in 1846, before Gauss became officially Riemann's mentor. In that year, in a letter to his father dated November 5 and translated in \cite{Riemann-Letters}, Riemann informs the latter about the courses he plans to follow,  and among them he mentions a course by Gauss on ``the theory of least squares."\footnote{The other courses are on the Cultural History of Greece and Rome, Theology, Recent Church history, General Physiology and Definite Integrals. Riemann had also the possibility to choose courses among Probability, Mineralogy and General Natural History. He adds: ``The most useful to me will be mineralogy. Unfortunately it conflicts with Gauss's lecture, since it is scheduled at 10 o'clock, and so I'd be able to attend only if Gauss moved his lecture forward, otherwise it looks like it won't be possible. General Natural History would be very interesting, and I would certainly attend, if along with everything else I had enough money."} During his two years stay in Berlin (1847--1849), Riemann\index{Riemann, Bernhard  (1826--1866)} continued to study thoroughly Gauss's\index{Gauss, Carl Friedrich (1777--1855)} papers. In another letter to his father, dated May 30, 1849, he writes (translation in \cite{Riemann-Letters}): 
\begin{quote}\small
Dirichlet\index{Dirichlet, Johann Peter Gustav Lejeune (1805--1859)} has arranged to me to have access to the library. Without his assistance, I fear there would have been obstacles. I am usually in the reading room by nine in the morning, to read two papers by Gauss that are not available anywhere else. I have looked fruitlessly for a long time in the catalog of the royal library for another work of Gauss,\index{Gauss, Carl Friedrich (1777--1855)} which won the Copenhagen prize, and finally just got it through Dr. Dale of the Observatory. I am still studying it.
\end{quote}

During the same stay in Berlin, Riemann\index{Riemann, Bernhard  (1826--1866)} followed lectures by  Dirichlet\index{Dirichlet, Johann Peter Gustav Lejeune (1805--1859)} on topics related to Gauss's works. He writes to his father (letter without date, quoted in \cite{Riemann-Letters}): 
 \begin{quote}\small
 My own course of specialization is the one with Dirichlet;\index{Dirichlet, Johann Peter Gustav Lejeune (1805--1859)} he lectures on an area of mathematics to which Gauss owes his entire reputation. I have applied myself very seriously to this subject, not without success, I hope.
 \end{quote}
 
Regarding his written production, Riemann endorsed Gauss' principle: \emph{pauca sed matura} (few but ripe).

 Riemann,\index{Riemann, Bernhard  (1826--1866)} as a child, liked history. In a letter to his father, dated May 3, 1840 (he was 14), he complains about the fact that at his Gymnasium there were fewer lessons on history than on \emph{Rechnen} (computing), cf. \cite{Riemann-Letters}. On August 5, 1841, he writes, again to his father, that he is the best student in history in his class. Besides history, Riemann was doing very well in Greek, Latin, and German composition (letters of February 1, 1845 and March 8, 1845). According to another letter to his father, dated April 30, 1845, it  is only in 1845 that Riemann started being really attracted by mathematics. In the same letter, Riemann declares that he plans to enroll the University of G\"ottingen to study theology, but that in reality he must decide for himself what to do, since otherwise he ``will bring nothing good to any subject." 
 
Besides Euler and Gauss, we shall mention several other mathematicians. Needless to say, it would have been unreasonable to try to be exhaustive in this paper; the subject would need a book, and even several books. We have tried  to present a few markers on the history of the major questions that were studied by Riemann, insisting only on the mathematicians whose works and ideas had an overwhelming impact on him.

 \medskip
 
   The content of the rest of this paper is the following.
 
Section \ref{s:complex} is essentially an excursion into the realm of Euler's\index{Euler, Leonhard (1707--1783)} ideas on the notion of function, with a stress on algebraic functions and functions of a complex variable. Algebraic functions are multivalued, and Euler included these functions\index{function!multi-valued}\index{multi-valued function} as an important element of the foundations of the field of analysis, which he laid down in his famous treatise \emph{Introductio in analysin infinitorum}\index{Euler!Introductio in analysin infinitorum}   (Introduction to the analysis of the infinite) \cite{Euler-Int-b}. Riemann's work on what became known as Riemann surfaces was largely motivated by the desire to find a domain of definition for an algebraic multi-valued function on which it becomes single-valued.\footnote{As a matter of fact, this is the origin of the use of the word ``uniformization" by Riemann.} The study of functions of a complex variable, which includes as a special case that of algebraic functions, is one of the far-reaching subjects of Riemann's investigations, and its development is one of the few most important  achievements of the nineteenth century (probably the most important one).

 Section  \ref{s:elliptic} is concerned with elliptic integrals. These integrals constitute a class of complex functions with new interesting properties, and the work described in this section is a natural sequel to that which is reviewed in \S \ref{s:complex}. We shall mention works done on this subject by Johann Bernoulli, Fagnano, Euler (who published thirty-three memoirs on elliptic integrals), 
Legendre, Abel and Jacobi.

 Section \ref{s:Abelian}  focusses on Abelian functions, a vast generalization of elliptic functions, which led to an important problem in which Riemann became interested, namely, the Jacobi inversion problem,\index{Jacobi inversion problem}\index{problem!Jacobi inversion} and  which he eventually solved using $\vartheta$ functions. In fact, Abelian integrals constitute one of the major topics that Riemann worked on. He started his investigation on this subject in his doctoral dissertation \cite{Riemann-Grundlagen} (1851), worked on it in his 1854 memoir \cite{Riemann-Abelian} whose title is quite rightly ``The theory of Abelian functions," and he never stopped working and lecturing on it during the few years that were left to him. Some lecture notes and memoirs by Riemann on Abelian functions were published posthumously. In particular, his memoir \emph{ \"Uber das Verschwinden der $\vartheta$-Functionen} (On the vanishing on theta functions) \cite{R2}, in which he gives a solution to Jacobi's problem of inversion for the general case of integrals of algebraic functions,  is analyzed in Chapter 4 of the present volume, written by Houzel \cite{Houzel}.

 Section \ref{s:hyper} is concerned with the so-called Gauss hypergeometric series. These series, in various forms, were studied by Euler in his \emph{Institutiones calculi integralis} (Foundations of integral calculus), a treatise in three volumes \cite{E342} (1768--1770), and in several other papers by him, and by Gauss. The hypergeometric series is a family of functions of the form 
 \[1+\frac{\alpha\beta}{1.\gamma}x
+\frac{\alpha(\alpha+1)\beta(\beta+1)}{1.2\gamma(\gamma+1)}x^2  
+\frac{\alpha(\alpha+1)(\alpha+2)\beta(\beta+1)(\beta+2)}{1.2.3\gamma(\gamma+1)(\gamma+2)}x^3+\ldots\] where
 $\alpha,\beta,\gamma$ are parameters and 
 where the variable is $x$. 
 
  Gauss considered that almost any transcendental function is obtained from a hypergeometric series by assigning special values to the parameters. By providing such a broad class of functions, the introduction of the hypergeometric series in the field of analysis opened up new paths.
 Besides Euler and Gauss, the predecessors of Riemann in this field include  Pfaff\index{Pfaff, Johann Friedrich (1765--1825)} and  Kummer\index{Kummer, Ernst (1810--1893)}.  
 
 In Section \ref{s:zeta}, we deal with the 
 zeta function. The history of this function is sometimes traced back to the work of Pietro Mengoli  (1625--1686) on the problem of finding the value of the infinite series of inverses of squares of  integers. Indeed, it is reasonable to assume that questions about this series were accompanied by questions about  the series of inverses of cubes and other powers.  But it was Euler again who studied $\displaystyle \sum_{1}^\infty \frac{1}{n^s}$ as a function of $s$ (for $s$ real), establishing the functional equation that it satisfies, and the relation with prime numbers. This was the starting point of Riemann's investigations on what became later known as the Riemann zeta function.
 
 In Section \ref{s:space}, we make a quick review of some works done by Riemann's predecessors on the notion of space. This is essentially a philosophical debate, but it has a direct impact on mathematics and in particular on Riemann's work on geometry, more especially on his habilitation dissertation.\index{Riemann! habilitation lecture}\index{habilitation lecture!Riemann} It is in his reflections on space that Riemann introduced in mathematics the notion of Mannigfaltigkeit, which he borrowed from the philosophical literature. This notion reflects Riemann's multi-faced view on space, and it is an ancestor of the modern notion of manifold. Our review of space is necessarily very sketchy, since this notion is one of the most fundamental notions of philosophy, and talking seriously about it would require a whole essay. In particular, there is a lot to say on the philosophy of space in the works of Newton, Euler and Riemann and the comparison between them, but it is not possible to do it in the scope of the present paper.  Our intent here is just to indicate some aspects of the notion of space as it appears in the works of these authors and those of some other philosophers, including Aristotle, Descartes and Kant, and, as much as possible in a short survey, to give some hints on the context in which they emerge. 
 
It is also important to say that the effect of this discussion on space goes far beyond the limits of philosophy. Euler's theories of physics are strongly permeated with his philosophical ideas on space. Gauss's differential geometry was motivated by his investigations on physical space, more precisely, on geodesy and astronomy, and, more generally, by his aspiration to understand the world around him. At a more philosophical level, Gauss was an enthusiastic reader of Kant, and he criticized the latter's views on space, showing that they do not agree with the recent discoveries -- his own and others' -- of geometry. Riemann, in this field, was an heir of Gauss. In his work, the curvature of space (geometric space) is the expression of the physical forces that act on it. These are some of the ideas that we try to convey in Section \ref{s:space} and in other sections of this paper.

 Section \ref{s:topo} is concerned with topology. Riemann is one of the main founders of this field  in the modern sense of the word, but several important topological notions may be traced back to Greek antiquity and to the later works of Descartes, Leibniz and Euler. We shall review the ideas of Leibniz, and consider in some detail the works of Descartes and Euler on the so-called Euler characteristic of a convex polyhedral surface, which in fact is nothing else but an invariant of the topological sphere, a question whose generalization is contained in Riemann's doctoral dissertation \cite{Riemann-Grundlagen} and his paper on Abelian functions \cite{Riemann-Abelian}, from where one can deduce the invariants of surfaces of arbitrary genus. 

 Section \ref{s:diff} is concerned with the differential geometry of surfaces. We review essentially the works of Euler, Gauss and Riemann, but there was also  a strong French school of differential geometry, operating between the times of Euler and Riemann, involving, among others, Monge and several of his students, and, closer to Riemann, Bonnet.\index{Bonnet, Pierre-Ossian (1819--1892)}    
    
  Section \ref{s:trigo} is a review of the history of trigonometric series and the long controversy on the notion of function that preceded this notion. In his Habilitation memoir, Riemann describes at length this important episode of eighteenth and nineteenth century mathematics which also led to his discovery of the theory of integration, which we discuss in the next section.
    
In  Section \ref{s:integration}, we review some of the history of the Riemann integral. From the beginning of integral calculus until the times of Legendre, passing through Euler, integration was considered as an antiderivative. Cauchy defined the integral by limits of sums that we call now Riemann sums, taking   smaller and smaller subdivisions of the interval of integration and showing convergence to make out of that a definition of the definite integral, but he considered only integrals of continuous functions, where convergence is always satisfied. It was Riemann who developed the first general theory of integration, leading to the notion of integrable and non-integrable function.

The concluding section, \S \ref{s:conclusion} contains a few remarks on the importance of returning to the texts of the old masters.
 \medskip

Some of the historical points in our presentation are described in more detail than others;  this reflects our personal taste and  intimate opinion on what is important in history and worth presenting in more detail in such a quick survey.    The reader will find at the end of this paper (before the bibliography) a table presenting in parallel some works of Euler and of Riemann on related matters.

 \section{Functions} \label{s:complex}

Vito Volterra,\index{Volterra, Vito (1860--1940)} in his 1900 Paris ICM plenary lecture  \cite{Volterra1900}, declared that the nineteenth century was ``the century of function theory."\footnote{The title of Volterra's lecture is: \emph{Betti,\index{Betti, Enrico (1923--1892)} Brioschi,\index{Brioschi, Francesco (1824--1897)} Casorati\index{Casorati, Felice (1835--1890)} :  Trois analystes italiens et trois mani\`eres d'envisager les questions d'analyse} (Betti, Brioschi, Casorati: Three Italian analysts and three manners of  addressing the analysis questions). In that lecture, Volterra presents three different ways of doing analysis, through the works of  Betti, Brioschi and Casorati, who are considered as the founders of modern Italian mathematics. The three mathematicians had very different personalities, and contrasting approaches to analysis, but in some sense they were complementing each other. In particular, Brioschi was capable of doing very long calculations, Betti was a geometer repugnant to calculations, and Casorati was an excellent teacher and an applied mathematician.}
In the language of that epoch, the expression ``function theory" refers, in the first place, to functions of a complex variable. One of the mottos, which was the result of a thorough experience in the domain, was that a function of a real variable acquires its full strength when it is complexified, that is, when it is extended to become a function of a complex variable. This idea was shared by Cauchy, Riemann Weierstrass, and others to whom we refer now.

On functions of a complex variable, we first quote a letter from Lagrange to Antonio Lorgna,\index{Lorgna, Antonio Maria (1735--1796)} an engineer and the governor of the military school at Verona who made important contributions to mathematics, physics and chemistry.
The letter is dated December 20, 1777. Lagrange writes (cf. Lagrange's \OE uvres, \cite{Lagrange-oeuvres}  t. 14, p. 261):
 \begin{quote}\small I consider it as one of the most important steps made by Analysis in the last period, that of not being bothered any more by imaginary quantities, and to be able to submit them to calculus, in the same way as the real ones.\footnote{Je regarde comme un des pas les plus importants que l'Analyse ait faits dans ces derniers temps, de n'\^etre plus embarrass\'ee des quantit\'es imaginaires et de pouvoir les soumettre au calcul comme les quantit\'es r\'eelles.}
\end{quote}

Gauss,\index{Gauss, Carl Friedrich (1777--1855)} who, among other titles he carried, was one of the main founders of the theory of functions of a complex variable, was also responsible for the introduction of complex numbers in several theories. In particular, he realized their power in number theory, and he used this in his \emph{Disquisitiones arithmeticae}\index{Gauss!Disquisitiones arithmeticae} (Arithmetical researches) \cite{Gauss-D} (1801), a masterpiece he wrote at the age of 24.  In his second paper on biquadratic residues \cite{Gauss-bi} (Sec. 30), he writes that ``number theory is revealed in its entire simplicity and natural beauty when the field of arithmetic is extended to the imaginary numbers." He explains that this means admitting integers of the form $a + bi$. ``Such numbers," he says, ``will be called complex integers."

In the same vein, Riemann, who had a marked philosophical viewpoint on things, writes, regarding complex functions, in \S 20 of his doctoral dissertation,\index{doctoral dissertation!Riemann}\index{Riemann!doctoral dissertation} \emph{Grundlagen f\"ur eine allgemeine Theorie der Functionen einer ver\"anderlichen complexen Gr\"osse} \cite{Riemann-Grundlagen}  (Foundations of a general theory of functions of a variable complex magnitude) \cite{Riemann-Grundlagen} (1851): ``Attributing complex values to the variable quantities reveals a harmony and a regularity which otherwise would remain hidden."  

Finally, let us quote someone closer to us, Jacques Hadamard, from his \emph{Psychology of invention in the mathematical field} \cite{Hadamard-essai}. A sentence by him which is often repeated is that ``the shortest and the best way between two truths of the real domain often passes by the imaginary one." We quote the whole passage (\cite{Hadamard-essai} p. 122--123):
\begin{quote}\small
It is Cardan, who is not only the inventor of a well-known joint which is an essential part of automobiles, but who has also fundamentally transformed mathematical science by the invention of imaginaries. Let us recall what an imaginary quantity is. The rules of algebra show that the square of any number, whether positive or negative, is a positive number: therefore, to speak of the square root of a negative number is mere absurdity. Now, Cardan deliberately commits that absurdity and begins to calculate on such ``imaginary" quantities.

One would describe this as pure madness; and yet the whole development of algebra and analysis would have been impossible without that fundament -- which, of course, was, in the nineteenth century, established on solid and rigorous bases.
It has been written that the shortest and the best way between two truths of the real domain often passes by the imaginary one.
\end{quote}

In the rest of this section, we review some markers in the history of functions, in particular functions of a complex variable and  algebraic functions, two topics which are at the heart of Riemann's work on Riemann surfaces, on Abelian functions, on the zeta function, and on other topics. Before that, we make a digression on the origin of the general notion of function.

It is usually considered that Euler's\index{Euler, Leonhard (1707--1783)} \emph{Introductio in analysin infinitorum}\index{Euler!Introductio in analysin infinitorum} \cite{Euler-Int-b} is the first treatise in which one can find the definition of a  function, according to modern standards, and where functions are studied in a systematic way. We take this opportunity to say a few words on Euler's treatise, to which we refer several times in the rest of this paper.

   The \emph{Introductio}\index{Euler!Introductio in analysin infinitorum}   is a treatise in two volumes, first published in 1748, which is concerned with a variety of subjects, including (in the first volume) algebraic curves, trigonometry, logarithms, exponentials and their definitions by limits, continued fractions, infinite products, infinite series and integrals. The second volume is essentially concerned with the differential geometry of curves and surfaces. The importance of the \emph{Introductio} lies above all in the fact that it made analysis the branch of mathematics where one studies functions. But the \emph{Introductio} is more than a treatise with a historical value. Two hundred and thirty years after the first edition appeared in print, Andr\'e Weil\index{Weil, Andr\'e (1906--1998)} considered that it was more useful for a student in mathematics to study that treatise rather than any other book on analysis. This is reported on by John Blanton who writes, in his English edition of the \emph{Introductio}\index{Euler!Introductio in analysin infinitorum} \cite{Euler-Intro}: 
   \begin{quote}\small
   In October, 1979, Professor Andr\'e Weil spoke at the University of Rochester on the life and work of Leonhard Euler. One of his remarks was to the effect that he was trying to convince the mathematical community that students of mathematics would profit much more from a study of Euler's \emph{Introductio in analysin infinitorum}, rather than the available modern textbooks.
   \end{quote}

   The importance of this work has also been highlighted by several other mathematicians. C. B. Boyer, in his 1950 ICM communication (Cambridge, Mass.) \cite{Boyer}, compares the impact of the \emph{Introductio} to that of Euclid's\index{Euclid!Elements}\index{Euclid!Elements} \emph{Elements} in geometry and to al-Khw\=arizm\=\i 's  \emph{Jabr} in algebra. He writes:

 \begin{quote}\small 
 The most influential mathematics textbooks of ancient times (or, for that matter, of all times) is easily named. The \emph{Elements} of Euclid,\index{Euclid!Elements} appearing in over a thousand editions, has set the pattern in elementary geometry ever since it was composed more than two and a quarter millenia ago. The medieval textbook which most strongly influenced mathematical development is not so easily selected;  but a good case can be made out of \emph{Al-jabr wal muq\=abala} of al-Khw\=arizm\=\i , just about half as old as the \emph{Elements}.\index{Euclid!Elements} From this Arabic work, algebra took its name and, to a great extent, its origin. Is it possible to indicate a modern textbook of comparable influence and prestige? Some would mention the \emph{G\'eom\'etrie} of Descartes,\index{Descartes, Ren\'e (1596--1650)}  or the \emph{Principia}\index{Newton!Principia} of Newton\index{Newton, Isaac (1643--1727)} or the \emph{Disquisitiones} of Gauss;\index{Gauss, Carl Friedrich (1777--1855)}  but in pedagogical significance these classics fell short of a work less known. [...] over these well known textbooks there towers another, a work which appeared in the very middle of the great textbook age and to which virtually all later writers admitted indebtedness. This was the \emph{Introductio in analysin infinitorum}\index{Euler!Introductio in analysin infinitorum} of Euler, published in two volumes in 1748. Here in effect Euler accomplished for analysis what Euclid and al-Khw\=arizm\=\i \  had done for synthetic geometry and elementary algebra, respectively.
 \end{quote}

Even though the \emph{Introductio} is generally given as the main reference for the introduction of functions in analysis, regarding the \emph{usage} of functions, one can go far back into history. Tables of functions exist since the Babylonians (some of their astronomical tables survive). Furthermore, in ancient  Greece, mathematicians manipulated functions, not only in the form of tables. In particular, the chord function (an ancestor of our sine function)\footnote{The relation between chord and sine is : $\sin x =\frac{1}{2}\mathrm{crd}\  2x$.} is used extensively in some Greek treatises. For instance, Proposition 67 of Menelaus'\index{Menelaus of Alexandria (c. 170--c. 140)} \emph{Spherics} (1st-2nd century A.D.) says the following \cite{Rashed-Menelaus}:
     \begin{quote}\small
  Let  $ABC$ and $DEG$ be two spherical triangles whose  angles $A$ and $D$ are equal, and where $C$ and $G$ are either equal or their sum is equal to two right angles. Then,
$$\displaystyle \frac{\mathrm{crd} \ 2AB}{\mathrm{crd} \ 2BC}=\frac{\mathrm{crd} \ 2DE}{\mathrm{crd} \ 2EG}.$$
\end{quote}

  Youschkevitch, in his interesting survey \cite{Y}, argues that the general idea of a dependence of a quantity upon another one is absent from Greek geometry. The author of the present paper declares that if in the above proposition of Menelaus one does not see the notion of function, and hence the general idea of a dependence of a quantity upon another one, then this author fails to know what mathematicians mean by the word function. 
 
 Leibniz\index{Leibniz, Gottfried Wilhelm (1646--1716)}  and Johann I Bernoulli,\index{Bernoulli, Johann (1667--1748)} who were closer to Euler,\index{Euler, Leonhard (1707--1783)} manipulated functions, even though the functions they considered were always associated with geometrical objects, generally, curves in the plane. For instance, in a memoir published in 1718 on the isoperimetry problem in the plane, \cite{Bernoulli1718} Bernoulli writes:
\begin{quote}\small
Here, we call \emph{function} of a variable magnitude, a quantity formed in whatever manner with that variable magnitude and constants.\footnote{On appelle ici \emph{Fonction} d'une grandeur variable, une quantit\'e compos\'ee de quelque mani\`ere que ce soit avec cette grandeur variable et des constantes. [The emphasis is Bernoulli's]}
\end{quote}
  The functions that Bernoulli considers in this memoir are associated to \emph{arbitrary} curves in the plane having the same perimeter, among which Bernoulli looks for the one which bounds the greatest area. This is an example of the general idea  that before Euler, analysis was tightly linked to geometry, and the study of functions consisted essentially in the study of curves associated to some geometric properties.  With the  \emph{Introductio}, things became different. Analysis started to release itself from geometry, and functions were studied for themselves. Let us now make a quick review of the part of this treatise which concerns us here.

The first paper is called \emph{On functions in general}. In this chapter, Euler\index{Euler, Leonhard (1707--1783)} states his general definition of a function, after a description of what is a variable quantity:
\begin{quote}\small
 A function of a variable quantity is an analytic\index{analytic function (in the sense of Euler)}\index{function!analytic (in the sense of Euler)} expression composed in any way whatsoever of the variable quantity and numbers or constant quantities.
 \end{quote}
  The word ``analytic" means in this context that the  function is obtained by some process that uses the four operations (addition, subtraction, multiplication and division), together with root extraction, exponentials, logarithms, trigonometric functions, derivatives and integrals. Analyticity in terms of being defined by a convergent power series is not intended by this definition. The meaning of the word ``analytic function" rather is ``a function used in (the field of) analysis."\index{analytic function (in the sense of Euler)}\index{function!analytic (in the sense of Euler)} Concerning the notion of variable, Euler writes (\S 3):\footnote{We are using the translation from Latin in \cite{Euler-Int-b}.} 
\begin{quote}\small [...] Even zero and complex numbers are not excluded from the signification of a variable quantity.
\end{quote}
Thus, functions of a complex variable are included in Euler's \emph{Introductio}.  We note however that in this treatise, Euler, in his examples, always deals with functions that are given by formulae: polynomials, exponentials, logarithms, trigonometric functions, etc. but also infinite products and infinite sums.

After the definition of a function, we find in the  \emph{Introductio}  the definition of an algebraic function.
In \S 7, Euler writes: 
\begin{quote}
 \small Functions are divided into algebraic and transcendental.\index{transcendental function}\index{function!transcendental} The former are those made up from only algebraic operations, the latter are those which involve transcendental operations. 
\end{quote}

And in \S 8: 
\begin{quote}\small Algebraic functions\index{algebraic function}\index{function!algebraic} are subdivided into non-irrational and irrational functions: the former are such that the variable quantity is in no way involved with irrationality; the latter are those in which the variable quantity is affected by radical signs. 
\end{quote}

Concerning irrational functions (\S\,9), he writes:
\begin{quote}\small
It is convenient to distinguish these into explicit and implicit irrational functions.
 
The explicit functions are those expressed with radical signs, as in the given examples. The implicit are those irrational functions which arise from the solution of equations. Thus $Z$ is an implicit irrational function of $z$ if it is defined by an equation such as $Z^7=az$ or $Z^2=bz^5$. Indeed, an explicit value of $Z$ may not be expressed even with radical signs, since common algebra has not yet developed to such a degree of perfection.
 
 \end{quote}

And in \S 10:
\begin{quote}
  \small Finally, we must make a distinction between single-valued and multi-valued functions.\index{multi-valued function}\index{function!multi-valued}

A single-valued function is one for which, no matter what value is assigned to the variable $z$, a single value of the function is determined. On the other hand, a multi-valued function is one such that, for some value substituted for the variable $z$, the function determines several values. Hence, all non-irrational functions, whether polynomial or rational, are single-valued functions, since expressions of this kind, whatever value be given to the variable $z$, produce a single value. However, irrational functions are all multi-valued, because the radical signs are ambiguous and give paired values. There are also among the transcendental\index{transcendental function}\index{function!transcendental} functions, both single-valued and multi-valued functions; indeed, there are infinite-valued functions. Among these are the arcsine of $z$, since there are infinitely many circular arcs with the same sine. 
\end{quote}

   Euler\index{Euler, Leonhard (1707--1783)} then gives examples of two-valued, three-valued and four-valued functions, and in \S 14 he writes:
\begin{quote}
 \small Thus $Z$ is a multi-valued function\index{function!multi-valued}\index{multi-valued function} of $z$ which for each value of $z$, exhibits $n$ values of $Z$ where $n$ is a positive integer. If $Z$ is defined by this equation
\[Z^n-PZ^n-1+QZ^{n-2}-RZ^{n-3}+SZ^{n-4}-\ldots=0\]
[...] Further it should be kept in mind that the letters $P,Q,R,S$, etc. should denote single-valued functions of $z$. If any of them is already a multi-valued function, then the function $Z$ will have many more values, corresponding to each value of $z$, than the exponent would indicate. It is always true that if some of the values are complex, then there will be an even number of them. From this we know that if $n$ is an odd number, there will be at least one real value of $Z$.  \end{quote}
 
 He then makes the following remarks:

\begin{quote}\small
If $Z$ is a multi-valued function\index{function!multi-valued}\index{multi-valued function} of $z$ such that it always exhibits a single real value, then $Z$ imitates a single-valued function of $z$, and frequently can take the place of  a single-valued function.
 
Functions of this kind are $P^{\frac{1}{3}}$, $P^{\frac{1}{5}}$, $P^{\frac{1}{7}}$, etc. which indeed give only one real value, the others all being complex, provided $P$ is a single-valued function of $z$. For this reason, an expression of the form $P^{\frac{m}{n}}$, whenever $n$ is odd, can be counted as a single-valued function, whether $m$ is odd or even. However, if $n$ is even then $P^{\frac{m}{n}}$ will have either no real value or two; for this reason, expressions of the form $P^{\frac{m}{n}}$, with $n$ even, can be considered to be two-valued functions, provided the fraction $\frac{m}{n}$ cannot be reduced to lower terms.
\end{quote}

  From this discussion we single out the fact that algebraic functions\index{algebraic function}\index{function!algebraic} are considered as functions, even though they are multi-valued. They are solutions of algebraic equations. Since we are talking about history, it is good to recall that the study of such equations\index{algebraic equation}\index{equation!algebraic} is an old subject that can be traced back to the work done on algebraic curves by the Greeks. In fact, Diophantus (3d century B.C.)\index{Diophantus of Alexandria (3d century B.C.)} thoroughly studied integral solutions of what is now called ``Diophantine equations."  They are examples of algebraic equations.\footnote{For what concerns Diophantus' \emph{Arithmetica}, 
 we refer the interested reader to the recent and definitive editions \cite{Diophante1}, \cite{Diophante2}, \cite{Diophante3}, and \cite{Rashed-Diophante} by R. Rashed.} 
 Algebraic equations are also present in the background of the geometric work of Apollonius (3d--2d century B.C.)\index{Apollonius of Perga (3d--2d century B.C.)} on conics. In that work, intersections of conics were used to find geometrical solutions of algebraic equations.\footnote{For a recent and definitive edition of Apollonius' \emph{Conics}, we refer the reader to the volumes \cite{Apollonius1}, \cite{Apollonius2}, \cite{Apollonius3}, \cite{Apollonius4}, and \cite{Apollonius5}, again edited by R. Rashed.}
 It is true however that in these works, there is no \emph{definition} of an algebraic function as we intend it today, and in fact at that time there was no definition of function at all.

     The multi-valuedness  of algebraic functions gave rise to tremendous developments by Cauchy \index{Cauchy, Augustin-Louis (1789--1857)} and Puiseux,\index{Puiseux, Victor-Alexandre (1820--1883)} and it was also a major theme in Riemann's\index{Riemann, Bernhard  (1826--1866)} work, in particular in his doctoral dissertation \cite{Riemann-Grundlagen} (1851) and his memoir on Abelian functions \cite{Riemann-Abelian} (1857). In fact, the main reason for which Riemann introduced the surfaces that we call today Riemann surfaces was to find ground spaces on which  multi-valued functions are defined and become single-valued. We discuss the works of Cauchy and Puiseux in relation with that of Riemann in Chapter 7 of the present volume, \cite{Papa-Puiseux}.

We note for later use that a definition of ``continuity" is given in Volume 2 of the \emph{Introductio}, where Euler says that a curve is continuous if it represents ``one determinate function," and discontinuous if it is decomposed into ``portions that represent different continuous functions." We shall see that such a notion was criticized by Cauchy (regardless of the fact that it is called ``continuity").\footnote{There are other imperfections in the \emph{Introductio}, even though this book is one of the most interesting treatises ever written on elementary analysis.}
    
We note finally that it is usually considered that the expression \emph{analysin infinitorum}\index{Euler!Introductio in analysin infinitorum} in the title of Euler's treatise does not refer to the field of infinitesimal analysis in the sense of Newton\index{Newton, Isaac (1643--1727)} or Leibniz,\index{Leibniz, Gottfried Wilhelm (1646--1716)} but, rather, to the use of infinity\index{infinity} (infinite series, infinite products, continued fractions expansions, integral representations, etc.) in analysis. Euler\index{Euler, Leonhard (1707--1783)} was also the first to highlight the zeta function,\index{zeta function}\index{Riemann!zeta function} the gamma function and elliptic integrals\index{elliptic integral} as functions.  However, it is good to recall that infinite sums were known long before Euler.\index{Euler, Leonhard (1707--1783)} For instance, Zeno of Elea (5th c. B.C.)\index{Zeno of Elea (5th c. B.C.)} had already addressed the question of convergence of infinite series, and to him are attributed several well-known paradoxes in which the role and the significance of infinite series and their convergence are emphasized (the paradox of Achilles and the tortoise, the arrow paradox,  the paradox of the grain of millet, etc.). However, infinite series are not considered as functions in these works. Zeno's paradoxes are commented in detail in Aristotle's\index{Aristotle (384--322 B.C.)}  \emph{Physics} \cite{Aristotle-Physics}, but also by mathematicians and philosophers from the modern period, including Bertrand Russell,\index{Russell, Bertrand (1872--1970)} Hermann Weyl\index{Weyl, Hermann (1885--1955)}, Paul Tannery and several others; cf.  \cite{Russell} p. 346--354, and \cite{Weyl-philo} and \cite{P.Tannery}.

We also recall that convergent series were used by Archimedes\index{Archimedes (c. 287 B.C.--c. 212 B.C.)} in his computations of areas and volumes. 

  Before leaving this book, let us mention that Euler\index{Euler, Leonhard (1707--1783)} establishes there a hierarchy among transcendental functions\index{transcendental function}\index{function!transcendental} by introducing a notion close to what we call today the transcendence degree of a function.

 In his later works, Euler\index{Euler, Leonhard (1707--1783)} dealt with much more general functions. For instance, in his 1755 memoir \cite{E213}, entitled \emph{Remarques sur les m\'emoires pr\'ec\'edents de M. Bernoulli} (Remarks on the preceding memoirs by Mr. Bernoulli),  any mechanical curve (that is any curve drawn by hand) is associated with a function.\footnote{One may recall here that the mathematicians of Greek antiquity (Archytas of Tarentum,\index{Archytas of Tarentum (428?347 B.C.)} Hippias, Archimedes, etc.) who examined curves formulated a mechanical definition. The curves with which they dealt were not necessarily defined by equations, they were ``traced by a moving point," sometimes (in theory) using a specific mechanical device. Of some interest here would be the connections between this subject and the theory of mechanical linkages,\index{mechanical linkage} which was extensively developed in the nineteenth century and became fashionable again in the twentieth century. A conjecture by Thurston says (roughly speaking) that any ``topological curve" is drawable by a mechanical linkage. This is a vast generalization of a result of Kempe stating that any bounded piece of an algebraic curve is drawable by some linkage, cf.  \cite{Kempe}.  We refer the reader to
Sossinsky's survey of this subject and its recent developments \cite{Sossinsky}, in particular the solution of Thurston's conjecture.} 
  In his \emph{Institutiones calculi differentialis cum eius usu in analysi finitorum ac doctrina serierum}
(Foundations of differential calculus, with applications to finite analysis and series)  \cite{E212}, also published in 1755, Euler\index{Euler, Leonhard (1707--1783)} gave again a very general definition of a function (p. vi):
\begin{quote}\small
Those quantities that depend on others in this way, namely, those that undergo a change when others change, are called functions of these quantities. This definition applies rather widely and includes all ways in which one quantity could be determined by another.
\end{quote}
Likewise, in his memoir \cite{Euler-rep-1777}, \emph{De repraesentatione superficiei sphaericae super plano} (On the representation of Spherical Surfaces onto the
Plane) (1777), Euler\index{Euler, Leonhard (1707--1783)} dealt with ``arbitrary mappings" between the sphere and the plane. He writes:\footnote{We are using George Heines' translation.}
\begin{quote}\small
 I take the word ``mapping" in the widest possible sense;
any point of the spherical surface is represented on the plane by any desired
rule, so that every point of the sphere corresponds to a specified point in the
plane, and inversely.
\end{quote}

We shall consider again the question of functions, from the epoch of Euler and until the work of Riemann, in \S \ref{s:trigo} concerned with trigonometric functions.

 Riemann,\index{Riemann, Bernhard  (1826--1866)} in his doctoral dissertation,\index{doctoral dissertation!Riemann}\index{Riemann!doctoral dissertation} \cite{Riemann-Grundlagen} (1851), also considers arbitrary functions. In fact, the dissertation starts as follows: ``If we designate by $z$ a variable magnitude, which may take successively all possible real values, then, when to each of these values corresponds a unique value of the indeterminate magnitude $w$, we say that $w$ is a function of $z$ [...]"
One may also refer to the beginning of \S XIX of the same dissertation, where Riemann states that the principles he is presenting are the bases of a general theory of functions which is independent of any explicit expression.

 The details of the seventeenth-century debate concerning functions are rather confusing if one does not include them in their historical context. For instance, the notion of ``continuity" which we alluded to and which  is referred to in the debate is different from what we intend today by this word. In fact, the word ``continuity," even restricted to the works of Euler,\index{Euler, Leonhard (1707--1783)} varied in the course of time. 
 
 Cauchy,\index{Cauchy, Augustin-Louis (1789--1857)}  the major figure standing between  Euler\index{Euler, Leonhard (1707--1783)} and Riemann for what concerns the notion of  function, in his \emph{M\'emoire sur les fonctions continues} (Memoir on continuous functions) \cite{Cauchy1844}, starts as follows: 
 \begin{quote}\small 
 In the writings of Euler\index{Euler, Leonhard (1707--1783)} and Lagrange,\index{Lagrange, Joseph-Louis (1736--1813)}  a function is termed \emph{continuous}\index{continuous  function (in the sense of Euler)}\index{function!continuous (in the sense of Euler)} or \emph{discontinuous}\index{discontinuous!function (in the sense of Euler)}\index{function!discontinuous (in the sense of Euler)} according to whether the various values of this function corresponding to various values of the variable follow or not the same law, or are given or not by only one equation. It is in these terms that the continuity of functions was defined by these famous geometers, when they used to say that ``the arbitrary functions, introduced by the integration of partial differential equations, may be continuous or discontinuous functions." However, the definition which we just recalled is far from offering mathematical accuracy [...] A simple change in notation will often suffice to transform a continuous fonction into a discontinuous one, and conversely.\footnote{Dans les ouvrages d'Euler et de Lagrange, une fonction est appel\'ee \emph{continue} ou \emph{discontinue}, suivant que les diverses valeurs de cette fonction, correspondantes \`a diverses valeurs de la variable, sont ou ne sont pas assujetties \`a une m\^eme loi, sont ou ne sont pas fournies par une seule \'equation. C'est en ces termes que la continuit\'e des  fonctions se trouvait d\'efinie par ces illustres g\'eom\`etres, lorsqu'ils disaient que ``les fonctions arbitraires, introduites par l'int\'egration des \'equations aux d\'eriv\'ees partielles, peuvent \^etre des fonctions continues ou discontinues." Toutefois, la d\'efinition que nous venons de rappeler est loin d'offrir une pr\'ecision math\'ematique [...] un simple changement de notation suffira souvent pour transformer une fonction continue en fonction discontinue, et r\'eciproquement.}
 \end{quote}
 
In fact, one might consider that Euler's\index{Euler, Leonhard (1707--1783)} definition of continuity is  just one definition that is different from the new definition which Cauchy had in mind (and which is the definition we use today). This would have been fine, and it would not be the only instance in mathematics where the same word is used for notions that are different, especially at different epochs. But Cauchy showed by an example that in this particular case Euler's definition is inconsistent, because the property it expresses depends on the parametrization that is used. Cauchy\index{Cauchy, Augustin-Louis (1789--1857)}  continues:
  \begin{quote}\small
  But the non-determinacy will cease if we substitute to Euler's definition the one I gave in Chapter II of the \emph{Analyse alg\'ebrique}. According to the new definition, a function of the variable $x$ will be \emph{continuous} between two limits $a$ and $b$ of this variable if between two limits the function has always a value which is unique and finite, in such a way that an infinitely small increment of this variable always produces an infinitely small increment of  the function itself.\footnote{Mais l'ind\'etermination cessera si \`a la d\'efinition d'Euler on substitue celle que j'ai donn\'ee dans le chapitre II de l'\emph{Analyse alg\'ebrique}. Suivant la nouvelle d\'efinition, une fonction de la variable r\'eelle $x$ sera continue entre deux limites $a$ et $b$ de cette variable, si, entre ces limites, la fonction acquiert constamment une valeur unique et finie, de telle sorte qu'un accroissement infiniment petit de la variable produise toujours un accroissement infiniment petit de la fonction elle-m\^eme.}
  \end{quote}
            
We quoted these texts in order to give an idea of the progress of the notion of continuity. We now come to the study of functions of a complex variable.

         In his memoir on Abelian functions, Riemann refers explicitly to Gauss for the fact that we represent a complex magnitude $z=x+iy$ by a point in the plane with coordinates $x$ and $y$.

It is not easy to know when the theory of functions of a complex variable started, and, in fact, the answer depends on whether one studies holomorphic functions, and what properties of holomorphic functions are meant (before the epoch of Riemann, they were not known to be equivalent): angle-preservation, power series expansion, the Cauchy-Riemann equation, etc. 

 Euler used complex variables and the notion of conformality (angle-preservation) in his memoirs on geographical maps. He wrote three memoirs on this subject, \emph{De repraesentatione superficiei sphaericae super plano} (On the representation of  spherical surfaces on a plane)
 \cite{Euler-rep-1777}, 
\emph{De proiectione geographica superficiei sphaericae} (On the geographical projections of spherical surfaces) \cite{Euler-pro-1777}, and 
 \emph{De proiectione geographica Deslisliana in mappa generali imperii russici usitata} (On Delisle's geographic projection used in the general map of the Russian empire)   \cite{Euler-pro-Desli-1777}.
 The three memoirs were published in 1777.
 In the development of the theory, he used complex numbers to represent angle-preserving maps. Lagrange also studied angle-preserving maps, in his memoir \emph{Sur la construction des cartes g\'eographiques} (On the construction of geographical maps)
  \cite{Lagrange1779}, published in 1779.
 
 In fact, the notion of angle-preserving map can be traced back to Greek antiquity, see the survey \cite{Papa-qc}. We already recalled that Euler, in his didactical treatise  \emph{Introductio in analysin infinitorum},\index{Euler!Introductio in analysin infinitorum} refers explicitly to functions in which the variable is a complex number.
 De Moivre,\index{Moivre@de Moivre, Abraham (1667--1754)} already in 1730,  considered polynomials defined on the complex plane, and it is conceivable that other mathematicians before him did the same \cite{Moivre}. Remmert, who, besides being a specialist of complex analysis, is a highly respected historian in this field, writes in his \emph{Theory of complex variables} \cite{Remmert-Theory} that the theory was born at the moment when Gauss\index{Gauss, Carl Friedrich (1777--1855)} sent a letter to Bessel, dated December 18, 1811, in which he writes:\footnote{The translation is Remmert's; cf. \cite{Remmert-Theory} p. 1.}
\begin{quote}\small
 At the beginning I would ask anyone who wants to introduce a new function in analysis to clarify whether he intends to confine it to real magnitudes (real values of the argument) and regard the imaginary values as just vestigial -- or whether he subscribes to my fundamental proposition that in the realm of magnitudes the imaginary ones $a+b\sqrt{-1}=a+bi$ have to be regarded as enjoying equal rights with the real ones. We are not talking about practical utility here; rather analysis is, to my mind, a self-sufficient science. It would lose immeasurably in beauty and symmetry from the rejection of any fictive magnitudes. At each stage truths, which otherwise are quite generally valid, would have to be encumbered with all sorts of qualifications.
 \end{quote}
  In fact, the letter also shows that at that time Gauss\index{Gauss, Carl Friedrich (1777--1855)} was already aware of the concept of complex integration, including Cauchy's integral theorem;  cf. \cite{Gauss} Vol. 8, p. 90--92.

  Cauchy, in his \emph{Cours d'analyse} \cite{Cauchy-cours} (1821), starts by defining  functions of real variables (p. 19), and then passes to complex variables. There are two distinct definitions in the real case, for functions of one or several variables:
\begin{quote}\small
When variable quantities are so tied to each other that, given the value of one of them, we can deduce the values of all the others, we usually conceive these various quantities expressed in terms of one of them, which then bears the name \emph{independent variable}; and the other quantities expressed in terms of the independent variable are what we call functions of that variable.

When variable quantities are so tied to each other that, given the values of some of them, we can deduce the values of all the others, we usually conceive these various quantities expressed in terms of several of them, which then bear the name \emph{independent variables}; and the remaining quantities expressed in terms of the independent variables, are what we call functions of these same variables.\footnote{Lorsque des quantit\'es variables sont tellement li\'ees entre elles que, la valeur de l'une d'elles \'etant donn\'ee, on puisse en conclure les valeurs de toutes les autres, on con\c coit d'ordinaire ces diverses quantit\'es exprim\'ees au moyen de l'une d'entre elles, qui prend alors le nom de \emph{variable ind\'ependante} ; et les autres quantit\'es exprim\'ees au moyen de la variable ind\'ependante sont ce qu'on appelle des \emph{fonctions} de cette variable.

Lorsque les quantit\'es variables sont tellement li\'ees entre elles que, les valeurs de quelques unes \'etant donn\'ees, on puisse en conclure celles de toutes les autres, on con\c coit ces diverses quantit\'es exprim\'ees au moyen de plusieurs d'entre elles, qui prennent alors le nom de \emph{variables ind\'ependantes} ; et les quantit\'es restantes, exprim\'ees au moyen des variables ind\'ependantes, sont ce qu'on appelle des \emph{fonctions} de ces m\^emes variables.}
\end{quote}

Talking about Cauchy's work on functions of a complex variable, one should also mention the Cauchy--Riemann equation\index{Cauchy--Riemann equations}\index{equations!Cauchy--Riemann} as a characterization of complex analycity, which Cauchy and Riemann introduced the same year, 1851, Cauchy in his papers  \cite{Cauchy-Sur-1851} and \cite{Cauchy-1851-mono} and Riemann in his doctoral dissertation \cite{Riemann-Grundlagen}. It is important to note also that the Cauchy--Riemann equations,\[\frac{\partial u}{\partial x}=\frac{\partial v}{\partial y}
 \ \ 
\hbox{and} \ \ 
 \frac{\partial u}{\partial y}=- \frac{\partial v}{\partial x}
,\] 
without the complex character, were used by d'Alembert\index{Alembert@d'Alembert, Jean le Rond (1717--1783)} in 1752, in his works on fluid dynamics, \emph{Essai d'une nouvelle th\'eorie de la r\'esistance des fluides} (Essay on a new theory of fluid resistance) \cite{YT} p. 27. D'Alembert showed later that functions $u$ and $v$ satisfying this pair of equations also satisfy Laplace's equation: $\Delta u=0$ and $\Delta v=0$.

The work of Cauchy is also reviewed in  the chapter \cite{Papa-Puiseux} in the present volume, written by the present author.

Riemann's\index{Riemann, Bernhard  (1826--1866)} doctoral dissertation \cite{Riemann-Grundlagen}\index{doctoral dissertation!Riemann}\index{Riemann!doctoral dissertation} is in some sense an essay on functions of a complex variable. Right at the beginning of the dissertation, Riemann states explicitly what he means by a function.  He starts with functions of a real variable: 
\begin{quote}\small If we designate by $z$ a variable magnitude, which may take successively all possible real values, then, if to each of these values corresponds a unique value of the indeterminate magnitude $w$, we say that $w$ is a function of $z$.
\end{quote}
He then talks about continuity of functions, in the modern sense of the word (as opposed to the sense that Euler\index{Euler, Leonhard (1707--1783)} gave to this word).\footnote{In the \emph{Introductio} Euler used the expression \emph{continuous function} for a function that is ``given by a formula." This is thoroughly discussed in \S \ref{s:trigo} of the present paper.} Then he writes: 
\begin{quote}\small 
This definition does not stipulate any law between the isolated values of the function, this is evident, because after this function has been dealt with for a given interval, the way it is extended outside this interval remains quite arbitrary. 
\end{quote}

Riemann\index{Riemann, Bernhard  (1826--1866)} then recalls that the possibility of using some ``mathematical law" that assigns to $w$ a value for a given value of $z$ was proper to the functions which Euler\index{Euler, Leonhard (1707--1783)} termed \emph{functiones continu\ae}. He writes that ``modern research has shown that there exist analytic expressions by which any continuous function on a given interval can be represented." He then declares that the case of functions of a complex variable is treated differently. In fact, Riemann considers only functions of a complex variable whose derivative does not depend on the direction, that is, holomorphic functions. He makes this property part of his definition of a function of a complex variable. Thus, when he talks about a function in the complex setting, he considers only conformal maps.

Regarding Riemann's dissertation, let us note that in a letter to his brother, dated November 26, 1851 \cite{Riemann-Letters}, after he submitted his doctoral dissertation manuscript, he writes that Gauss\index{Gauss, Carl Friedrich (1777--1855)} took it home to examine it for a few days, and that before reading it, Gauss told him: 
\begin{quote}\small
[Riemann speaking]  for years he had been preparing an essay, on which today he is still occupied, whose subject is the same or at least in part the same as that covered by me. Already in his doctoral dissertation now 52 years ago he actually expressed the intention to write on this subject. 
\end{quote}

This is an instance where Gauss was aware of a theory, or part of it, long before its author; we shall mention several other such instances in what follows.

\section{Elliptic integrals}\label{s:elliptic}

In the huge class of integrals of functions, the integrals of algebraic functions constitute the simplest and the most natural class to work with. The class of elliptic integrals (and their Abelian generalizations) which deal with such functions soon turned out to be enough tractable and at the same time very rich from the point of view of the problems that they posed. These integrals led to a huge amount of work by several prominent mathematicians, as we shall see  in this section.

Riemann had several reasons to work on Abelian integrals. Motivated by  lectures by Dirichlet,\index{Dirichlet, Johann Peter Gustav Lejeune (1805--1859)} Jacobi and others, he worked on the open problems that these functions presented, in particular the Jacobi inversion problem.\index{Jacobi inversion problem}\index{problem!Jacobi inversion}

When Riemann started his work on integrals as functions of a complex variable, this subject was already well developed. An important challenging problem that he tackled was the so-called Jacobi inversion problem\index{Jacobi inversion problem}\index{problem!Jacobi inversion} which we mention below. Most of all, these functions constituted for Riemann an interesting class of non-necessarily algebraic functions of a complex variable. The double periodicity of these integrals, the multi-valuedness of their inverses, the operations that one can perform on them, constituted a treasure of examples of new functions of a complex variable, and a context in which his theory of Riemann surfaces may naturally be used. 

We start by summarizing some of the main ideas and problems that concern elliptic functions that were addressed since the time of Euler.
\begin{enumerate}

\item The study of definite integrals representing arcs of conics and of lemniscates,  and the comparison of their properties with those of integrals representing arcs of circles, which are computable in terms of the trigonometric functions or their inverses. We recall, by way of comparison, that whereas the integral $\displaystyle\int_0^x\frac{dt}{\sqrt{1-t^2}}$ represents arc length along a circle centered at the origin, the integral $\displaystyle\int_0^x\frac{dt}{\sqrt{1-t^4}}$ represents arc length along the lemniscate of polar equation $r^2=\cos 2\theta$.

\item The search for sums and product formulae for such integrals, in the same way as there are formulae for sums and products of trigonometric functions.
 
\item The study of periods, again, in analogy with those of trigonometric functions.

\end{enumerate} 
      
        In fact, some of the first questions concerning elliptic integrals can be traced back to Johann I Bernoulli\index{Bernoulli, Johann (1667--1748)}  who tried to use the newly discovered integral calculus to obtain formulae for lengths of arcs of conic sections and some other curves. Bernoulli found the first addition formulae for such integrals. Finding general addition theorems for elliptic integrals remained one of the major problems for the following hundred years, involving the works of several major figures including Euler, Legendre, Abel, Jacobi and Riemann.  Bernoulli also discovered that  the lengths of some curves, expressed using integrals, may be expressed using infinite series \cite{Bernoulli1698}.                 

Johann Bernoulli was Euler\index{Euler, Leonhard (1707--1783)}'s teacher, and it is not surprising that the latter became interested in these problems early in his career. In his first paper on the subject,  \emph{Specimen de constructione aequationum differentialium sine indeterminatarum separatione} (Example of the construction of differential equations without separation of variables)  \cite{E28} written in 1733, Euler\index{Euler, Leonhard (1707--1783)} gives a formula for arc lengths of ellipses. He obtains them by first writing a differential equation satisfied by these arcs.  Generally speaking, Euler systematically searched for differential equations that describe the various situations that he was studying.

     Between the work of Bernoulli and that of Euler, we must mention that of Fagnano\index{Fagnano, Giulio Carlo de' Toschi di (1682--1766)}, who, around the year 1716, in a study he was carrying on the lemniscate, discovered some results which Euler considered several years later as outstanding. These results included an addition formulae for a class of elliptic integrals \cite{Fagnano-m}, and the fact that on an ellipse or a hyperbola, one may find infinitely many pairs of arcs whose difference is expressible by algebraic means. The word used by Euler and others for such arcs (or differences of arcs) is that they are ``rectifiable."  Fagnano managed to reduce the rectifiability of the lemniscate to that of the ellipse and hyperbola. A few words on Fagnano are in order.

   Giulio Carlo de' Toschi di  Fagnano (1682--1766) was a noble Italian interested in science, who worked during several decades in isolation, away from any scientific environment. Weil's\index{Weil, Andr\'e (1906--1998)} authoritative book on the history of number theory \cite{Weil-Hammurapi} starts with the following:
      \begin{quote}\small
      Accroding to Jacobi, the theory of elliptic functions was born between the twenty-third of December 1751 and the twenty-seventh of January 1752. On the former date, the Berlin Academy of Sciences handed over to Euler the two volumes of Marchese Fagnano's \emph{Produzioni Mathematiche}, published in Pesaro in 1750 and just received from the author; Euler was requested to examine the book and draft a suitable letter of thanks. On the latter date, Euler, referring explicitly to Fagnano's work on the lemniscate,  read to the Academy the first of a series of papers, eventually proving in full generality the addition and multiplication theorems for elliptic integrals.
      \end{quote}
     On p. 245 of the same treatise, Weil\index{Weil, Andr\'e (1906--1998)} writes:
              \begin{quote}\small
               On 23 December 1751 the two volumes of Fagnano's \emph{produzioni Mathematiche}, just published, reached the Berlin Academy and were handed over to Euler; the second volume contained reprints of pieces on elliptic integrals which appeared between 1714 and 1720 in an obscure Italian journal and had remained totally unknown. On reading these few pages Euler\index{Euler, Leonhard (1707--1783)} caught fire instantly; on 27 January 1752 he was presenting to the Academy a memoir \cite{E252} with an exposition of Fagnano's main results, to which he was already adding some of his own.
               
               The most striking of Fagnano's results concerned transformations of the ``lemniscate  differential"
               \[w(z)=\frac{dz}{\sqrt{1-z^4}};\]
               how he had reached them was more than even Euler could guess. ``Surely his discoveries would shed much light on the theory of transcendental functions," Euler wrote in 1753, ``if only his procedure supplied a sure method for pursuing these investigations further; but it rests upon substitutions of a tentative character, almost haphazardly applied ..."\footnote{The reference is to Euler's memoir   
         \emph{Specimen novae methodi curvarum quadraturas et rectificationes aliasque quantitates transcendentes inter se comparandi}
(An example of a new method for the quadurature and rectificaiton of curves and of comparing other quantities which are transcendentally related to each other)   \cite{E263}.}
              \end{quote}

              In a letter dated October 17, 1730 (\cite{Euler-Goldbach} p. 624), well before being aware of Fagnano's work, Euler informed Goldbach that ``even admitting logarithms,"\footnote{The reference to logarithms comes from the fact that $\frac{dt}{t}$ and some more general rational functions can be integrated using logarithms.} he could by no means compute the integral $\displaystyle \int \frac{a^2dx}{\sqrt{a^4-x^4}}$, that  ``expresses the curve element of the rectangular elastic curve, or rectify this ellipse." Fagnano, instead of giving explicit values, established equalities between such integrals which paved the way to a new series of results by Euler\index{Euler, Leonhard (1707--1783)} and others. In a letter to Goldbach, dated May 30, 1752, that is, about six months after reading Fagnano's work, Euler writes  (see \cite{Euler-Goldbach} p. 1064): ``Recently some curious integrations occurred to me." He first notes that three differential equations
          \[\frac{dx}{\sqrt{1-x^2}}=\frac{dy}{\sqrt{1-y^2}},\]
            \[\frac{dx}{\sqrt{1-x^4}}=\frac{dy}{\sqrt{1-y^4}},\]
            and   \[\frac{dx}{\sqrt{1-x^3}}=\frac{dy}{\sqrt{1-y^3}}\]
            can be integrated explicitly, and lead respectively to
            \[y^2+x^2=c^2+2xy\sqrt{1-c^2},\]
            \[y^2+x^2=c^2+2xy\sqrt{1-c^4}-c^2y^2\]
            and
              \[y^2+x^2+c^2x^2y^2=4c-4c^2(x+y)+2xy-2cxy(x+y).\]
         He adds that from these and other formulae of the same kind,  he deduced the following theorem (see Figure \ref{elliptic-letter}):          
         \begin{quote}\small
       If, in the quadrant $ACB$ of an ellipse, the tangent $VTM$ at an arbitrary point $M$ is drawn which meets one of the axes, $CB$, at $T$, if $TV$ is taken equal to $CA$ and from $V$, $VN$ is drawn parallel to $CB$, and if finally $CP$ is the perpendicular on the tangent through the center $C$, then I say the difference of the arcs $BM$ and $AN$ will be rectifiable, namely, $BM-AN=MP$.
         \end{quote}
            
%            \medskip
%            
%              \begin{figure}
%\centering
%\includegraphics[width=0.5\linewidth]{elliptic.pdf}    \caption{\small {A picture from Euler's calculation of the length of a segment of an ellipse, from Euler's letter to Goldbach dated May 30, 1752.}}   \label{elliptic-letter}  
%\end{figure}
% 

              In the following letter to Goldbach, dated June 3rd, 1752, Euler gave a proof of  this theorem and clarified a formula that Fagnano had given in his 1716 paper \cite{Fagnano1716}.

About five weeks after Euler\index{Euler, Leonhard (1707--1783)} received the work of Fagnano, he presented to the Berlin Academy a memoir entitled  \emph{Observationes de comparatione arcuum curvarum irrectificibilium}
(Observations on the comparison of arcs of irrectifiable curves) \cite{E252} in which he expands on what he had announced in his correspondence with Goldbach, generalizing Fagnano's duplication result on the lemniscate  to a general multiplication result and giving examples of arcs of an ellipse, hyperbola and lemniscate whose differences are rectifiable. 
 This was the beginning of a systematic study by Euler of elliptic integrals. The year after, he presented to the Saint Petersburg Academy of Sciences a 
memoir entitled 
 \emph{De integratione aequationis differentialis $\displaystyle \frac{mdx}{\sqrt{1-x^4}}=\frac{ndy}{\sqrt{1-y^4}}$}
(On the integration of the differential equation $\displaystyle \frac{mdx}{\sqrt{1-x^4}}=\frac{ndy}{\sqrt{1-y^4}}$) \cite{E251} which starts with the sentence:\footnote{The translation is by S. G. Langton.} 
\begin{quote}\small
When, prompted by the illustrious Count Fagnano, I first considered this equation, I found indeed an algebraic relation between the variables $x$ and $y$ which satisfied the equation.
\end{quote} Several years later, in his famous treatise \emph{Institutiones calculi integralis} \cite{E342}, Euler\index{Euler, Leonhard (1707--1783)} included a section on the addition and multiplication of integrals of the form $$\displaystyle \int\frac{Pd Z}{\sqrt{A+2BZ+CZ^2+2DZ^3+EZ^4}}.$$

Fagnano's works, in three volumes, were edited in 1911--1912 by Gambioli, Loria and Volterra \cite{Fagnano1911}.

          Among the large number of memoirs that Euler wrote on elliptic integrals,\footnote{The Euler\index{Euler, Leonhard (1707--1783)} archive lists thirty-three memoirs by him under the heading ``Elliptic integrals," published between 1738 and 1882. It is sometimes hard to know exactly the year where Euler\index{Euler, Leonhard (1707--1783)} wrote his memoirs. For several of them, the date of publication was much later that the date of writing, and there are several reasons for that, including the huge backlog of the publication department of the Academies of Sciences of Saint Petersburg and Berlin, the main reason being that Euler used to send them too many papers.} we mention the short memoir \cite{E211}, \emph{Problema, ad cuius solutionem geometrae invitantur; theorema, ad cuius demonstrationem geometrae invitantur} 
(A Problem, to which a geometric solution is solicited; a theorem, to which a geometric proof is solicited),  published in 1754, containing his  result on the rectification of the difference of two arcs of an ellipse. We also mention the memoir \cite{E264}, \emph{Demonstratio theorematis et solutio problematis in actis erud. Lipsiensibus propositorum} 
(Proof of a theorem and solution of a theorem proposed in the Acta Eruditorum of Leipzig)  \cite{E264}, in which he\index{Euler, Leonhard (1707--1783)} studies  the division by 2 of an arc of ellipse.  The memoir \cite{E345}, entitled \emph{Integratio aequationis $\frac{d x}{\sqrt{\alpha+\beta x+\gamma x^2+\delta x^3+\epsilon x^4}}=\frac{d y}{\sqrt{\alpha+\beta y+\gamma y^2+\delta y^3+\epsilon y^4}}$} (The integration of the equation $\frac{d x}{\sqrt{\alpha+\beta x+\gamma x^2+\delta x^3+\epsilon x^4}}=\frac{d y}{\sqrt{\alpha+\beta y+\gamma y^2+\delta y^3+\epsilon y^4}}$),  written in 1765 and published in 1768, is mentioned by Jacobi\index{Jacobi Carl Gustav Jacob (1804--1851)} in a letter to Legendre\index{Legendre, Adrien-Marie (1752--1833)} which we quote below.
    
             Besides Euler, one may mention d'Alembert.\index{Alembert@d'Alembert, Jean le Rond (1717--1783)} In a letter dated December 29, 1746, Euler\index{Euler, Leonhard (1707--1783)} writes to his Parisian colleague (see \cite{Euler-Alembert} p. 251):
            
             \begin{quote}\small
             I read with as much profit as satisfaction your last piece with which you honored our Academy. [...] But what pleased me most in your piece is the reduction of several integral formulae to the rectification of the ellipse and the hyperbola; a matter to which I had also already given my thoughts, but I was not able to get entirely to the formula
             \[\frac{dx}{\sqrt{\alpha+\beta x+\gamma x^2+\delta x^3+ \epsilon x^4}}\]
             and I regard your formula as a masterpiece of your expertise.\footnote{J'ai lu avec autant de fruit que de satisfaction votre derni\`ere pi\`ece dont vous avez honor\'e notre acad\'emie. [...] Mais ce qui m'a plu surtout dans votre pi\`ece c'est la r\'eduction de plusieurs formules int\'egrales \`a la rectification de l'ellipse et de l'hyperbole ; mati\`ere \`a laquelle j'avais aussi d\'ej\`a pens\'e, mais je n'ai pu venir \`a bout de la formule
             \[\frac{dx}{\sqrt{\alpha+\beta x+\gamma x^2+\delta x^3+ \epsilon x^4}}\]
             et je regarde votre formule comme un chef-d'\oe uvre de votre expertise.}
             \end{quote}

  Lagrange,\index{Lagrange, Joseph-Louis (1736--1813)} whose name is associated with that of Euler in several contexts, studied elliptic integrals in his famous \emph{Th\'eorie des fonctions analytiques} (Theory of analytic functions) \cite{Lagrange-theorie} (first edition 1797).  In particular, he discovered a relation between Euler's addition formula and a problem in spherical trigonometry.

    After Euler\index{Euler, Leonhard (1707--1783)}, d'Alembert and Lagrange, we  must talk about Legendre,\index{Legendre, Adrien-Marie (1752--1833)}  who investigated these integrals for almost forty years. He wrote two famous  treatises on the subject, his \emph{Exercices de calcul int\'egral sur divers ordres de transcendantes et sur les quadratures} (Exercises of integral calculus on various orders of transcendence and on the quadratures) \cite{Legendre-Exercices}  (1811--1816)  and his \emph{Trait\'e des fonctions elliptiques et des int\'egrales eul\'eriennes} (Treatise of elliptic functions and Eulerian integrals) \cite{Legendre-Traite} (1825--1828), both in three volumes. In the introduction to the latter (p. 1ff.), Legendre\index{Legendre, Adrien-Marie (1752--1833)} makes a brief history of the subject, from its birth until the moment he started working on it. According to his account, elliptic functions were first studied by MacLaurin\index{Maclaurin, Colin (1698--1764)} and d'Alembert\index{Alembert@d'Alembert, Jean le Rond (1717--1783)} who found several formulae for integrals that represent arcs of ellipses or arcs of hyperbolas.\footnote{See e.g. \cite{MacLaurin} and \cite{Alembert1846}.} Legendre declares that their results were too disparate to form a theory. He then mentions Fagnano, \index{Fagnano, Giulio Carlo de' Toschi di (1682--1766)} recalling that his work was the starting point of the profound analogy between elliptic integrals and trigonometric functions. After describing Fagnano's work, Legendre\index{Legendre, Adrien-Marie (1752--1833)} talks about some of the main contributions of Euler\index{Euler, Leonhard (1707--1783)}, Lagrange\index{Lagrange, Joseph-Louis (1736--1813)} and Landen on the subject. His treatise starts with a detailed study of integrals of the form $\displaystyle \int\frac{Pdx}{R}$ investigated by Euler, where $P$ is an arbitrary rational function of $x$ and $R=\displaystyle\sqrt{\alpha+\beta x+\gamma x^2+\delta x^3+\epsilon x^4}$. The expression \emph{Eulerian integral} contained in the title of Legendre's treatise was coined by him. He writes:
    \begin{quote}\small
    Although Euler's name is attached to almost all the important theories  of integral calculus, I nevertheless thought that I was allowed to give more especially the name \emph{Eulerian integral} to two sorts of transcendants whose properties constituted the subject of several beautiful memoirs of Euler, and form the most complete theory on definite integrals which exists up to now [...]
    \footnote{Quoique le nom d'Euler soit attach\'e \`a presque toutes les th\'eories importantes du calcul int\'egral, cependant j'ai cru qu'il me serait permis de donner plus sp\'ecialement le nom d'\emph{Int\'egrales Eul\'eriennes} \`a deux sortes de transcendantes dont les propri\'et\'es ont fait le sujet de plusieurs beaux m\'emoires d'Euler, et forment la th\'eorie la plus compl\`ete que l'on connaisse jusqu'\`a pr\'esent sur les intagrales d\'efinies [...]}
    \end{quote}
    
After Legendre, and among the immediate predecessors of Riemann on elliptic functions, we find Abel,\index{Abel, Nils Henrik (1802--1829)} Jacobi,\index{Jacobi Carl Gustav Jacob (1804--1851)} and Gauss.\index{Gauss, Carl Friedrich (1777--1855)} The last two were his teachers in Berlin and G\"ottingen respectively. With this work, the emphasis in the study of elliptic integrals shifted to that of their inverses. Considering inverses is naturally motivated by the analogy with trigonometric functions, as one may see by recalling that the integral $\displaystyle \int_0^x\frac{dt}{\sqrt{1-t^2}}$ represents the arcsine function, and therefore, its inverse is the more tractable sine function. The periodic behavior of inverses of elliptic integrals like $\displaystyle \int_0^x\frac{dt}{\sqrt{1-t^4}}$ and others, which became later one of the main questions in that theory, is in some sense a generalization of that of trigonometric functions.

Abel\index{Abel, Nils Henrik (1802--1829)} and Jacobi\index{Jacobi Carl Gustav Jacob (1804--1851)} developed simultaneously the theory of elliptic integrals, and separating their results has always been a difficult task. It is also well established that Gauss  discovered  several results of Abel\index{Abel, Nils Henrik (1802--1829)} and Jacobi\index{Jacobi Carl Gustav Jacob (1804--1851)} before them, but never published them. This is attested in his notebook and in his correspondence, published in his Collected Works. Gauss started his notebook in 1796, at the age of 19, and he wrote his last note there in 1814. The notes consist of 146 statements, most of them very concise, and they fill up a total of 20 pages in his  \emph{Collected Works} (vol. 10). This  edition of the notebook published in Gauss's\index{Gauss, Carl Friedrich (1777--1855)} Collected Works is accompanied by detailed comments by Bachmann, Brendel, Dedekind, Klein, L\oe wy, Schlesinger and St\"ackel. There is a French translation of the notebook \cite{Gauss-cahier}. Among the notes contained in this diary, several concern elliptic functions. For instance, in Notes 32 and 33, Gauss studies the inverse of the lemniscate integral $\int \frac{dx}{\sqrt{1-x^4}}$, as a particular case of the elliptic integral  $\int \frac{dx}{\sqrt{1-x^n}}$. In Note 53, he mentions that he is studying the general integral  $\int \frac{dx}{\sqrt[n]{1-x^n}}$, which was already considered by Euler in his \emph{Institutiones calculi integralis}. In Note 54, he states that he has an easy method for determining the integral  $\int \frac{x^ndx}{1+x^m}$, again an integral that was considered by Euler. There are several other notes on elliptic integrals in Gauss's\index{Gauss, Carl Friedrich (1777--1855)} notebook.

Jacobi\index{Jacobi Carl Gustav Jacob (1804--1851)} read Euler's works while he was in high school. He obtained his PhD at the age of 21, and at the age of 22, he started a correspondence with Legendre,\index{Legendre, Adrien-Marie (1752--1833)} who was 74, informing him about his results on elliptic integrals.  This correspondence became famous. It is reproduced in Crelle's Journal\footnote{Crelle's Journal, 80 (1875), p. 205--279}  and in  Jacobi's \emph{Collected Works}.\footnote{\emph{Collected Works}, t. I, p. 385--46.} The beginning of this  correspondence is touching.  Jacobi\index{Jacobi Carl Gustav Jacob (1804--1851)} sends his first letter to Legendre \index{Legendre, Adrien-Marie (1752--1833)} on August 5, 1827, expressing his great respect for the work of his older French colleague. He writes  (\cite{Jacobi-Collected} vol. 1, p. 390): 
\begin{quote}\small
A young geometer dares to present you a few discoveries in the theory of elliptic functions, to which he was led by a diligent study of your beautiful writings.  It is to you, Sir, that this brilliant part of analysis owes the highest degree of perfection to which it has been elevated, and it is only in following the footsteps of such a great master that the geometers will be able to push it beyond limits which have been so far prescribed. Thus, it is to thee that I must offer the following, as a fair tribute of admiration and gratefulness.\footnote{Un jeune g\'eom\`etre ose vous pr\'esenter quelques d\'ecouvertes faites dans la th\'eorie des fonctions elliptiques, auxquelles il a \'et\'e conduit par la lecture assidue de vos beaux \'ecrits. C'est \`a vous, Monsieur, que cette partie brillante de l'analyse doit le haut degr\'e de perfectionnement auquel elle a \'et\'e port\'ee, et ce n'est qu'en marchant sur les vestiges d'un si grand ma\^\i tre, que les g\'eom\`etres pourront parvenir \`a la pousser au-del\`a des bornes qui lui ont \'et\'e prescrites jusqu'ici. C'est donc \`a vous que je dois offrir ce qui suit comme un juste tribut d'admiration et de reconnaissance.}
\end{quote}

In his response, dated November 30, 1827, Legendre,\index{Legendre, Adrien-Marie (1752--1833)} referring to one of the theorems that Jacobi\index{Jacobi Carl Gustav Jacob (1804--1851)} communicated to him, writes (\cite{Jacobi-Collected} vol. 1, p. 396):
\begin{quote}\small
I checked this theorem by my own methods and I found it perfectly correct. Even though I regret that this discovery escaped me, the joy I experienced was most vivid when I saw the significant improvement that was added to the beautiful theory of which I am the creator and which I developed almost alone during more than forty years.\footnote{J'ai v\'erifi\'e ce th\'eor\`eme par les m\'ethodes qui me sont propres et je l'ai trouv\'e parfaitement exact. En regrettant que cette d\'ecouverte m'ait \'echapp\'ee je n'en ai pas moins \'eprouv\'e une joie tr\`es vive de voir un perfectionnement si notable ajout\'e \`a la belle th\'eorie, dont je suis le cr\'eateur, et que j'ai cultiv\'e presque seul depuis plus de quarante ans.}
\end{quote}

In another letter, sent on January 12, 1828, Jacobi\index{Jacobi Carl Gustav Jacob (1804--1851)} informs Legendre\index{Legendre, Adrien-Marie (1752--1833)} about Abel's discoveries, in particular on the division of the lemniscate  (\cite{Jacobi-Collected}, vol. 1, p. 401): 
\begin{quote} \small
Since my last letter, researches of the highest importance were published on elliptic functions by a young geometer, who may be personally known to you.\footnote{Depuis ma derni\`ere lettre, des recherches de la plus grande importance ont \'et\'e publi\'ees sur les fonctions elliptiques de la part d'un jeune g\'eom\`etre, qui peut-\^etre vous sera connu personnellement.}
\end{quote}

Legendre\index{Legendre, Adrien-Marie (1752--1833)} sent his response on February 9, informing his correspondent that he knew about Abel's work, but that he was happy to see it summarized in a language which was closer to his own.\footnote{[J'avais d\'ej\`a connaissance du beau travail de M. Abel\index{Abel, Nils Henrik (1802--1829)} ins\'er\'e dans le \emph{Journal de Crelle}. Mais vous m'avez fait beaucoup de plaisir de m'en donner une analyse dans votre langage qui est plus rapproch\'e du mien.] (\cite{Jacobi-Collected}, t. 1, p.  407).}

Regarding Gauss's work on elliptic functions, we mention an excerpt of the first letter from Jacobi to Legendre \cite{Jacobi-Collected} p. 393--394:

\begin{quote}\small
These researches were born only very recently. However, they are not the only ones that are conducted in Germany on the same object. Mr. Gauss, when he learned about them, informed me that he had developed, already in 1808, the cases of  3 sections, 5 sections and 7 sections, and that he found at the same time the corresponding new scales of modules. It seems to me that this information is very interesting.\footnote{Il n'y a que tr\`es peu de temps que ces recherches ont pris naissance. Cependant elles ne sont pas les seules entreprises en Allemagne sur le m\^eme objet. M. Gauss, ayant appris de celles-ci, m'a fait dire qu'il avait d\'evelopp\'e d\'ej\`a en 1808 les cas de 3 sections, 5 sections et de 7 sections, et trouv\'e en m\^eme temps les nouvelles \'echelles de modules qui s'y rapportent. Cette nouvelle, \`a ce qui me para\^\i t, est bien int\'eressante.}
\end{quote}

Legendre\index{Legendre, Adrien-Marie (1752--1833)} was outraged by Gauss's reaction. In his response to Jacobi, dated November 30, 1827, he writes (\cite{Jacobi-Collected} p. 398):
\begin{quote}\small 
How is it possible that Mr. Gauss dared telling you that most of your theorems were known to him and that he discovered them back in 1808? This excess of impudence is unbelievable from a man who has enough personal merit so as he does not need to appropriate the discoveries of others... But this is the same man who, in 1801, wanted to attribute to himself the law of reciprocity published in 1785 and who wanted, in 1809, to take hold of the method of least squares that was published in 1805.\footnote{Comment se fait-il que M. Gauss\index{Gauss, Carl Friedrich (1777--1855)} ait os\'e vous dire que la plupart de vos th\'eor\`emes lui \'etaient connus et qu'il en avait fait la d\'ecouverte d\`es 1808 ? Cet exc\`es d'impudence n'est pas croyable de la part d'un homme qui a assez de m\'erite personnel pour n'avoir pas besoin de s'approprier les d\'ecouvertes des autres... Mais c'est le m\^eme homme qui en 1801 voulut s'attribuer la d\'ecouverte de la loi de r\'eciprocit\'e publi\'ee en 1785 et qui voulut s'emparer en 1809 de la m\'ethode des moindres carr\'es publi\'ee en 1805.
}
\end{quote}

It was only at the publication of Gauss's \emph{Collected Works},\footnote{Gauss's collected works, \emph{Carl Friedrich Gauss' Werke}, in twelve volumes, were published between 1863 and 1929.}  containing in particular his famous notebook,  that it became clear that Gauss's assertion concerning the fact that he had discovered before Abel\index{Abel, Nils Henrik (1802--1829)} most of the properties of elliptic functions, including their double periodicity, was correct. 
One of the first results of Abel\index{Abel, Nils Henrik (1802--1829)} concerns integrals of arcs of lemniscate,  a curve which he showed to be divisible by ruler and compass into $n$ equal parts, for the same values of $n$ for which the circle is divisible into $n$ equal parts. The same result was stated without proof in Gauss's\index{Gauss, Carl Friedrich (1777--1855)} \emph{Disquisitiones arithmeticae} \cite{Gauss-D}.\index{Gauss!Disquisitiones arithmeticae} 

Abel's first major results on elliptic functions are contained in his 1827 paper \emph{Recherches sur les fonctions elliptiques} (Researches on elliptic functions) \cite{Abel-1827}. He explains there the double periodicity of these functions, as well as their multiplication and division properties.  The analogy with circular functions is again highlighted. At the beginning of his paper, Abel\index{Abel, Nils Henrik (1802--1829)} talks about his famous predecessors, Euler\index{Euler, Leonhard (1707--1783)}, Lagrange and Legendre.\index{Legendre, Adrien-Marie (1752--1833)} He writes (p. 101):
\begin{quote}\small
The first idea of these [elliptic] functions were given by the immortal Euler, who showed that the separable equation
\[\frac{\partial x}{\sqrt{\alpha+\beta x+\gamma x^2+\delta x^3+\epsilon x^4}}+\frac{\partial y}{\sqrt{\alpha+\beta y+\gamma y^2+\delta y^3+\epsilon y^4}}=0\] is algebraically integrable. After Euler\index{Euler, Leonhard (1707--1783)}, Lagrange\index{Lagrange, Joseph-Louis (1736--1813)} added  something, when he gave his elegant theory of the transformation of the integral \[\int\frac{R dx}{\sqrt{(1-p^2x^2)(1-q^2x^2)}},\] where $R$ is a rational function of $x$. But the first, if I am not mistaken, who went thoroughly into the nature of these functions, is Mr. Legendre,\index{Legendre, Adrien-Marie (1752--1833)} who, first in a memoir on elliptic functions, and then in his excellent \emph{Exercices de math\'ematiques}, developed numerous elegant properties of these functions, and showed their usefulness.\footnote{La premi\`ere id\'ee de ces fonctions a \'et\'e donn\'ee par l'immortel Euler, en d\'emontrant que l'\'equation s\'epar\'ee
\[\frac{\partial x}{\sqrt{\alpha+\beta x+\gamma x^2+\delta x^3+\epsilon x^4}}+\frac{\partial y}{\sqrt{\alpha+\beta y+\gamma y^2+\delta y^3+\epsilon y^4}}=0\]
est int\'egrable alg\'ebriquement. Apr\`es Euler\index{Euler, Leonhard (1707--1783)}, Lagrange y a ajout\'e quelque chose, en donnant son \'el\'egante th\'eorie de la transformation de l'int\'egrale \[\displaystyle\int\frac{R dx}{\sqrt{(1-p^2x^2)(1-q^2x^2)}},\] o\`u $R$ est une fonction rationnelle de $x$. Mais le premier, et si je ne me trompe, le seul, qui ait approfondi la nature de ces fonctions, est M. Legendre, qui d'abord dans un m\'emoire sur les fonctions elliptiques, et ensuite dans ses excellents exercices de math\'ematiques, a d\'evelopp\'e nombre de propri\'et\'es \'el\'egantes de ces fonctions, et a montr\'e leur application.
}
\end{quote}

Riemann was already interested in elliptic functions while he was a student in Berlin. Klein, in his \emph{Development of mathematics in the 19th century} \cite{Klein-development} (Chapter VI) writes that the latter, since the end of the 1840s, was interested in elliptic functions because this subject was fashionable in Germany. From a letter to his father, dated May 30, 1849, we know that Riemann was following in Berlin Jacobi's\index{Jacobi Carl Gustav Jacob (1804--1851)} and Eisenstein's\index{Eisenstein, Ferdinand Gotthold Max (1823--1852)}  lectures on elliptic functions. He writes (cf. \cite{Riemann-Letters}): ``Jacobi has just begun a series of lectures in which he leads off once again with the entire system of the theory of elliptical functions in the most advanced, but elementary way." In another letter (without date), also written in Berlin, Riemann writes: ``I enrolled with five other students into a private class (Privatissumum) with Eisenstein,\index{Eisenstein, Ferdinand Gotthold Max (1823--1852)} who was promoted in the course of this semester to a Privatdozent with a paper on the theory of elliptic functions."

We already mentioned Euler's impact on Jacobi. Eisenstein is another prominent mathematician on which Euler exerted a crucial influence. In his biography of Eisenstein \cite{Schmitz}, M. Schmitz writes that during  the  period 1837--1842, while he was a Gynmasium pupil,  Eisenstein  attended  lectures  by Dirichlet at  the University  of  Berlin,  and that he studied on his own  Gauss's
\emph{Disquisitiones Arithmeticae} as well as  
 papers and books by Euler 
and Lagrange. We quote Eisenstein, from his autobiography translated in  \cite{Schmitz}: 
\begin{quote} \small After I  had  acquired  the  fundamentals  by  private  study  (I  never  had  a  private  
tutor)  I  proceeded  to  advanced mathematics and studied, besides other books containing advanced material, the brilliant work of Euler and Lagrange about differential and integral calculus. I was able to commit this material securely to my memory and to master it entirely, because I made it a rule to 
compose every theory in writing as 
soon as I understood it. 
\end{quote}

In his ICM communication \cite{Weil-ICM-H}, Weil\index{Weil, Andr\'e (1906--1998)} declares (p. 233) that ``Eisenstein fell in love with mathematics at an early age by reading Euler and Lagrange."

          We shall conclude this section with two other quotes of Weil. Before that, let us recall that elliptic integrals are studied in number theory in relation with the theory of elliptic curves. Weil\index{Weil, Andr\'e (1906--1998)} writes in an essay on the history of number theory, \cite{Weil-Ens}, p. 15, that Fermat, in his work on number theory, had already dealt with elliptic curves (without the name), in particular in his proof of the non-existence of integer solutions for the equation $x^4-y^4=z^2$. We quote him from his book on the history of number theory that we already mentioned (\cite{Weil-Hammurapi} p. 242):
              \begin{quote}\small
              What we call now ``elliptic curves" (i.e. algebraic curves of genus 1) were considered by Euler\index{Euler, Leonhard (1707--1783)} under two quite different aspects without ever showing an awareness of the connection between them, or rather of their substantial identity. 
              On the one hand, he must surely have been familiar, from the very beginning of his career, with the traditional methods for handling Diophantine equations of genus 1. [...] On the other hand he had inherited from his predecessors, and notably from Johann Bernoulli, a keen interest in what we know as ``elliptic integrals" because the rectification of the ellipse depends upon integrals of that type; they were perceived to come next to the integrals of rational functions in order of difficulty.
              \end{quote}

         Eisenstein\index{Eisenstein, Ferdinand Gotthold Max (1823--1852)} and Dirichlet\index{Dirichlet, Johann Peter Gustav Lejeune (1805--1859)} were mostly interested in elliptic functions because of their use in number theory, contrary to Riemann, who, even though he was introduced to elliptic functions through Eisenstein's lectures, was not excited by that field.  Weil\index{Weil, Andr\'e (1906--1998)} writes in his essays \cite{Weil-Ens}, p. 21:  
              \begin{quote}\small
              [...] The case of Riemann is more curious. Of all the great mathematicians of the last century, he is outstanding for many things, but also, strangely enough, for his complete lack of interest for number theory and algebra. This is really striking, when one reflects how close he was, as a student, to Dirichlet\index{Dirichlet, Johann Peter Gustav Lejeune (1805--1859)} and Eisenstein,\index{Eisenstein, Ferdinand Gotthold Max (1823--1852)} and, at a later period, also to Gauss and to Dedekind who became his most intimate friend. During Riemann's student days in Berlin, Eisenstein tried (not without some success, he fancied) to attract him to number theory. ln 1855, Dedekind was lecturing in G\"ottingen on Galois theory, and one might think that Riemann, interested as he was in algebraic functions, might have paid some attention.  But there is not the slightest indication that he ever gave any serious thoughts to such matters.
              \end{quote}

We shall mention the work of Dirichlet on number theory (in particular on the prime number theorem) in \S \ref{s:zeta} below.
              In Chapter 8 \cite{Papa-Riemann3} of the present book, we report on several treatises on elliptic functions that were published in France during the few decades that followed Riemann's early work on the subject. In the next section, we review the more general Abelian functions.

\section{Abelian functions} \label{s:Abelian}

A few years before Riemann started his work on elliptic functions and elliptic integrals, the general interest moved towards the more general Abelian integrals, and their inversion.\index{Jacobi inversion problem}\index{problem!Jacobi inversion} The term \emph{Abelian function}, first introduced by Jacobi in honor of Abel, 
 is generally given to the functions obtained by inverting an  arbitrary algebraic integral or a combination of such integrals.\index{Abelian integral} An algebraic integral is an integral of the form $\int R(x,y) dx$ where $R$ is a rational function of the two variables $x$ and $y$ and where $x$ and $y$ satisfy furthermore a polynomial equation $f(x,y)=0$. In his 1826 memoir submitted to the Paris Academy, Abel extended Euler's addition formula for elliptic integrals to Abelian integrals. He  proved that the sum of an arbitrary number of such integrals can be written as the sum of $p$ linearly independent integrals, to which is added an algebraic-logarithmic expression. Here $p$ is the so-called \emph{genus} of the algebraic curve defined by the equation $f(x,y)=0$. After he learned about Abel's work, Jacobi formulated a generalized inversion\index{Jacobi inversion problem}\index{problem!Jacobi inversion} problem for a system of $p$ hyperelliptic integrals. His ideas were pursued by several mathematicians, and in particular by Riemann, who gave a solution to the inversion problem in terms of $\vartheta$ functions.\index{Jacobi inversion problem}\index{problem!Jacobi inversion} 

 Abel\index{Abel, Nils Henrik (1802--1829)} also discovered that the inverse functions of elliptic integrals\index{elliptic integral} are doubly periodic functions defined on the complex plane. This property was at the basis of the later introduction of group theory in the theory of elliptic curves.

In the passage from elliptic functions to Abelian functions, one must also mention Galois. The day before his death, Galois sent a letter to his friend Auguste Chevalier in which he described his thoughts, saying that one could write a memoir based on his ideas on integrals. The letter is analyzed by Picard in his article \cite{Picard-Galois}.\footnote{This article constituted the preface to the Collected Works of Galois which were published shortly after.} Picard writes:
\begin{quote}\small
All what we know about these researches is contained in what he says in this letter. Several points remain obscure in some statements of Galois; however, we can have a precise idea of some of the results he reached in the theory of integrals of algebraic functions. We thus acquire the certainty that he possessed the most essential results on Abelian integrals\index{Abelian integral} that Riemann was led to obtain twenty-five years later. We see without surprise Galois talking about the periods of an Abelian integral relative to an arbitrary algebraic function [...] The statements are precise; the famous author makes the classification of Abelian integrals into three kinds, and he declares that if $n$ denotes the number of linearly independent integrals of the first kind, the number of periods is $2n$. The theorem relative to the parameter inversion in the integrals of the third type is clearly marked, as well as the relations with the periods of Abelian integrals. Galois also talks about a generalization of Legendre's\index{Legendre, Adrien-Marie (1752--1833)} classical equation where the periods of elliptic integrals appear, a generalization which probably led him to the important relation that was discovered later on by Weierstrass and Mr. Fuchs.\footnote{Nous ne connaissons de ces recherches que ce qu'il en dit dans cette lettre ; plusieurs points restent obscurs dans quelques \'enonc\'es de Galois, mais on peut cependant se faire une id\'ee pr\'ecise de quelques-uns des r\'esultats auxquels il \'etait arriv\'e dans la th\'eorie des int\'egrales de fonctions alg\'ebriques. On acquiert ainsi la conviction qu'il \'etait en possession des r\'esultats les plus essentiels sur les int\'egrales ab\'eliennes que Riemann devait obtenir vingt-cinq ans plus tard. Nous voyons sans \'etonnement Galois parler des p\'eriodes d'une int\'egrale ab\'elienne relative \`a une fonction alg\'ebrique quelconque [...] Les \'enonc\'es sont pr\'ecis ; l'illustre auteur fait la classification en trois esp\`eces des int\'egrales ab\'eliennes, et affirme que, si $n$ d\'esigne le nombre des int\'egrales de premi\`ere esp\`ece lin\'eairement ind\'ependantes, les p\'eriodes seront en nombre $2n$. Le th\'eor\`eme relatif \`a l'inversion du param\`etre dans les int\'egrales de troisi\`eme esp\`ece est nettement indiqu\'e, ainsi que les relations entre les p\'eriodes des int\'egrales ab\'eliennes ; Galois parle aussi d'une g\'en\'eralisation de l'\'equation classique de Legendre,\index{Legendre, Adrien-Marie (1752--1833)} o\`u figurent les p\'eriodes des int\'egrales elliptiques, g\'en\'eralisation qui l'avait probablement conduit \`a l'importante relation d\'ecouverte depuis par Weierstrass et par M. Fuchs.}
\end{quote}

              In his paper on Abelian functions \cite{Riemann-Abelian}, Riemann establishes existence results for Abelian functions and more generally their determination in terms of the points of discontinuity and the information on the ramification at these points. It is in that paper that Riemann introduces the notion of birational equivalence and number of moduli both of which played an essential role in mathematics. In the same paper, he presents Abel's addition theorem for elliptic integrals, and he solves Jacobi's inversion problem in terms of $p$ variable magnitudes, for a $(2p+2)$-connected surface.\index{Jacobi inversion problem}\index{problem!Jacobi inversion} It is also in this paper that Riemann gives his well known classification of Abelian integrals into three types, a classification which depends on the existence and the nature of the singularities (poles or logarithmic). Riemann mentions in his paper several works on the inversion problem, in particular the successful attempt by Weierstrass in  the case of hyperelliptic integrals.

   On the work of Riemann on Abelian integrals\index{Abelian integral}, the reader is also referred to Chapter 4, by Houzel, in the present volume \cite{Houzel}. For a comprehensive survey on the work of Abel,\index{Abel, Nils Henrik (1802--1829)} the interested reader is referred to the article \cite{Houzel-Abel} by Houzel.

\section{Hypergeometric series}\label{s:hyper}

The theory of the hypergeometric series is another topic which Riemann tackled and whose roots involve in an essential way the works of Euler and Gauss. Riemann's main paper on the subject is  \emph{Beitr\"age zur Theorie der durch die Gauss'sche Reihe $F(\alpha,\beta,\gamma,x)$ darstellbaren Functionen}  (Contribution to the theory of functions representable by Gauss's series $F(\alpha,\beta,\gamma,x)$) \cite{Riemann-Beitrage}, published in 1857. The work in this paper was used by Riemann later in his development of the theory of analytic differential equations. There are  also fragments on the same subject published in Riemann's \emph{Collected works}. 

The hypergeometric series is a function of the form 
\[F(\alpha,\beta,\gamma,x)\]
\[=1+\frac{\alpha\beta}{1.\gamma}x
+\frac{\alpha(\alpha+1)\beta(\beta+1)}{1.2\gamma(\gamma+1)}x^2  
+\frac{\alpha(\alpha+1)(\alpha+2)\beta(\beta+1)(\beta+2)}{1.2.3\gamma(\gamma+1)(\gamma+2)}x^3+\ldots
 \]
where $x$ is the variable.

The term ``hypergeometric series" appears in Euler's \emph{Institutiones calculi integralis} \cite{E342} (1769), Chapter XI. 
 The series is a solution of the so-called Euler hypergeometric differential equation which appears in Chapters VIII and XI of the same treatise. As a matter of fact, this name was given to several different but closely related objects. Euler, in one of his earliest memoir \emph{ De progressionibus transcendentibus seu quarum termini generales algebraice dari nequeunt}
(On transcendental progressions, that is, those whose general terms cannot be given algebraically)  \cite{E19}, published in 1738,  starts by mentioning Wallis's\index{Wallis, John (1616--1703)} ``hypergeometric series" $1! + 2! + 3! + 4! + ...$ (without the factorial notation). The terminology here refers to the fact that in analogy with the case of  geometric progressions, where each term is obtained from the preceding one by multiplying it by a constant, one defined a hypergeometric progression as a progression in which each term is obtained from the preceding one by multiplying it by a factor which increases by a unit at each step. Wallis's papers on this subject include \cite{Wallis1} (1, Scholium to Proposition 190) and \cite{Wallis2} (p. 315).

Gauss mentions a hypergeometric series in his doctoral dissertation\index{doctoral dissertation!Gauss}\index{Gauss!doctoral dissertation} \emph{Demonstratio nova theorematis omnem functionem algebraicam rationalem integram unius variabilis in factores reales primi vel secundi gradus resolvi posse} \cite{Gauss-Demostratio} (New proof of the theorem that every rational integral algebraic function of one variable can be resolved into real factors of the first or second degree) (1799).

We refer the reader to the paper \cite{Dutka} for a comprehensive history of the hypergeometric series.

At the beginning of his announcement of his memoir  \cite{Riemann-Beitrage}, Riemann states: ``This memoir treats a class of functions which are useful to solve various problems in mathematical physics." As a matter of fact, these functions are still commonly used today in mathematical physics. 
Riemann notes that the name \emph{hypergeometric series}\index{hypergeometric series} was first proposed by Pfaff,\index{Pfaff, Johann Friedrich (1765--1825)} for a more general series, whereas Euler, after Wallis,\index{Wallis, John (1616--1703)} used such a name for a series which is slightly different. 
Pfaff was Gauss's friend, and had been his teacher. He studied this function in his book \emph{Disquisitiones analyticae maxime ad calculum integralem et doctrinam serierum pertinentes} (Analytic investigations most relevant for integral calculus and the doctrine of series) \cite{Pfaff-D} (1797).  Gauss has a series of unpublished results on the hypergeometric series, which he communicated to the astronomer Bessel, who was also his friend, in a  letter dated September 3, 1805. The results were used by Gauss in his later works. In his writings on the subject, Gauss used continued fractions in his study of the quotient of two hypergeometric series. He developed these ideas in his paper
\emph{Disquisitiones generales circa serium 
$1+\frac{\alpha\beta}{1.\gamma}x
+\frac{\alpha(\alpha+1)\beta(\beta+1)}{1.2\gamma(\gamma+1)}xx  
+\frac{\alpha(\alpha+1)(\alpha+2)\beta(\beta+1)(\beta+2)}{1.2.3\gamma(\gamma+1)(\gamma+2)}x^3+\ldots$ etc.} (General investigations on the series 
$1+\frac{\alpha\beta}{1.\gamma}x
+\frac{\alpha(\alpha+1)\beta(\beta+1)}{1.2\gamma(\gamma+1)}xx  
+\frac{\alpha(\alpha+1)(\alpha+2)\beta(\beta+1)(\beta+2)}{1.2.3\gamma(\gamma+1)(\gamma+2)}x^3+\ldots$ etc.)
 \cite{Gauss-Hy}. The same year, he wrote another paper on the same subject which he never published but which is contained in his \emph{Collected Works} edition \cite{Gauss-H}. 
 Riemann, in his paper \cite{Riemann-Beitrage}, proved that these fractions converge in the complex plane cut along the subset $[2,+\infty]$ of the $x$-axis. In the same memoir, he introduced in the study of the hypergeometic functions a new method, which applies to all functions that satisfy  linear differential equations with  algebraic coefficients. 
He recalls in the announcement of that memoir, published in the \emph{G\"ottinger Nachrichten}, No. 1, 1857,  that Euler and Gauss made a thorough study of these functions from the theoretical point of view.

In the introduction to his paper \cite{Gauss-Hy} Gauss declares that practically any transcendental function that appears in analysis may be obtained as a special case of the hypergeometric series.  In fact, it is known that functions like $\log (1+z)$, $\arcsin z$ and several orthogonal polynomials, including Legendre polynomials and Chebyshev polynomials, can be expressed using hypergeometric functions.
The so-called confluent hypergeometric function (or Kummer's function) is a limit of the hypergeometric function. 

 The introduction of the hypergeometric series brought a whole new class of new functions to the field of analysis which, at least in the times of Euler, consisted in the study of functions.

              \section{The zeta function} \label{s:zeta}

This section is concerned with Riemann's article
 \emph{\"Uber die Anzahl der Primzahlen unter einer gegebenen Gr\"osse}  (On the number of primes less than
 a given magnitude) \cite{Riemann-primes}. This memoir, which is only 8 pages long, changed the course of mathematics.  Riemann wrote it at the occasion of his election to the Berlin Academy of Sciences, on August 11, 1859. Every newly elected member at that academy was asked to report on his most recent research, and Riemann chose this topic.  A short history of the subject will show that the list of predecessors of Riemann in this field includes names which are familiar to us now: Euler, as always, then Legendre,\index{Legendre, Adrien-Marie (1752--1833)} Dirichlet\index{Dirichlet, Johann Peter Gustav Lejeune (1805--1859)} and Gauss.
  
 Riemann starts his memoir by recalling that Gauss and Dirichlet had been interested in this subject several years before him. He displays the following formula, which he recalls was noted by Euler, and which was his own departure point :
              \[\prod \frac{1}{1-\frac{1}{p^s}}=\sum\frac{1}{n^s}.\]
           Here, $p$ takes all the prime values and $n$ all the integer values.  Riemann    considers the function represented by these two expressions as a function of a complex variable $s$ as long as the two series converge, and he denotes this function by $\zeta(s)$.\footnote{Even though the notation $\zeta(s)$ and the name zeta function\index{zeta function}\index{Riemann!zeta function} first appear in Riemann's paper, we shall commit the usual anachronism of using the notation $\zeta(s)$ for the series $\sum_{n=1}^\infty \frac{1}{n^s}$ even when we talk about the work done on this series before Riemann.}  He then gives an integral  formula for this function\index{zeta function}\index{Riemann!zeta function}, and he notes that this integral is  ``uniform" (uni-valued), that it is defined and finite  for any  value of $s$ except for $s=1$  and  that it vanishes when $s$ is  a negative odd integer.

              The distribution of primes, which is the subject of Riemann's paper, may be traced back to Greek antiquity.   The reader may recall that there are several results on prime numbers in Euclid's \emph{Elements}.\index{Euclid!Elements} In particular, Proposition 20 of Book IX says that there are infinitely many primes. It is also known that the Greeks had a method to list effectively the sequence of primes (Eratosthenes sieve). 
Without any doubt, the general question of the distribution of primes kept busy the  mathematicians of that epoch.
 It is also good to recall, right at the beginning, that Euler, in his paper \emph{Variae observationes circa series infinitas} (Various observations about infinite series) \cite{E72}, showed that the series of inverses of primes,
              \[\frac{1}{2}+\frac{1}{3}+\frac{1}{5}+\frac{1}{7}+\frac{1}{11}+\ldots
              ,\] diverges, which in some sense is a wide generalization of the fact that the number of primes is infinite.

 Euler was fascinated by the question of the distribution of primes. We quote him from a paper entitled \emph{D\'ecouverte d'une loi tout extraordinaire des nombres par rapport \`a la somme de leurs diviseurs} (Discovery of a very extraordinary law of numbers in relation to the sum of their divisors) \cite{E175}, written in 1747 and published in 1751: 
 \begin{quote}\small    

Mathematicians tried in vain, until now, to discover some or other order in the sequence of prime numbers, and we have reasons to think that this is a mystery which human mind will never be able to penetrate. To be convinced, it suffices to take a look at the tables of prime numbers, that a few persons have taken the trouble to continue beyond one hundred thousand: one will primarly notice that there is no order and no rule there.\footnote{Les math\'ematiciens ont t\^ach\'e jusqu'ici en vain \`a d\'ecouvrir un ordre quelconque dans la progression des nombres premiers, et on a lieu de croire, que c'est un myst\`ere auquel l'esprit humain ne saurait jamais p\'en\'etrer. Pour s'en convaincre, on n'a qu'\`a jeter les yeux sur les tables des nombres premiers, que quelques personnes se sont donn\'e la peine de continuer au-del\`a de cent mille : et on s'apercevra d'abord qu'il ne r\`egne aucun ordre ni r\`egle.}

\end{quote}

Let us return now to the zeta function.

The history of the zeta function\index{zeta function}\index{Riemann!zeta function} in Euler's works naturally starts with the question of the value of the sum of the series of reciprocals of squares, $\zeta(2)=\displaystyle \sum_{1}^\infty \frac{1}{n^2}$. Before Euler, this series was known to be convergent, and the determination of its value was an open question whose formulation can be traced back at least to Pietro Mengoli\index{Mengoli, Pietro (1625--1686)}  in his treatise \emph{Novae quadrature arithmeticae, seu de additione fractionum}\footnote{Mengoli's  treatise is entirely devoted to the theory of infinite series, despite the word \emph{quadrature} (that is, computation of areas) in the title.}  (New arithmetic quadratures, or the addition of fractions) \cite{Mengoli} (1650). Several mathematicians worked on the problem, including Wallis,\index{Wallis, John (1616--1703)} Leibniz, Stirling, de Moivre,\index{Moivre@de Moivre, Abraham (1667--1754)} Goldbach and several Bernoullis.  In fact, the question of computing infinite sums was already a fashionable subject at that epoch. Mengoli, Huygens and Leibniz\index{Leibniz, Gottfried Wilhelm (1646--1716)} independently computed the sum of reciprocals of the triangular numbers, that is, numbers of the form $\frac{(n)(n+1)}{2}$.   Leibniz's  computation of the series of inverses of triangular numbers uses the classical ``telescopic method" known to students, so its level of difficulty has nothing to do with Euler's computation of $\zeta(2)$.  The problem of finding the value of $\zeta(2)$ became widely known among mathematicians after it was asked explicitly by Jakob Bernoulli\index{Bernoulli, Jakob (1654--1705)} in his series of papers \emph{Positiones de seriebus infinitis} (Positions of an infinite series) (1689).\footnote{See the comments of this work of Bernoulli in Weil's article \cite{Weil1989} p. 4.} In the same work, Bernoulli considered the series for an arbitrary rational number $s$.

Euler published several papers on various aspects of the zeta function.\index{zeta function}\index{Riemann!zeta function} In particular, he was the first to discover a formula establishing a relation between this series and prime numbers.  It is interesting to recall that Euler has been investigating the convergence of infinite series and infinite products
since his early days as a mathematician.\footnote{One should note that power series representations of functions already appear in the works of Newton,\index{Newton, Isaac (1643--1727)} in the 1660s.} His first letter addressed to Goldbach, dated October 13, 1729,  concerns the $\Gamma$ function, a function that interpolates the factorials. Goldbach had asked the opinion of several mathematicians on that problem. Euler writes \cite{Euler-Goldbach}:\footnote{In this volume of the \emph{Opera Ominia}, the letters are translated into English.}
\begin{quote}\small
When lately I came across a few ideas that apparently could contribute to the
interpolation of series having a variable law -- as you are wont to call it -- I took
a closer look and discovered many things regarding that subject. As Mr. Bernoulli
hinted that these results might please you, Sir, I decided to write to you and
submit them to your judgment. For the series 1, 2, 6, 24, 120, \ldots , which you
have treated extensively, as I see,  I have found the general term [...]
\end{quote}

The letter ends with:
\begin{quote}\small
You, Sir, who have already enriched the theory of series by
so many important discoveries, will therefore judge for yourself what else may be
expected from this novel way to deal with series. It would certainly acquire its
greatest utility and perfection if you could bring yourself to investigate how the
differential calculus can be most conveniently applied to these questions. For up to
now my method has the drawback that I cannot find what I want, but rather have
to be content with wanting what I find. 
\end{quote}

      In his paper \emph{De summatione innumerabilium progressionum} (The summation of an innumerable progression) \cite{E20},  Euler starts by giving a 7-digit approximate value of $\zeta(2)$, namely, 1.644934.  Needless to say, such a computation needed from his part a large amount of computing, because the series converges very slowly. Before that, Wallis\index{Wallis, John (1616--1703)} had given, in his  \emph{Arithmetica infinitorum} (Arithmetic of the infinite), 1655, a 3-digit approximation of that series.   Goldbach and Daniel Bernoulli also gave 3-digit approximations, in 1728.\footnote{In a letter to Goldbach, sent in 1728, Daniel Bernoulli writes that the value of the series $\zeta(2)$ ``is very nearly $8/5$, and Goldbach answers that $\zeta(2)-1$ lies between $16233/25200$ and $30197/46800$; cf. Weil \cite{Weil-Hammurapi}  p. 257 for more details on this history.} The reader may find interesting information on that subject in the correspondence between Euler, Bernoulli and Goldbach.

      In 1735, Euler,\index{Euler, Leonhard (1707--1783)} who was 28 years old, obtained the summation formula for $\zeta(2)$ and, more generally, for the infinite series  $\zeta(2\nu)=\displaystyle \sum_{1}^\infty \frac{1}{n^{2\nu}}$ for any positive integer $\nu$. He found the values $\zeta(2)= \pi^2/6$ and $\zeta(2\nu)=r_\nu \pi^{2\nu}$, where $r_\nu$ are rational numbers which are closely related to the Bernoulli numbers. In the introduction to his memoir \emph{De summis serierum reciprocarum} (On the sums of series of reciprocals) \cite{E41} (1735), he writes:\footnote{The translation from the Latin is by Andr\'e Weil,\index{Weil, Andr\'e (1906--1998)} \cite{Weil-Hammurapi} p. 261.}
\begin{quote}\small
So much work has been done on the series $\zeta(n)$ that it seems hardly likely that anything new about them may still turn up ... I too, in spite of repeated efforts, could achieve nothing more than approximate values for their sums ... Now, however, quite unexpectedly, I have found an elegant formula for $\zeta(2)$, depending upon the quadrature of the circle.\footnote{Weil adds: [i.e., upon $\pi$].}
\end{quote}

      Euler's discovery made him famous, perhaps for the first time, among mathematicians in all Europe. When the news of Euler's discovery reached the city of Basel, the first reaction of his teacher, Johann Bernoulli,\index{Bernoulli, Johann (1667--1748)} was to exclaim that the most burning desire of his deceased older brother Jakob \index{Bernoulli, Jakob (1654--1705)} was now fulfilled.   
Seen all the work he has done on the subject, there is no doubt that throughout his life, Euler tried (without success) to find a formula for $\zeta(s)$ for $s$ an odd
integer.

   It was not unusual for Euler to publish several proofs of the same result, and his result on the convergence on $\zeta(2)$ is one instance of this fact. In particular, there are proofs of this fact in his memoirs \cite{E41} (presented to the Saint Petersburg Academy on December 5, 1735 and published in 1740) and \cite{E72} (presented to the Saint Petersburg Academy on April 25, 1727 and published in 1744), and an account is given in his \emph{Introductio} \cite{Euler-Int-b} (first edition 1748).

              In a letter to Goldbach dated August 28, 1742 (Letter 54 in \cite{Euler-Goldbach}),  Euler expresses $\zeta(2)$ in terms of dilogarithms. 
               We recall that the dilogarithm function\footnote{This name was still not given to that function in the work of Euler mentioned.} is defined as 
               \[\mathrm{Li}(x)=\sum_{k=1}^\infty \frac{x^k}{k^2}.\]
               We have $\mathrm{Li}(1)=\zeta(2)$. In his paper \cite{E20}, presented to the Saint Petersburg Academy in 1731 and published in 1738, Euler had already used the dilogarithm function to find numerical approximations for $\zeta(2)$.

In his memoir \emph{Remarques sur un beau rapport entre les s\'eries des puissances tant directes que r\'eciproques}
(Remarks on a beautiful relation between direct as well as reciprocal power series), \cite{E352}, written in 1749 and published in 1768, Euler\index{Euler, Leonhard (1707--1783)} found the  functional equation satisfied by the zeta function.\index{zeta function}\index{Riemann!zeta function} The relation is not explicitly written by Euler but it follows from a relation he writes, as pointed out by Weil\index{Weil, Andr\'e (1906--1998)} in \cite{Weil-Ens} p. 10, who deduces it immediately from the following formula which Euler\index{Euler, Leonhard (1707--1783)} writes: 
\[\frac{1-2^{n-1}+3^{n-1}-4^{n-1}+5^{n-1}-6^{n-1}+ \ \mathrm{etc}.}{1-2^{-n}+3^{-n}-4^{-n}+5^{-n}-6^{-n}+ \ \mathrm{etc}.}
\]
\[
=\frac{-1.2.3...(n-1)(2^n-1)}{(2^{n-1}-1)\pi^n}\cos \frac{n\pi }{2}
.\]
Weil comments on this formula:
\begin{quote}\small
In the left hand side, we have formally the quotient $\zeta(1-n)\zeta(n)$, except that Euler had written alternating signs to make the series more tractable; the effect of this is merely to multiply $\zeta(n)$ by $1-2^{1-n}$, and $\zeta(1-n)$ by $1-2^n$. In the right hand side we have the gamma function, which Euler\index{Euler, Leonhard (1707--1783)} had invented. Euler proves the formula for every positive integer $n$ (using the so-called  Abel summation to give a meaning to the divergent series in the numerator of the left hand side), and conjectures its validity for all $n$.
\end{quote}

 It was Riemann who showed later on that this equation is valid for any real number $\not=0,1$.

In  his paper \emph{Variae observationes circa series infinitas} which we already mentioned, \cite{E72}, Euler\index{Euler, Leonhard (1707--1783)} found, for $s  >1$, the formula
\[\zeta(s)=\frac{1}{\Gamma(s)}\int_0^\infty \frac{x^{s-1}dx}{e^x-1}.\]
Here $\Gamma$ is the Euler gamma function, which is an extension of the factorial:
\[\Gamma(s)=\int_0^\infty e^{-u}u^{s-1}du.\]
In the same paper, he obtained the following formula, valid for real $s>1$: 
\[\zeta(s)= \prod_p\frac{1}{1-\frac{1}{p^s}}\]
where the product is over all prime numbers $p$. (Weil explains Euler's derivation of this formula in \cite{Weil-Hammurapi} p. 265-266.)  This equality was the starting point of Riemann's investigations in his paper \cite{Riemann-primes}, and it became at the basis of the field called ``analytic number theory." Incidentally, it gives a new proof of the fact that there are infinitely many prime numbers (taking $s=1$ in the formula). We note by the way that Euler\index{Euler, Leonhard (1707--1783)} gave another proof of the existence of infinitely many prime numbers, using the divergence of the harmonic series $\sum \frac{1}{n}$.

After Euler, the next substantial work on the zeta function,\index{zeta function}\index{Riemann!zeta function} $\zeta(s)$, was done more than a century later, by Riemann. Indeed, in the history of number theory that he wrote, Weil\index{Weil, Andr\'e (1906--1998)} considers (see \cite{Weil-Hammurapi} p. 278) that after Euler, the subject was dead, and that Riemann resurrected it. He conjectures that in 1859, Riemann started working on this subject after he seized a remark by Eisenstein,\index{Eisenstein, Ferdinand Gotthold Max (1823--1852)} see  \cite{Weil1989}  for the details. Let us summarize some of the major ideas that Riemann brought in his short paper:

\begin{enumerate}

\item Using analytic continuation, Riemann showed that the zeta function can be extended to a holomorphic function defined on the complex plane, except at the point 1 where the function has a simple pole with residue 1.

\item He discovered the relation between the zeros of the zeta function and the asymptotic distribution of prime numbers. In fact, Riemann gave the principal term in the asymptotic law of the so-called counting function $\pi(x)$ which measures the number of prime numbers $\leq x$. 
More precisely, Riemann gave the formula 
\[\pi(x)\sim\frac{x}{\log x}, \ x\to\infty\]
with a sketch of a proof.
The result became known as the  ``prime number theorem."
 Complete proofs of this theorem
 were given later by Hadamard\index{Hadamard, Jacques (1865--1963)} and de la Vall\'ee Poussin\index{Vall\'ee@de la Vall\'ee Poussin, Charles Jean (1866--1962)} in 1896.

\item Starting from the functional equation discovered by Euler\index{Euler, Leonhard (1707--1783)} -- and of which Riemann provided two new proofs adapted to the newly  extended function -- Riemann showed that the set of zeros of the zeta function contains the even negative integers, and conjectured that all the other zeros are situated on the line $\mathrm{Im}(s)=\frac{1}{2}$. This is the famous Riemann hypothesis.\index{Riemann hypothesis}

\item Riemann obtained a new functional equation satisfied by the zeta function:\index{zeta function}\index{Riemann!zeta function} 

\[
\zeta(s)=2^s\pi^{s-1}\sin (\frac{\pi s}{2})\Gamma(1-s)\zeta(1-s)\] 
for $s\not=0,1$.

\end{enumerate}

               Finding the asymptotic behavior of the prime counting function $\pi(x)$ was, at the epoch of Riemann, one of the major problems in number theory.  Legendre,\index{Legendre, Adrien-Marie (1752--1833)} Gauss and Dirichlet\index{Dirichlet, Johann Peter Gustav Lejeune (1805--1859)} had  already investigated this problem, and more precisely, they worked on a conjecture saying that $\displaystyle \pi(x)$ is asymptotic to a function of the size of $\frac{x}{\ln x}$.   Riemann's  main contribution was the introduction of complex analysis in this study, and his intuition that the distribution of primes is related to the zeros of the zeta function extended to the complex plane. The works by de la Vall\'ee Poussin and Hadamard rely heavily on Riemann's ideas, and the outlines of their proofs are based on his sketch. We talk about Hadamard's work on the zeta function and the prime number theorem in Chapter 8 of the present volume, \cite{Papa-Riemann3}.  Let us add here a few historical notes on the counting function; it will give us the occasion to mention again the work of Legendre.\index{Legendre, Adrien-Marie (1752--1833)}

               In 1798, Legendre\index{Legendre, Adrien-Marie (1752--1833)} published his \emph{Essai sur la th\'eorie des nombres} (Essay on number theory) \cite{Legendre1798}, a  long essay (about 472 pages without the tables) in which, based on numerical evidence, he proposed  a conjecture on the form of the counting function $\pi(x)$. He writes (p. 19):
               \begin{quote}\small
               Moreover, it is likely that the rigorous formula which gives the value of $b$ when $a$ is very large is of the form $\displaystyle b=\frac{a}{A\log a+B}$, $A$ and $B$ being constant coefficients, and $\log a$ denoting a hyperbolic logarithm. The exact determination of these coefficients would be a curious problem, worth of training the expertise of the analysts.\footnote{Au reste, il est vraisemblable que la formule rigoureuse qui donne la valeur de $b$ lorsque $a$ est tr\`es grand, est de la forme $\displaystyle b=\frac{a}{A\log a+B}$, $A$ et $B$ \'etant des coefficients constants, et $\log a$ d\'esignant un logarithme hyperbolique. La d\'etermination exacte de ces coefficients serait un probl\`eme curieux et digne d'exercer la sagacit\'e des Analystes.}
               \end{quote}
              Legendre also gave an approximate value of the constant $A(x)$. 
          Let us note incidentally that Legendre,\index{Legendre, Adrien-Marie (1752--1833)} in the preface to his essay, makes a short history of the development of number theory, starting with the Greeks (Euclid and Diophantus), and passing by Vi\`ete, Bachet, Fermat, Euler and Lagrange.

               In the second edition of his essay (1808), Legendre formulated another conjecture, saying that there are infinitely many primes in any arithmetic progression, that is, primes of the form $l+kn$ for any natural integer $n$.  This conjecture is at the foundations of the theory of Dirichlet series,\index{Dirichlet series} and it was at the basis of several approaches on the prime number theorem. The conjecture was proved by Dirichlet\index{Dirichlet, Johann Peter Gustav Lejeune (1805--1859)} in 1837 \cite{Dirichlet1837}, in a paper which brought new tools on how to approach the prime number theorem. In particular, Dirichlet introduced in this paper his famous $L$-function.\index{Dirichlet $L$-function}\index{L-function@$L$-function}
                              
               Besides Dirichlet and Legendre, one has to mention Gauss,\index{Gauss, Carl Friedrich (1777--1855)} who, at the age of 15 or 16, started an extensive investigation on the distribution of prime numbers. Based mostly on empirical data (tables of prime numbers that he compiled), he observed  that the density of prime numbers around a fixed number $x$ is inversely proportional to $\log x$, and he deduced that the counting function $\pi(x)$ should be well approximated by the integral $\int_2^x \frac{dt}{t}dt$. Gauss never published this work, but he described it in an 1849 letter to his friend and former student, the astronomer J. F. Encke.\index{Encke, Johann Franz (1791--1865)}  Gauss,\index{Gauss, Carl Friedrich (1777--1855)}  in that letter, makes a comparison between his results and those of Legendre. The letter is included in Gauss's correspondence, edited in his Complete Works, and it is also translated and commented in the article \cite{Goldstein} by L. J. Goldstein.

   Finally, one has to mention the work of Chebyshev\index{Chebyshev, Pafnouti Lvovitch (1821--1894)} in his two papers \cite{Cheb-1} and \cite{Cheb-2}, done slightly before Riemann (the papers are published in 1851 and 1852), in which he gave precise approximate values for the prime number counting function, making use of the zeta function in the study of the counting function, as Riemann did in his 1859 paper.  
 Chebyshev's paper \cite{Cheb-2} contains the proof of the so-called Bertrand postulate\index{Bertrand postulate} stating that for any integer $n \geq 3$, there  exists a prime number $p$ satisfying $n<p<2N$.\footnote{The work of Chebyshev deserves to be much more developed than in these few lines. Like his famous Swiss-Russian predecessor Leonhard Euler, Chebyshev published on most of the fields of pure and applied mathematics.  In 1852, he made a stay in France, whose aim was essentially to visit factories and industrial plants, but during his stay he also met several French mathematicians and discussed with them. The list includes  Bienaym\'e, Cauchy, Liouville,  Hermite, Lebesgue, Poulignac, Serret and others. A detailed report on this stay, written by Chebyshev himself, is contained in his \emph{Collected works} \cite{T-oeuvres}. Chebyshev used to  published in French journals and his relations with French mathematicians remained constant over the years. In 1860, he was elected  corresponding member of the Paris Academy of Sciences, and in 1874 foreign member. We learn from his report that at the end of his 1852 stay in France, on his way back to Russia, Chebyshev stopped in Berlin and had several discussions with Dirichlet.\index{Dirichlet, Johann Peter Gustav Lejeune (1805--1859)} It is conceivable that during that meeting the two mathematicians talked about the problems related to  the prime number counting function. We refer the reader to the article \cite{2016-Tchebyshev} where some of Chebyshev's works are compared with works of Euler.}

The question of the zeros of the zeta function\index{zeta function}\index{Riemann!zeta function} was proposed by Hilbert in one of the problems he offered at the Paris 1900 ICM.

 Riemann's memoir \cite{Riemann-primes} had a major influence on several later mathematicians, including Weil,\index{Weil, Andr\'e (1906--1998)} Siegel,\index{Siegel, Carl
Ludwig (1896--1981)} and Selberg.

 We conclude this section by quoting Weil,\index{Weil, Andr\'e (1906--1998)} from an obituary article by A. Knapp \cite{Knapp}:

 \begin{quote} \small
A substantial portion of Weil's research was motivated by an effort to prove the Riemann hypothesis concerning the zeroes of the Riemann zeta function. He was continually looking for new ideas from other fields that he could bring to bear on a proof. He commented on this matter in a 1979 interview:\footnote{Pour la Science, November 1979.} Asked what theorem he most wished he had proved, he responded, ``In the past it sometimes occurred to me that if I could prove the Riemann hypothesis, which was formulated in 1859, I would keep it secret in order to be able to reveal it only on the occasion of its centenary in 1959. Since in 1959, I have felt that I am quite far from it, I have gradually given up, not without regret."\footnote{Autrefois, il m'est quelquefois venu \`a l'esprit que, si je pouvais d\'emontrer l'hypoth\`ese de Riemann, laquelle avait \'et\'e formul\'ee en 1859, je la garderais secr\`ete pour ne la r\'ev\'eler qu'\`a l'occasion de son centenaire en 1959. Comme en 1959, je m'en sentais encore bien loin, j'y ai peu \`a peu renonc\'e, non sans regret.}
\end{quote}
One of the famous Weil\index{Weil, Andr\'e (1906--1998)} conjectures is known as the ``Riemann hypothesis over finite fields."

\section{On space}\label{s:space}
 
Riemann's habilitation lecture contains a discussion on the nature of physical space and its relation with geometry. The concepts on which Riemann dwells there make it clear that the theme of space belongs to his profound thought. One of the main ideas on which he stresses is the possibility that physical space is different from the space of Euclidean geometry, a point of view that makes Riemann in some sense a predecessor of modern physics.

In speculating on space, Riemann follows a long tradition which includes the Greeks,  Newton, Descartes, Kant and many others, a tradition which survived until the modern period; one may mention, among the mathematicians of the post-Riemannian period, Hermann Weyl,\index{Weyl, Hermann (1885--1955)} Ren\'e Thom,\index{Thom, Ren\'e (1923--2002)} Alexandre Grothendieck,\index{Grothendieck, Alexandre (1928--2014)} and there are many others. It is therefore natural to have, in this paper, a section on space, in which, not only we review Riemann's ideas -- this is done in several chapters of the present volume-- but where we mention some of the ideas on this subject that were expressed by his predecessors. Our exposition will necessarily be succinct. Writing a serious essay on the notion of space needs a whole volume.

Space\index{space} is one of the first very few basic philosophico-epistemological notions. It appears at several places in the works of Aristotle:\index{Aristotle (384--322 B.C.)} there are sections on space in the \emph{Categories},\index{Categories (Aristotle)}
\cite{Aristotle-Categories}, the \emph{Physics} \cite{Aristotle-Physics},  the  \emph{Metaphysics} \cite{Aristotle-Metaphysics}, the treatise \emph{On the heavens} \cite{Aristotle-Heavens}, etc. Furthermore, like for many other subjects, we learn from Aristotle's works the opinions of his predecessors on space: the Meletians, the Pythagoreans, Plato,\index{Plato (5th--4th c. B.C.)} etc.

 In the \emph{Categories}\index{Categories (Aristotle)} (5a, 8-14), Aristotle explains that space, like time, belongs to the category of \emph{continuous quantity}.\footnote{In the \emph{Categories}, (4b 20-5b 11) Aristotle distinguishes seven different types of quantities, which he classifies as continuous and discrete. Discrete quantity comprises number and speech.  Continuous quantity comprises the line, the surface, the body, time, and space. Needless to say, although this classification may appear limited from a modern point of view, it has the great merit of existing, may be for the first time. Aristotle asked the pertinent questions.} In Book IV of his \emph{Physics}, he writes about the difference between ``space" and ``place."  This is a fundamental distinction, with an impact in physics, and it had a huge influence on later thinkers.\footnote{This theme of space and its relation to place was particularly expanded by Aristotle's commentators. We mention in particular the medieval Andalusian polymath Averroes (1126--1198).\index{Averroes (Ibn Rushd of Cordoba) (1126--1198)} The third chapter of Rashed's book \emph{Les math\'ematiques
infinit\'esimales du IX\`eme au XI\`eme si\`ecle} \cite{Rashed4} contains a critical edition together with a translation and commentaries of the treatise \emph{On space}  by the Arabic scientist Ibn al-Haytham\index{Ibn al-Haytham,  Ab\=u `Al\=\i \ al-\d{H}asan ibn al-\d{H}asan (965--1039)} (known in the West under the name al-Hazen) in which this author criticizes Aristotle's theory of space developed in his \emph{Physics}, and where he defines subsets
of space by metric properties. 
There is also a rich discussion onthe notion of space in Greek philosophy in the multi-volume encyclopedic work of  P. Duhem \cite{Duhem}, see in particular vol. I, p. 197ff.} The question has also implications in the history of topology. The Greek origin for the word place is \emph{topos}  ($\tau \acute{o} \pi o \varsigma$), and  is translated into Latin by \emph{situs}. The expression \emph{analysis situs}\index{analysis situs@\emph{analysis situs}}, which was used by Leibniz and the Western founders of topology,  finds its origin there.

 Among the Western thinkers whose work on the theme of space emerges amid the classical philosophical monuments, we mention Galileo, Newton, Descartes, Leibniz, Huygens\index{Huygens, Christiaan (1629--1695)} and Kant.\index{Kant, Immanuel (1724--1804)} Most of them are quoted by Riemann.

We start by quoting a text from Greek antiquity. This is a fragment by Archytas of Tarentum\index{Archytas of Tarentum (428?347 B.C.)} which is often referred to in the literature on Pythagorean philosophy, to show the kind of questions on space and on place that the ancient Greeks addressed, e.g., whether space is bounded or not, and the paradoxes to which this question leads (see \cite{Huffman} p. 541):
\begin{quote}\small
``But Archytas," as Eudemus says, ``used to propound the argument
in this way: `If I arrived at the outermost edge of the
heaven [that is to say at the fixed heaven], could I extend my
hand or staff into what is outside or not?' It would be paradoxical
not to be able to extend it. But if I extend it, what is outside
will be either body or place. It doesn't matter which, as we will
learn. So then he will always go forward in the same fashion to
the limit that is supposed in each case and will ask the same
question, and if there will always be something else to which
his staff [extends], it is clear that it is also unlimited. And if it is
a body, what was proposed has been demonstrated. If it is
place, place is that in which body is or could be, but what is
potential must be regarded as really existing in the case of eternal things, and thus there would be unlimited body and space."
(Eudemus, Fr. 65 Wehrli, Simplicius, In Ar. Phys. iii 4; 541)
\end{quote}

The most basic question that was addressed by many of the philosophers of the modern period that we mentioned is probably the following: Does space have an objective existence or is it only a construction of human mind? Before trying to answer this question, or to have an opinion on it, it is helpful to make it precise what notion of space it refers to:  three-dimensional physical space? the three-dimensional space of Euclidean geometry? an abstract notion of space? Other related questions are: Is Euclid's three-dimensional geometry a pure logical construction or is it a mathematical formulation of the properties of external nature? Is the space of (theoretical) physics the same as the mathematicians' space? Does void exist, and what function does it have? These are some of the questions which obviously obsessed Riemann, and before him, many others.

In Descartes' doctrine, space depends on matter, therefore void cannot exist. Leibniz\index{Leibniz, Gottfried Wilhelm (1646--1716)} and Euler after him shared the same opinion. Newton had a notion of ``absolute space" and ``relative space." Furthermore, following the ancient Greeks,  Descartes made  a difference between space and place.  We quote some passages fom his  \emph{Principes de la philosophie} (Principles of philosophy) \cite{Descartes-Principes} (1644).
\begin{quote}\small
\emph{Principle XIV. How place and space differ:} However, place and space are different in names, because place indicates more expressly situation than magnitude or figure, and that on the contrary, we think about that one when we talk about space; for we say that a  thing entered at the place of another, even though it does not have exactly neither the same magnitude nor figure, and for that we do not mean that it occupies the same space that this other thing occupies; and when the situation is changed, we say that the place has also changed, even though it has the same magnitude and figure than before: in this sort, if we say that a thing is in some place, we only mean that it is situated in such a way with respect to other things; but if we add that it occupies a certain space, or place, then we mean  that it has such magnitude and figure that it can occupy it exactly.\footnote{\emph{Principe XIV. Quelle difference il y a entre le lieu et l'espace :}
Toutefois le lieu et l'espace sont diff\'erents en 
 leurs noms, parce que le lieu nous marque plus 
 express\'ement la situation que la grandeur ou la 
figure, et qu'au contraire nous pensons plut\^ot \`a 
celles-ci lorsqu'on nous parle de l'espace ; car nous 
disons qu'une chose est entr\'ee en la place d'une 
autre, bien qu'elle n'en ait exactement ni la grandeur ni la figure, et n'entendons point qu'elle occupe pour cela le m\^eme espace qu'occupait cette autre chose ; et lorsque la situation est chang\'ee, 
nous disons que le lieu est aussi chang\'e, quoiqu'il soit de m\^eme grandeur et de m\^eme figure 
qu'auparavant : de sorte que si nous disons qu'une 
chose est en un tel lieu, nous entendons seulement qu'elle est situ\'ee de telle façon \`a l'\'egard de 
quelques autres choses ; mais si nous ajoutons 
qu'elle occupe un tel espace, ou un tel lieu, nous 
entendons outre cela qu'elle est de telle grandeur et de telle figure qu'elle peut le remplir tout 
justement.}
\end{quote}
\begin{quote}\small 
\emph{Principle XV: How the surface surrounding a body can be taken as its exterior place:} Thus, we never make a distinction between space and extent, for what regards length, width and depth; but we sometimes consider place as if it were within the thing which is placed, and sometimes also as if it were outside it. By no means the interior differs from space; but sometimes we take the exterior to be either the surface surrounding immediately the thing which is placed (and one has to notice that by surface we must not intend any part of the body surrounding it but only the extremity which is between the body which surrounds and the one which is surrounded which is only a mode or a way), or to be the surface in general, which is not part of a body rather than another one, and which always seems to be the same, provided it has the same magnitude and the same figure; because even if we see that the body that surrounds another body passes somewhere else with its surface, we are not used to say that what was surrounded by it has changed its place for this reason, it stays at the same situation regarding the other bodies that we consider as still. Thus, we say that a boat which is carried away by the stream of a river, and which is at the same time pushed away by the wind by a force which is so equal that it does not change its situation regarding the shores, stays at the same place, even though we see that all the surface that surrounds it changes permanently.\footnote{\emph{Principe XV. Comment la superficie qui environne un corps peut \^etre prise pour son lieu 
exterieur :}
Ainsi nous ne distinguons jamais l'espace d'avec l'\'etendue 
en longueur, largeur et profondeur ; mais nous consid\'erons 
quelquefois le lieu comme s'il \'etait en la chose qui est plac\'ee, et quelquefois aussi comme s'il en \'etait dehors. L'int\'erieur ne diff\`ere en aucune façon de l'espace ; mais nous 
prenons quelquefois l'ext\'erieur ou pour la superficie qui environne imm\'ediatement la chose qui est plac\'ee (et il est \`a remarquer que par la superficie on ne doit entendre aucune 
partie du corps qui environne, mais seulement l'extr\'emit\'e 
qui est entre le corps qui environne et celui qui est environn\'e, qui n'est rien qu'un mode ou une fa\c con), ou bien pour la superficie en g\'en\'eral, qui n'est point partie d'un 
corps plut\^ot que d'un autre, et qui semble toujours la m\^eme, 
tant qu'elle est de m\^eme grandeur et de m\^eme figure ; car
encore que nous voyions que le corps qui environne un autre 
corps passe ailleurs avec sa superficie, nous n'avons pas coutume de dire que celui qui en \'etait environn\'e ait pour cela chang\'e de place lorsqu'il demeure en la m\^eme situation \`a 
l'\'egard des autres corps que nous consid\'erons comme immobiles. Ainsi nous disons qu'un bateau qui est emport\'e par 
le cours d'une rivi\`ere, et qui en m\^eme temps est repouss\'e 
par le vent d'une force si \'egale qu'il ne change point de situation \`a l'\'egard des rivages, demeure en m\^eme lieu, bien 
que nous voyions que toute la superficie qui l'environne 
change incessamment.}
\end{quote}

 Euler had also a strong philosophical background and, needless to say, a tendency for abstraction. We recall that the subject of his first public lecture, delivered  at the University of Basel at the occasion of his graduation, was the comparison between the philosophical systems of Newton and Descartes. The notions of space, of motion and of force are discussed in several of his papers on physics. His\index{Euler, Leonhard (1707--1783)} most important work related to these matters is his \emph{Mechanica}, in two volumes of 500 pages each, \cite{E15}  with its systematic use of analysis (differential equations) in the field of mechanics, as opposed to Newton's geometric point of view developed in his \emph{Principia}.  In his memoir  \emph{Recherches sur l'origine des forces}
(Research on the origin of forces) \cite{E121} (1750), Euler uses an argument involving a notion of ``impenetrability of bodies" from which he deduces the law of shock of bodies.  
We also mention his \emph{Anleitung zur Naturlehre, worin die Gründe zu Erklärung aller in der Natur sich ereignenden Begebenheiten und Veränderungen festgesetzet wedren} (Introduction to natural science establishing the fundamentals for the
explanation of the events and changes that occur in nature), \cite{E842}, a long memoir written in 1745, but never completed and published in 1862. Hermann Weyl\index{Weyl, Hermann (1885--1955)} says (\cite{Weyl-philo}  p. 42) about this memoir that Euler ``in magnificent clarity summarizes the foundations of the philosophy of  nature of his time." 
In this memoir, Euler discusses notions like the \emph{extent} of material bodies, the infinite divisibility of these bodies, motion, space, place magnitude, aether and  gravity. 
His memoir   \emph{Essai d'une d\'emonstration m\'etaphysique du principe g\'en\'eral de l'\'equilibre} (Essay on a metaphysical demonstration of the
general principle of equilibrium) \cite{E200} concerns  again,  force, equilibrium, motion and gravity.
In his memoir \emph{R\'eflexions sur l'espace et le temps}
(Reflections on space and time)  \cite{E149}, he makes a comparison between the mathematicians' and the philosophers' (which he calls the ``metaphysicians") points of view. He describes \emph{position} as the relation of a body
with other bodies around it. He declares that the metaphysicians are wrong in claiming that the notions of space and place are abstract constructions of the mind, and he argues to show  the reality of space and time. He claims that both absolute space and time, as mathematicians
represent them, are real and exist beyond human imagination.  He discusses inertia and the relativity of motion, the ideas of place and position, supported by notions from mechanics.
 
 Euler's philosophical ideas, and their impact on Riemann, have not yet been seriously discussed in the literature.

Immanuel Kant is among the commanding figures that preceded Riemann on the subject of philosophy of space. As a matter of fact, space was already a major theme in Kant's\index{Kant, Immanuel (1724--1804)}  \emph{Inaugural dissertation} (1770).\index{doctoral dissertation!Kant}\index{Kant!doctoral dissertation} Kant expresses there his doctrine of the a priori nature of space and of geometric objects, that is, the belief that they are not derived from an outside experience. The following excerpt contains an expression of this point of view, which, as we shall recall,  Gauss criticized later (\cite{Kant-Dissertation} \S\,15, A--D):
 \begin{quote}\small
 The concept of space is not abstracted from external sensations. For I am unable to conceive of anything posited without me unless by representing it as in a place different from that in which I am, and of things as mutually outside of each other unless by locating them in different places in space. Therefore the possibility of external perceptions, as such, presupposes and does not create the concept of space, so that, although what is in space affects the senses, space cannot itself be derived from the senses.

  The concept of space is a singular representation comprehending all things in itself, not an abstract and common notion containing them under itself. What are called several spaces are only parts of the same immense space mutually related by certain positions, nor can you conceive of a cubic foot except as being bounded in all directions by surrounding space.

 The concept of space, therefore, is a pure intuition, being a singular concept, not made up by sensations, but itself the fundamental form of all external sensation. This pure intuition is in fact easily perceived in geometrical axioms, and any mental construction of postulates or even problems. That in space there are no more than three dimensions, that between two points there is but one straight line, that in a plane surface from a given point with a given right line a circle is describable, are not conclusions from some universal notion of space, but only discernible in space as in the concrete. Which things in a given space lie toward one side and which are turned toward the other can by no acuteness of reasoning be described discursively or reduced to intellectual marks. There being in perfectly similar and equal but incongruous solids, such as the right and the left hand, conceived of solely as to extent, or spherical triangles in opposite hemispheres, a difference rendering impossible the coincidence of their limits of extension, although for all that can be stated in marks intelligible to the mind by speech they are interchangeable, it is patent that only by pure intuition can the difference, namely, incongruity, be noticed. Geometry, therefore, uses principles not only undoubted and discursive but falling under the mental view, and the obviousness of its demonstrations -- which means the clearness of certain cognition in as far as assimilated to sensual knowledge -- is not only greatest, but the only one which is given in the pure sciences, and the exemplar and medium of all obviousness in the others. For, since geometry considers the relations of space, the concept of which contains the very form of all sensual intuition, nothing that is perceived by the external sense can be clear and perspicuous unless by means of that intuition which it is the business of geometry to contemplate. Besides, this science does not demonstrate its universal propositions by thinking the object through the universal concept, as is done in intellectual disquisition, but by submitting it to the eyes in a single intuition, as is done in matters of sense.

  Space is not something objective and real, neither substance, nor accident, nor relation; but subjective and ideal, arising by fixed law from the nature of the mind like an outline for the mutual co-ordination of all external sensations whatsoever. Those who defend the reality of space either conceive of it as an absolute and immense receptacle of possible things, an opinion which, besides the English, pleases most geometricians, or they contend for its being the relation of existing things itself, which clearly vanishes in the removal of things and is thinkable only in actual things, as besides Leibniz, is maintained by most of our countrymen. The first inane fiction of the reason, imagining true infinite relation without any mutually related things, pertains to the world of fable. But the adherents of the second opinion fall into a much worse error. Whilst the former only cast an obstacle in the way of some rational or monumental concepts, otherwise most recondite, such as questions concerning the spiritual world, omnipresence, etc., the latter place themselves in fiat opposition to the very phenomena, and to the most faithful interpreter of all phenomena, to geometry. For, not to enlarge upon the obvious circle in which they become involved in defining space, they cast forth geometry, thrown down from the pinnacle of certitude, into the number of those sciences whose principles are empirical. If we have obtained all the properties of space by experience from external relations only, geometrical axioms have only comparative universality, such as is acquired by induction. They have universality evident as far as observed, but neither necessity, except as far as the laws of nature may be established, nor precision, except what is arbitrarily made. There is hope, as in empirical sciences, that a space may some time be discovered endowed with other primary properties, perchance even a rectilinear figure of two lines.
\end{quote}
  
The reader will notice that Kant talks about ``geometrical axioms," and mentions axioms of Euclidean geometry such as the fact that ``between two points there is but one straight line." Kant was by no means a mathematician, but he had a sufficient knowledge, as a philosopher, of several basic principles of mathematics.

It appears from Gauss's correspondence, published in Volume VII of his \emph{Collected Works} (p. 200ff.) that Gauss he meditating on the nature of space since a very young age, probably from the age of 16. It is from these meditations that he became interested in the parallel postulate and in non-Euclidean geometry, spherical and (the hypothetical) hyperbolic. Unlike most of the geometers that preceded him, Gauss was convinced, at a very early stage of his life, that the parallel postulate was not a consequence of the others, and he spent a lot of time and energy pondering on the principles of hyperbolic geometry, a geometry resulting from the negation of the postulate.

Gauss was also thoroughly interested in philosophy, and, in particular, he read Kant. 
He became very critical of the latter's conception of space, exemplified in the text we just quoted as being ``not something objective and real, neither substance, nor accident, nor relation, but subjective and ideal, arising by fixed law from the nature of the mind." On Kant,\index{Kant, Immanuel (1724--1804)} Gauss had the advantage of being a mathematician. In a letter to his friend Bessel,\index{Bessel, Friedrich (1784--1846)} dated April 9, 1830, Gauss writes (translation from \cite{Carus} p. 13):
\begin{quote}\small
We must confess in all humility that a number is \emph{solely} a product of our mind. Space, on the other hand, possesses also a reality outside of our minds, the laws of which we cannot fully prescribe \emph{a priori}.
\end{quote}

In another letter, sent to Wolfgang Bolyai\index{Bolyai, Wolfgang (1775--1856)} on March 6, 1832 and published in his \emph{Collected Works}, Gauss writes, concerning the two hypotheses on the angle sum in a triangle, that it is precisely in the difficulty of this decision that ``lies the clearest proof that Kant was wrong in asserting that space is just a form of our perception."
 
  Gauss was also very critical of Kant's argument based on symmetries in the text we quoted above (``There being in perfectly similar and equal but incongruous solids, such as the right and the left hand, conceived of solely as to extent,\ldots it is patent that only by pure intuition can the difference, namely, incongruity, be noticed"). We further discuss this in Chapter 6 of the present volume \cite{Papa-Riemann1}. 
 
  It is not surprising that Riemann declares, in his habilitation lecture\index{Riemann! habilitation lecture},\index{habilitation lecture!Riemann} that, concerning his ideas on space,  he is influenced by Gauss.
  
Riemann's ideas on space were discussed by Clifford, the first mathematician who translated into English Riemann's habilitation text, cf. \cite{Clifford-space}.

\section{Topology}\label{s:topo}

 Poincar\'e, who is certainly the major founder of the modern field of topology,\footnote{We may quote P. S. Alexandrov, who declared in a talk he gave at a celebration of the centenary of Poincar\'e's birth \cite{Alex}: ``To the question of what is Poincar\'e's relationship to topology, one can reply in a single sentence: he created it." On Poincar\'e and Riemann, Alexandrov, in the same talk, says the following: ``The close connection of the theory of functions of a complex variable, which Riemann has observed in embryonic form, was first understood in all its depth by Poincar\'e."} declares in his ``Analysis of his own works" (\emph{Analyse des travaux scientifiques de Henri Poincar\'e faite par lui-m\^eme}), \cite{Poin-Acta} p. 100,  that he has two predecessors in the field, namely, Riemann and Betti. The latter,\index{Betti, Enrico (1923--1892)} in his correspondence with his friend and colleague Placido Tardy\index{Tardy, Placido (1816--1914)} reports on several conversations he had with Riemann on topology. Two letters from Betti to Tardy on this subject are reproduced and translated in the book \cite{Pont1974} by Pont, in the article \cite{Weil1979a} by Weil, and prior to them, by Loria in his obituary on Tardy \cite{Loria-Tardy}.
 
The first of these two letters by Betti, dated October 6, 1863, starts with the following (Weil's translation): ``I have newly talked with Riemann about the connectivity of spaces,
 and have formed an accurate idea of the matter," and he goes on explaining to his friend the notion of connectivity and that of the order of connectivity. Betti then writes: 
 \begin{quote}\small 
 What gave Riemann the idea of the
 cuts was that Gauss defined them to him, talking about other matters, in a private conversation. In his writings one finds that analysis situs, that is, this
 consideration of quantities independently from their measure, is "wichtig"; in the last years of his life he has been much concerned with a problem in analysis situs, namely: given a winding thread and knowing, at every one of its self-intersections, which part is above and which below, to find whether it can be unwound without making knots; this problem he did not succeed in solving except in special cases ...
 \end{quote} 
   The second letter, dated October 16, 1863, starts with: ``Riemann proves quite easily that every space can be reduced to an SC space by means of 1-cuts and SC 2-cuts." In the same letter, Betti\index{Betti, Enrico (1923--1892)} elaborates on this subject, giving many examples in $n$ dimensions. He concludes the letter by noting that the number of line sections is equal to the number of periodicity moduli of an $(n-1)$-integral, the number of simply connected surface sections to the number of periodicity moduli of an $(n-2)$-integral,  and so on. 
 
 This should make clear the parentage, for what concerns topology, from Riemann to Poincar\'e, potentially including Betti. 
In this section, we go further back in the history of topological ideas, and we review some of the important works done before Riemann in this field. 

Ren\'e Thom\index{Thom, Ren\'e (1923--2002)} considers that topology was born in ancient Greece. He expanded on this idea in several articles, cf. \cite{Thom2} and \cite{Thom1}. This is a  perfectly reasonable theory. In fact, the question depends on what sense we give to the word ``topology."  If the matter concerns the notions of limit and convergence,  then the roots of this field are indeed in Greek antiquity, and more especially, in the writings of Zeno,\index{Zeno of Elea (5th c. B.C.)}  which do not survive, but which were quoted by his critics and commentators, including Plato,\index{Plato (5th--4th c. B.C.)} Aristotle\index{Aristotle (384--322 B.C.)} and Simplicius.\index{Simplicius (6th c. A.D.)} Likewise, if the question concerns the notion of space, and the related notion of place, then  the roots also are in Greek science. We already alluded to this fact in the previous section. The Greeks made a distinction between space and place and the notion of place (situs) is at the basis of topology. The three words place, situs and $\tau \acute{o} \pi o \varsigma$ are synonyms. To the best of our knowledge, a systematic investigation of the origin of topology in Greek antiquity has  never been conducted. A whole book may be written on that subject. Failing to do this now, we shall start our exposition of the roots of topology with Leibniz, as it is usually done. Indeed, it is commonly accepted that the first explicit mention of topology as a mathematical field was made by him.\index{Leibniz, Gottfried Wilhelm (1646--1716)} 

Even though no purely topological result can be attributed to Leibniz, he had the privilege to express for the first time, back in the seventeenth century, the need for a new branch of mathematics, which would be ``a geometry that is more general than the rigid Euclidean geometry and the analytic geometry of Descartes."\index{Descartes, Ren\'e (1596--1650)} Leibniz describes his geometry as purely qualitative and concerned with the study of figures independently of their metrical properties. 
In a letter to Christiaan Huygens, sent on September 8, 1679 (cf. \cite{Leibniz-corresp} p.~578--569 and \cite{Huygens} vol.  VIII n${}^{\mathrm{o}}$ 2192), he writes:
\begin{quote}\small
After all the progress I have made in these matters, I am still not happy with Algebra, because it provides neither the shortest ways nor the most beautiful constructions of Geometry. This is why when it comes to that, I think that we need another analysis which is properly geometric or linear, which expresses to us directly \emph{situm}, in the same way as algebra expresses \emph{magnitudinem}. And I think that I have the tools for that, and that we might represent figures and even engines and motion in character, in the same way as algebra represents numbers in magnitude.\footnote{Apr\`es tous les progr\`es que j'ai faits en ces mati\`eres, je ne suis pas encore content de l'Alg\`ebre, en ce qu'elle ne donne ni les plus courtes voies, ni les plus belles constructions de G\'eom\'etrie. C'est pourquoi lorsqu'il s'agit de cela, je crois qu'il nous faut encore une autre analyse proprement g\'eom\'etrique ou lin\'eaire, qui nous exprime directement \emph{situm}, comme l'alg\`ebre exprime \emph{magnitudinem}. Et je crois d'en avoir le moyen, et qu'on pourrait repr\'esenter des figures et m\^eme des machines et mouvements en caract\`eres, comme l'alg\`ebre repr\'esente les nombres en grandeurs. [We have modernized the French.]}
\end{quote}
In the same letter (\cite{Leibniz-corresp} p. 570), Leibniz adds:
 
\begin{quote}\small
I found the elements of a new characteristic,\index{characteristic (Leibniz)} completely different from Algebra and which will have great advantages for the exact and natural mental representation,  although without figures, of everything that depends on  the imagination.
Algebra is nothing but the characteristic of undetermined numbers or magnitudes. But it does not directly express the place, angles and motions, from which it follows that it is often difficult to reduce, in a computation, what is in a figure, and that it is even more difficult to find geometrical proofs and constructions which are enough practical even when the Algebraic calculus is all done.\footnote{J'ai trouv\'e quelques \'el\'ements d'une nouvelle caract\'eristique, tout \`a fait diff\'erente de l'Alg\`ebre, et qui aura de grands avantages pour repr\'esenter \`a l'esprit exactement et au naturel, quoique sans figures, tout ce qui d\'epend de l'imagination. L'Alg\`ebre n'est autre chose que la caract\'eristique des nombres ind\'etermin\'es ou des grandeurs. Mais elle n'exprime pas directement la situation, les angles et les mouvements, d'o\`u vient qu'il est souvent difficile de r\'eduire dans un calcul ce qui est dans la figure, et qu'il est encore plus difficile de trouver des d\'emonstrations et des constructions g\'eom\'etriques assez commodes lors m\^eme que le calcul d'Alg\`ebre est tout fait.}
\end{quote}

Together with his letter to Huygens, Leibniz included the manuscript of an essay he wrote on the new subject. He writes, in the same letter (\cite{Leibniz-corresp} p. 571):
\begin{quote}\small
But since I don't see that anybody else has ever had the same thought, which makes me fear that it might be lost if I do not get enough time to complete it, I will add here an essay which seems to me important, and which will suffice at least to rendre my aim more credible and easier to conceive, so that if something  prevents its realization now, it will serve as a monument for posterity and give the possibility to somebody else to finish it.\footnote{Mais comme je ne remarque pas que quelqu'autre ait jamais eu la m\^eme pens\'ee, ce qui me fait craindre qu'elle ne se perde, si je n'y ai pas le temps de l'achever, j'ajouterai ici un essai qui me para\^\i t consid\'erable, et qui suffira au moins \`a rendre mon dessein plus croyable et plus ais\'e \`a concevoir, afin que si quelque hasard en emp\^eche la perfection \`a pr\'esent, ceci serve de monument \`a la post\'erit\'e, et donne lieu \`a quelque autre d'en venir \`a bout.}
\end{quote}

He then explains in more detail his vision of this new domain of mathematics, and where it stands with respect to algebra and geometry, giving several examples of a formalism to denote loci, showing how this formalism expresses statements such that the intersection of two spherical surfaces is a circle, and the intersection of two planes is a line.

 Leibniz' letter\index{Leibniz, Gottfried Wilhelm (1646--1716)} ends with the words 
  (\cite{Leibniz} p. 25): 
\begin{quote}\small
I have only one remark to add, namely, that I see that it is possible to extend the characteristic\index{characteristic (Leibniz)} to things which are not subject to imagination. But this is too important and it would lead us too far for me to be able to explain myself on that in a few words.\footnote{Je n'ai qu'une remarque \`a ajouter, c'est que je vois qu'il est possible d'\'etendre la caract\'eristique jusqu'aux choses, qui ne sont pas sujettes \`a l'imagination ; mais cela est trop important et va trop loin pour que je me puisse expliquer l\`a-dessus en peu de paroles.}
\end{quote}

When Leibniz started his correspondence with Huygens,  the latter was already a well established scientist whose achievements were behind him, and it was not easy to convince him of the usefulness of a new theory. Huygens thought that the theory was too abstract and he remained skeptical about it.
He was above all a geometer working on concrete geometrical  problems.

One may recall that when Leibniz sent him the above letter,  Huygens was considered as a world authority in geometry and physics. He was settled in Paris since 15 years, and he was a leading member of the \emph{Acad\'emie Royale des Sciences}. 
  Leibniz had studied mathematics with Huygens, who was seventeen years older than him, and he considered him as his mentor.  Huygens responded to Leibniz  in a letter dated November  22, 1679 (\cite{Leibniz-corresp} p. 577):
\begin{quote}\small
I have examined carefully what you are asking me regarding your new characteristic,\index{characteristic (Leibniz)} and to be frank with you, I cannot not conceive the fact that you have so much expectations from what you spread on me. Because your example of places concerns only realities which were already perfectly known, and the proposition saying that the intersection of a plane and a spherical surface makes the circumference of a circle does not follow clearly. Finally, I cannot see in what way you can apply your characteristic to which you seem you want to reduce all these different matters, like the quadratures, the invention of curves by the properties of tangents, the irrational roots of equations, Diophantus' problems, the shortest and the most beautiful constructions of the geometric problems. And what still appears to me stranger than anything else, the invention and the explanation of machines. I say it to you unsuspiciously, in my opinion this is only wishful thinking, and I need other proofs in order to believe that there could be some reality in what you present. I would nevertheless restrain myself from saying that you are mistaken, knowing the subtlety and the deepness of your mind. I only beg you that the magnificence of the things you are searching won't let you postpone from giving us those which you already found, like this Arithmetic Quadrature you discovered, concerning the roots of the equations beyond the cubical, if you are still satisfied with it.\footnote{J'ai examin\'e attentivement ce que vous me demandez touchant votre nouvelle caract\'eristique, mais pour vous l'avouer franchement, je ne con\c cois pas parce que vous m'en \'etalez, que vous y puissiez fonder de si grandes esp\'erances. Car votre exemple des Lieux ne regarde que des v\'erit\'es qui nous \'etaient d\'ej\`a fort connues, et la proposition de ce que l'intersection d'un plan et d'une surface sph\'erique fait la circonf\'erence d'un cercle, s'y conclut assez obscur\'ement. Enfin, je ne vois point de quel biais vous pourriez appliquer votre caract\'eristique \`a toutes ces choses diff\'erentes qu'il semble que vous y vouliez r\'eduire, comme les quadratures, l'invention des courbes par la propri\'et\'es des tangentes, les racines irrationnelles des \'Equations, les probl\`emes de Diophante, les plus courtes et plus belles constructions des probl\`emes g\'eom\'etriques. Et ce qui me para\^\i t encore le plus \'etrange, l'invention et l'explication des machines. Je vous le dis ing\'enument, ce ne sont l\`a \`a mon avis que de beaux souhaits, et il me faudrait d'autres preuves pour croire qu'il y e\^ut de la r\'ealit\'e dans ce que vous avancez. Je n'ai pourtant garde de dire que vous vous abusiez, connaissant d'ailleurs la subtilit\'e et profondeur de votre esprit. Je vous prie seulement que la grandeur des choses que vous cherchez ne vous fasse point diff\'erer de nous donner celles que vous avez d\'ej\`a trouv\'ees, comme est cette Quadrature Arithm\'etique et que vous avez d\'ecouvert pour les racines des \'equations au-del\`a du cube, si vous en \^etes content vous-m\^eme.}  \end{quote}

In another letter dated January 11, 1680 (\cite{Leibniz-corresp} p. 584) Huygens writes:

\begin{quote}\small
For what concerns the effects of your characteristic,\index{characteristic (Leibniz)} I see that you insist on being persuaded of them, but as you say yourself, the examples will be more important than reasonings. This is why I am asking you much simpler examples, but capable of overcoming my incredulity, because that of the places, I confess, does not seem to me of that sort.\footnote{Pour ce qui est des effets de votre caract\'eristique, je vois que vous persistez \`a en \^etre persuad\'e, mais, comme vous dites vous-m\^eme, les exemples toucheront plus que les raisonnements. C'est pourquoi je vous en demande des plus simples, mais propres \`a convaincre mon incr\'edulit\'e, car celui des lieux, je l'avoue, ne me para\^\i t pas de cette sorte.}
\end{quote}

The essay that Leibniz sent did not obtain Huygens' backing  and it remained hidden among other manuscripts in Huygens' estate. It was published for the first time in 1833, and drew the attention of several nineteenth-century mathematicians, including Grassmann (1809--1877), the founder of the theory of vector spaces, who realized its importance for the new field of topology. There are two recent editions of this text, both included  in doctoral dissertations, by J. Acheverr\'\i a  \cite{Leibniz-Characteristica} (1995),  in France,  and by de Risi, \cite{Risi} (2007), in Germany. The two dissertations contain other texts by Leibniz on the same subject.

Leibniz\index{Leibniz, Gottfried Wilhelm (1646--1716)} used several names for the new field, including\index{analysis situs@\emph{analysis situs}} \emph{analysis situs, geometria situs, characteristica situs, characteristica geometrica, analysis geometrica, speciosa situs}, etc.

The first mathematician who worked consciously on topological questions is Euler. These questions include the definition and the invariance of the Euler characteristic of a convex polyheron,\index{Euler characteristic}\index{theorem!Euler characteristic} the problem known as that of the K\"onigsberg seven bridges, another question related to the Knight's tour on the chessboard, and a musical question concerning a graph known as the \emph{speculum musicum}\index{speculum@\emph{speculum musicum}}. This graph was introduced in Euler's \emph{Tentamen novae theoriae musicae ex certissimis harmoniae
principiis dilucide expositae} (A attempt at a new theory of music, exposed in all clearness
according to the most well-founded principles of harmony) \cite{E33}. Its vertices are the twelve notes of the chromatic scale, and the edges connect two elements which differ by a fifth or a major third with the property that one may traverse all the edges of the graph passing exactly once by each note.   The article \cite{Papa-Topo} is a detailed survey of the work of Euler on these questions. In the present section,
 we start by reviewing the work of  Euler on the question of the seven bridges of K\"onigsberg. This work shows that  Euler considered himself as the direct heir of Leibniz for what concerns the field of topology. We shall then describe in detail the works of Euler and Descartes on the Euler characteristic,\index{Euler characteristic}\index{theorem!Euler characteristic}  a question which is directly related to the topological classification of surfaces, which was one of Riemann's major achievements in topology. We recall that Euler formulated this result for a surface which is the boundary of a convex polyhedron having  $F$ faces, $A$ edges and $S$ vertices; the formula is then: 
 \[F-A+S=2.\]

We start with the problem of the K\"onigsberg bridges.

In the eighteenth century, the city of K\"onigsberg\footnote{Today, the city of K\"onigsberg, called Kaliningrad, is part of a Russian exclave between Poland and Lithuania on the Baltic Sea.} consisted of four quarters separated by branches of the river Pregel and related by seven bridges. The famous ``problem of the seven bridges of K\"onigsberg" asks for a path  in that city that starts at a given point and returns to the same point after crossing once and only once each of the seven bridges. At the time of Euler, this was a popular question among the inhabitants of K\"onigsberg.
 
Euler showed that such a path does not exist. He presented his solution to the Saint Petersburg Academy of Sciences on August 26, 1735, and in the same year he wrote a memoir on the solution of a more general 
problem entitled \emph{Solutio problematis ad geometriam situs pertinentis} (Solution of a problem relative to the geometry of position) \cite{Euler-Solutio}. Euler learned about the problem from a letter, dated March 7, 1736, sent to him by Carl Leonhard Gottlieb Ehler,\index{Ehler, Carl Leonhard Gottlieb  (1685--1753)} one of his friends who was the mayor of Danzig\footnote{Today, Danzig is the city of Gdansk, in Poland.} and who had worked as an astronomer in Berlin.  Euler solved the problem just after he received the letter.  In a letter dated March 13, 1736, written to Giovanni  Marioni,\index{Marioni, Giovanni Jacopo  (1676--1755)} an Italian astronomer working at the court of Vienna,  Euler declares that he became interested in this question because he realized that the problem could not be solved using geometry, algebra or combinatorics, and that therefore he wondered whether ``it belonged to the `geometry of position,'   (\emph{geometria situs})  which 
 Leibniz has so much sought for."
 In the same letter, Euler announced to Marioni that after some thought, he found a proof which applies not only to that case, but to all similar problems.

   In the introduction of his paper, Euler writes (translation from \cite{Briggs}):
 \begin{quote}\small
 In addition to that branch of geometry which is concerned with magnitudes, and
which has always received the greatest attention, there is another branch, previously almost unknown, which Leibniz first mentioned, calling it the
geometry
of position. This branch is concerned only with the determination of position
and its properties; it does not involve measurements, nor calculations made
with them. It has not yet been satisfactorily determined what kind of problems
are relevant to this geometry of position, or what methods should be used in
solving them. Hence, when a problem was recently mentioned, which seemed
geometrical but was so constructed that it did not require the measurement of
distances, nor did calculation help at all, I had no doubt that it was concerned
with the geometry of position, especially as its solution involved only position,
and no calculation was of any use. I have therefore decided to give here the
method which I have found for solving this kind of problem, as an example of
the geometry of position.
 \end{quote}
 
Euler's work on this problem is commented in several articles and books.

 We now come to Euler's polyhedron formula, and we start with Descartes.

Long before Euler came out with his formula $F-A+S=2$  relating the faces ($F$), edges ($A$) and vertices ($S$) of a convex polyhedron,
 Descartes\index{Descartes, Ren\'e (1596--1650)}  obtained an equivalent result,  with a geometric proof, involving the solid angles and the dihedral angles between the faces.  Described in modern terms, Descartes' proof consists in computing in two different manners the total curvature of the boundary of the polyhedron.     Descartes\index{Descartes, Ren\'e (1596--1650)} wrote that proof at the age of 25, but did not publish it. 
The story of  Descartes'\index{Descartes, Ren\'e (1596--1650)} manuscript is interesting and we recall it now.

Descartes' manuscript was discovered in Hanover, among Leibniz's\index{Leibniz, Gottfried Wilhelm (1646--1716)} estate. The latter had copied Descartes'\index{Descartes, Ren\'e (1596--1650)} proof during a stay in Paris, in 1675 or 1676, presumably with the intention of publishing it. The original manuscript of Descartes,\index{Descartes, Ren\'e (1596--1650)} which carries the title \emph{Progymnasmata de solidorum elementis} (Preparatory exercises to the elements of solids) \cite{Costabel} is mentioned  in an inventory of papers which Descartes left in some chests in  
     Stockholm, the city where he died. The copy, made by  Leibniz, 
     carries the same title, with the additional mention   \emph{excerpta ex manuscripto Cartesii} (Excerpt from a  manuscript of Descartes).\index{Descartes, Ren\'e (1596--1650)}\index{Descartes, Ren\'e (1596--1650)} After the manuscript was discovered, a French translation was published by Foucher de Careil\index{Foucher de Careil, Louis-Alexandre (1826--1891)} in 1859, in a volume of unpublished works of Descartes. This publication contained errors, because
      Foucher, who did it, was not a mathematician. The edition is nevertheless interesting, and in the introduction to the volume \cite{Foucher}, Foucher recalls the story of the discovery. The story is also told by Adam\index{Adam, Charles (1857--1896)} in the commentaries of the volume of the Adam-Tannery edition of  Descartes' works containing this theorem (\cite{Descartes-Tannery}  tome X, p. 257--263).

In 1890, Jonqui\`eres\index{Fauque de Jonqui\`eres, E.-J.-Ph. (1820--1901)} presented to the Paris Acad\'emie des Sciences two Comptes Rendus notes entitled \emph{Sur un point fondamental de la th\'eorie des poly\`edres} (On a fundamental property of the theory of polyhedra) \cite{Jonc0} and  \emph{Note sur le th\'eor\`eme d'Euler dans la th\'eorie des poly\`edres} (Note on the theorem of Euler\index{Euler, Leonhard (1707--1783)} on the theory of polyhedra) \cite{Jonc1}, without being aware of the work of Descartes\index{Descartes, Ren\'e (1596--1650)} on this subject. After Jordan pointed out the existence of the work of Descartes\index{Descartes, Ren\'e (1596--1650)} in Foucher's edition, Jonqui\`eres  published other  Comptes Rendus notes on the work of Descartes, cf. \cite{Jonc2} \cite{Jonc3}  \cite{Jonc4}.  Poincar\'e,\index{Poincar\'e, Henri (1854--1912)} in his celebrated first memoir on \emph{Analysis situs} \cite{Poin1} attributes to  Jonqui\`eres the generalization of Euler's theorem to non-necessarily convex polyhedral surfaces. 
 
 There is a relatively recent  (1987) critical edition of Descartes' \emph{Progymnasmata} with a French translation, with notes and commentaries, by P. Costabel \cite{Costabel}.

Euler reported on his work on polyhedra in his correspondence with Goldbach. In a letter dated November 14, 1750, Euler\index{Euler, Leonhard (1707--1783)} informs his friend of the following two results which he refers to as Theorems 6 and 11 respectively:

\begin{quote}\small 

6. In any solid enclosed by plane surfaces the sum of the number of faces and the number of solid angles is greater by 2 than the number of edges.

11. The sum if all planar angles equals four times as many right angles as the number of solid angles, decreased by 8.
\end{quote}

The term ``solid enclosed by plane surfaces" refers to a convex polyhedron. The first result is the Euler characteristic formula,\index{Euler characteristic}\index{theorem!Euler characteristic} and the second one is a form of the Gauss-Bonnet theorem for the sphere. Euler writes:

\begin{quote}\small
I am surprised that these general properties in stereometry have not been noticed by anybody, as far as I know, but still more that the most important of them, viz., Th. 6 and Th. 11, are so hard to prove; indeed I still cannot prove them in a way that satisfies me.
\end{quote}

In the same letter, Euler gives several examples where the two theorems are satisfied.

  In his memoir \cite{E231}, entitled \emph{Demonstratio nonnullarum insignium proprietatum, quibus solida hedris planis inclusa sunt praedita}
(Proof of some of the properties of solid bodies enclosed by planes) and written one year after \cite{E230}, Euler gave proofs of the two results. In the introduction of \cite{E231}, he declares that his polyhedron formula \index{Euler characteristic}\index{theorem!Euler characteristic} is part of a more general research he is conducting on polyhedra. In fact, in the letter to Goldbach we mentioned, Euler\index{Euler, Leonhard (1707--1783)} announces a result on volumes of simplices in terms of their side lengths (a three-dimensional analogue of Heron's formula for the area of triangles), which he proves later in his  paper \emph{Demonstratio nonnullarum insignium proprietatum, quibus solida hedris planis inclusa sunt praedita} \cite{E231}. Euler writes:\footnote{Translation by C. Frances and D. Richeson.}
    \begin{quote}\small
     Although I had uncovered many properties which are common to all bodies enclosed by plane faces and which seemed to be completely analogous to those which are commonly included among the first principles of rectilinear plane figures, still, not without a great deal of surprise did I realize that the most important of those principles were so recondite that all the time and effort spent looking for a proof of them had been fruitless. Nor, when I consulted my friends, who are otherwise extremely versed in these matters and with whom I had shared those properties, were they able to shed any light from which I could derive these missing proofs. After the consideration of many types of solids I came to the point where I understood that the properties which I had perceived in them clearly extended to all solids, even if it was not possible for me to show this in a rigorous proof. Thus, I thought that those properties should be included in that class of truths which we can, at any rate, acknowledge, but which it is not possible to prove.
 \end{quote}

  One advantage of Euler's proof, compared to the one of Descartes, is that it shows in a clear way the combinatorial aspect of the problem, highlighting the notion of edges and faces of the polyhedron.

The proof that Euler\index{Euler, Leonhard (1707--1783)}  gives in \cite{E231} is based on an induction on the number of solid angles, reducing them by one at each step. He writes:
\begin{quote}\small
These proofs are in no way
inferior to those proofs used in Geometry except that here due to the
nature of solids one must use more imagination, in as much as solids
are being depicted on a flat surface.
\end{quote}

At the same time, Euler was laying down the foundations of combinatorial topology. He writes (Scholion to Proposition 4):
\begin{quote}\small 
  I admit that I have thus brought to
light only the first principles of Solid Geometry, on which this science
should be built as it develops further. No doubt it contains many outstanding qualities of solids of which we are so far completely ignorant.[...]
\end{quote}

Legendre,\index{Legendre, Adrien-Marie (1752--1833)} in his \emph{\'El\'ements de g\'eom\'etrie} (Note XII) published in 1794 \cite{Legendre}, gave a complete proof of Euler's theorem based on geometry. This proof is considered as one of the simplest, and it is repeated in more modern works, e.g. in Hopf \cite{Hopf}.

 A large number of  mathematicians commented on Euler's polyhedron formula, expanding some arguments in Euler's proofs, giving new proofs, and sometimes  comparing Euler's work with that of Descartes. To show the diversity of these works, we mention the papers by Andreiev \cite{Andreev}, Bertrand \cite{Bertrand1860}, Bouga\"\i ev \cite{Bougaiev}, Brianchon \cite{Brianchon}, Catalan \cite{Catalan1}, Cauchy \cite{Cauchy1813}, \cite{Cauchy1813-2}, Feil \cite{Feil},  Gergonne, \cite{Gergonne}, Grunert,  \cite{Grunert}, Jonqui\`eres \cite{Jonc0},  \cite{Jonc1}, \cite{Jonc2}, \cite{Jonc3}, \cite{Jonc4}, \cite{Jordan3}, \cite{Jordan4}, \cite{Jordan44}, Jordan \cite{Jordan1}, \cite{Jordan2},  Lebesgue \cite{Lebesgue}, Lhuillier \cite{Lhuillier}, 
 Poincar\'e \cite{Poin-CR-Euler},  \cite{Poin1}, Poinsot    \cite{Poinsot1}, \cite{Poinsot2},  Prouhet \cite{Prouhet1},  \cite{Prouhet2}, Steiner, \cite{Steiner}, Valat \cite{Valat}  and Thiel \cite{Thiel}.
We shall quote some of these works below.  We also mention that in 1858, the Paris \emph{Acad\'emie des Sciences} proposed as a subject for the 1861 \emph{Grand prix}: ``To improve, in some important point, the geometric theory of polyherda." M\"obius participated and presented a memoir (but did not get the prize).

In  1811, Cauchy\index{Cauchy, Augustin-Louis (1789--1857)} brought out a purely combinatorial proof of that theorem.\index{Euler characteristic}\index{theorem!Euler characteristic} In this proof, one starts by deleting a face of the polyhedron and reduces the problem to another one concerning a planar polygon.\footnote{A similar proof is given by Hilbert and Cohn-Vossen \cite{HCV} p. 290.}  In his article \emph{Recherches sur les poly\`edres} (Researches on polyhedra) \cite{Cauchy1813}, published in 1813, Cauchy writes: 
\begin{quote}\small
Euler\index{Euler, Leonhard (1707--1783)} has determined, in the Petersburg \emph{M\'emoires}, year 1758, the relation that exists between the various elements that compose the surface of a polyhedron; and Mr. Legendre, in his  \emph{\'El\'ements de G\'eom\'etrie}, proved in a much  simpler manner Euler's theorem, by considerations of spherical polygons. Having been led by some researches to a new proof of that theorem, I reached a theorem which is much more general than the one of Euler\index{Euler, Leonhard (1707--1783)}, whose statement is the following: 
\\
Theorem. If we decompose\footnote{Cauchy ``decomposes" the polyhedron by taking new vertices in the interior of the three-dimensional polyhedron (and not on the boundary surface)} a polyhedron in as many others as we wish, by taking at will new vertices in the interior, and if we represent by $P$ the number of new polyhedra thus formed, by $S$ the total number of vertices, including those of the initial polyhedron, by $F$ the total number of faces, and by $A$ the total number of edges, then we will have  \[S+F=A+P-1,\] that is, the sum of the number of vertices and that of faces will overpass by one the sum of the number of edges that of polyhedra.\footnote{Euler a d\'etermin\'e, dans les \emph{M\'emoires} de P\'etersbourg, ann\'ee 1758, la relation qui existe entre les diff\'erents \'el\'ements qui composent la surface d'un poly\`edre ; et M. Legendre,\index{Legendre, Adrien-Marie (1752--1833)} dans ses \emph{\'El\'ements de G\'eom\'etrie}, a d\'emontr\'e d'une mani\`ere beaucoup plus simple le th\'eor\`eme d'Euler, par la consid\'eration des polygones sph\'eriques. Ayant \'et\'e conduit par quelques recherches \`a une nouvelle d\'emonstration de ce th\'eor\`eme, je suis parvenu \`a un th\'eor\`eme plus g\'en\'eral que celui d'Euler et dont voici l'\'enonc\'e :
  \\
Th\'eor\`eme. Si l'on d\'ecompose un poly\`edre en tant d'autres que l'on voudra, en prenant \`a volont\'e dans l'int\'erieur de nouveaux sommets ; que l'on repr\'esente par $P$ le nombre de nouveaux poly\`edres ainsi form\'es, par $S$ le nombre total de sommets, y compris ceux du premier poly\`edre, par $F$ le nombre total de faces, et par $A$ le nombre total des ar\^etes, on aura \[S+F=A+P-1,\] c'est-\`a-dire que la somme faite du nombre des sommets et de celui des faces surpassera d'une unit\'e la somme faite du nombre des ar\^etes et de celui des poly\`edres.}
\end{quote}

Poinsot\index{Poinsot, Louis (1777-1859)} \cite{Poinsot1}, in 1858, published a proof of Euler's formula\index{Euler formula} using some of Cauchy's arguments.  He writes: ``This relation,
which Euler\index{Euler, Leonhard (1707--1783)} was the first to prove, does not hold only for convex polyhedra, as one might think, but for polyhedra of any kind." In fact, this statement needs some explanation. We are used today to the fact that Euler's formula is valid for polyhedra which are homeomorphic to a sphere. This notion did not exist at that time, neither the word, nor the idea. 
 One had to wait for that to the work of Jordan,\index{Jordan, Camille (1838--1922)} who set up the precise hypotheses under which Euler's formula is valid. In his article \cite{Jordan4}, he writes that Euler's theorem is valid for polyhedra which he calls ``simple," or ``Eulerian," that is, polyhedra (\cite{Jordan4} p. 35) ``such that any contour drawn on the surface which does not traverse itself  divides this surface into two separate regions; a category that contains as a particular case convex polyhedra."\footnote{[...] tels que tout contour ferm\'e trac\'e sur leur surface et ne se traversant pas lui-m\^eme divise cette surface en deux r\'egions s\'epar\'ees ; cat\'egorie qui renferme comme cas particulier les poly\`edres convexes.}  A few pages later (p. 38), Jordan makes the following commentary: ``It would have been easy to show that if we can draw on a polyhedron $\lambda$ different contours which do not intersect each other and which do not divide the surface into separate parts, we would have $S+F=A+2-2\lambda$."\footnote{Il serait ais\'e de d\'emontrer que si l'on peut tracer sur un poly\`edre $\lambda$ contours diff\'erents, ne se coupant pas mutuellement et ne divisant pas la surface en parties s\'epar\'ees, on aura $S+F=A+2-2\lambda$.} In fact, Jordan had extracted the notion we call today ``topological surface of finite type," to which the general theory applies, cf. \cite{Jordan3} p. 86: 
 \begin{quote}\small
 A surface is said to be of type $(m,n)$  if it is bounded by $m$ closed contours and if furthermore we can draw on it $n$ closed contours
 that do not intersect themselves nor mutually, without dividing it into two distinct regions.\footnote{Une surface sera dite d'esp\`ece $(m,n)$ si elle est limit\'ee par $m$ contours ferm\'es et si l'on peut d'autre part y tracer $n$ contours ferm\'es ne se coupant pas eux-m\^emes ni mutuellement, sans la partager en deux r\'egions distinctes.}
 \end{quote}
  Then Jordan makes the relation with the polyhedra to which Euler's formula applies: ``The polyherda of kind $(0,0)$ are nothing but those which I called \emph{Eulerian}."\footnote{Les poly\`edres de l'esp\`ece $(0,0)$ ne sont autres que ceux que j'ai appel\'es \emph{eul\'eriens}.}

  It is interesting to read Lebesgue's\index{Lebesgue, Henri (1875--1941)} comments on 
  some proof of Euler's theorem, because it gives us some hints of how the subject of topology was viewed in those days.  Lebesgues' comments are written in 1924 (\cite{Lebesgue} p. 319): 
 \begin{quote}\small
 I don't agree at all with those who pretend to attribute Euler's theorem to Descartes\index{Descartes, Ren\'e (1596--1650)}. Descartes did not state the theorem; he did not see it. Euler\index{Euler, Leonhard (1707--1783)} perceived it and he fully understood its character. For Euler, the description of the form of a polyhedron must precede the use of the measures of its elements, and this is why he set his theorem as a fundamental theorem. 
 For him, like for us, this is a theorem of enumerative \emph{Analysis situs};\index{analysis situs@\emph{analysis situs}} therefore he tried to find it by considerations independent of any metrical theory, that in effect belong to what we call the field of \emph{Analysis situs}. And this is why he left to Legendre\index{Legendre, Adrien-Marie (1752--1833)} the honor of finding a rigorous proof. None of us who had read a little bit of Euler\index{Euler, Leonhard (1707--1783)} and who were amazed by his prodigious technical masterliness will doubt, even for one second, that if Euler had thought of putting aside his theorem and deducing it from one of its metric corollaries, he would have easily succeeded. (It should be noted that Euler does not at all restrict his researches to convex polyhedra.) It seems to me, on the contrary, that the fact that Descartes\index{Descartes, Ren\'e (1596--1650)} passed so closely to the theorem without seeing it, emphasizes Euler's credit.  (At least, this is what we believe, because Descartes\index{Descartes, Ren\'e (1596--1650)} employed in his notebook some algebraic characters which he used before knowing Vi\`ete's characters.) But Leibniz,\index{Leibniz, Gottfried Wilhelm (1646--1716)} who found Descartes' notebook enough interesting to copy it, who realized that Descartes' geometry does not apply to questions involving order and position relations, who dreamed of constructing the algebra of these relations and who in advance gave it the name  \emph{Analysis situs},\index{analysis situs@\emph{analysis situs}}  did not notice, in Descartes' notebook Euler's theorem which is so fundamental in  \emph{Analysis situs}. This theorem really belongs to Euler\index{Euler, Leonhard (1707--1783)}. As for the proof, one could, may be with a little bit of unfairness, call it the proof of Legendre\index{Legendre, Adrien-Marie (1752--1833)} and Descartes.\index{Descartes, Ren\'e (1596--1650)} This proof is metrical, and it is fair to blame it for the fact that it uses notions that are foreign to  \emph{Analysis situs}. But one should not exaggerate the value of this grievance.\footnote{Je ne suis pas du tout d'accord avec ceux qui pr\'etendent attribuer \`a Descartes\index{Descartes, Ren\'e (1596--1650)} le th\'eor\`eme d'Euler. Descartes n'a pas \'enonc\'e le th\'eor\`eme ; il ne l'a pas vu. Euler\index{Euler, Leonhard (1707--1783)} l'a aper\c cu et en a bien compris le caract\`ere. Pour Euler, la description de la forme d'un poly\`edre doit pr\'ec\'eder l'utilisation des mesures de ses \'el\'ements et c'est pourquoi il a pos\'e son th\'eor\`eme comme th\'eor\`eme fondamental. C'est, pour lui comme pour nous, un th\'eor\`eme d'\emph{Analysis situs} \'enum\'erative ; aussi a-t-il cherch\'e \`a le d\'emontrer par des consid\'erations ind\'ependantes de toute propri\'et\'e m\'etrique, appartenant bien \`a ce que nous appelons le domaine de l'\emph{Analysis situs}. Et c'est pourquoi il a laiss\'e \`a Legendre\index{Legendre, Adrien-Marie (1752--1833)} l'honneur d'en trouver la preuve rigoureuse ;\index{Legendre, Adrien-Marie (1752--1833)} aucun de ceux qui ont quelque peu lu Euler, et qui ont \'et\'e stup\'efaits de sa prodigieuse virtuosit\'e technique, ne doutera un seul instant que si Euler\index{Euler, Leonhard (1707--1783)} avait pens\'e \`a faire passer son th\'eor\`eme au second plan et \`a le d\'eduire d'un de ses corollaires m\'etriques, il n'y e\^ut facilement r\'eussi. (Il convient d'ajouter qu'Euler ne restreint nullement ses recherches aux poly\`edres convexes.) Que Descartes\index{Descartes, Ren\'e (1596--1650)} soit pass\'e si pr\`es du th\'eor\`eme sans le voir me para\^\i t au contraire souligner le m\'erite d'Euler. Encore peut-on dire que Descartes \'etait jeune quand il s'occupait de ces questions. (C'est du moins ce que l'on croit, parce que Descartes\index{Descartes, Ren\'e (1596--1650)} a employ\'e dans son cahier certains caract\`eres cossiques qu'il utilisait avant de conna\^\i tre les notations de Vi\`ete.) Mais Leibniz qui a trouv\'e le cahier de Descartes assez int\'eressant pour le copier, qui a reconnu que la g\'eom\'etrie de Descartes ne s'appliquait pas aux questions o\`u interviennent des relations d'ordre et de position, qui a r\^ev\'e de construire l'alg\`ebre de ces relations et l'a nomm\'ee \`a l'avance \emph{Analysis situs}, n'a pas aper\c cu, dans le cahier de Descartes,\index{Descartes, Ren\'e (1596--1650)} le th\'eor\`eme d'Euler\index{Euler, Leonhard (1707--1783)} si fondamental en \emph{Analysis situs}. Le th\'eor\`eme appartient bien \`a Euler ; quant \`a la d\'emonstration, on pourrait, un peu injustement peut-\^etre, la d\'enommer d\'emonstration de Legendre\index{Legendre, Adrien-Marie (1752--1833)} et Descartes. Cette d\'emonstration est m\'etrique ; il est juste de lui reprocher de faire appel \`a des notions \'etrang\`eres \`a l'\emph{Analysis situs}. Mais il ne faudrait pas s'exag\'erer la valeur de ce grief.}
\end{quote}

 We now give a quick review of some work of Gauss on topology, another field in which his impact on Riemann was huge.
  
Gauss was interested in applications of \emph{Geometria situs} (a term he used in his writings), in particular in astronomy, geodesy and electromagnetism. In astronomy, he addressed the question of whether orbits of celestial bodies may be linked (cf. his short treatise entitled \emph{\"Uber die Grenzen der geocentrischen Orter der Planeten}).
From his work on geodesy, we mention his letter to Schumacher, 21 Nov. 1825, (from Gauss's \emph{Werke} vol. VIII, p. 400):
\begin{quote}\small

Some time ago I started to take up again a part of my general investigations on curved surfaces, which shall become the foundation of my projected work on higher geodesy. [...] Unfortunately, I find that I will have to go very far afield [...]. One has to follow the tree down to all its root threads, and some of this costs me week-long intense thought. Much of it even belongs to geometria situs, an almost unexploited field.

\end{quote}

From Gauss's \emph{Nachla\ss}, we know that he worked on a combinatorial theory of knot projections, during the year 1825, and again in 1844. (Gauss's \emph{Werke},  Vol. VIII, p. 271--286).   We already mentioned at the beginning of this section, that we learn from a 
 letter sent by Betti\index{Betti, Enrico (1923--1892)} to Tardy that the idea of analyzing a surface by performing successive cuts was given to Riemann by Gauss, in a private conversation.  
Besides Riemann, Gauss had two students who worked on topology and who were certainly influenced by him: Listing and M\"obius.

  Riemann introduced the fundamental topological notions for surfaces: connectedness, degree of connectivity, the classification of closed surfaces by their genus. He developed this theory for the purpose of using it in his work on the theory of functions of a complex variable.
    In his memoir on Abelian functions, he talks about \emph{analysis situs}, referring to Leibniz:
\begin{quote}\small
In the study of functions obtained by the integration of exact differentials, a few theorems of analysis situs are almost essential. Under that name, which was used by Leibniz, although may be in a slightly different sense, it is permitted to designate the theory of continuous magnitudes which studies these magnitudes, not as independent of their position and measurable with respect to each other, but by disregarding all idea of measure and studying them only for what regards their relation of position and inclusion. I intend to treat this subject later, in a way that is completely independent of any measure.
\end{quote}

In his habilitation dissertation\index{Riemann! habilitation lecture},\index{habilitation lecture!Riemann} Riemann mentions the possibility of working in the new field of topology, talking about the notion of ``place." We quote this cryptic passage:
\begin{quote}\small
Measure consists in the superposition of the magnitudes to be compared; it therefore requires a means of using one magnitude as the standard for another. In the absence of this, two magnitudes can only be compared when one is a part of the other; in which case we can only determine the more or less and not the how much. The researches which can in this case be instituted about them form a general division of the science of magnitude in which magnitudes are regarded not as existing independently of position and not as expressible in terms of a unit, but as regions in a manifoldness. Such researches have become a necessity for many parts of mathematics, e.g., for the treatment of many-valued analytical functions; and the want of them is no doubt a chief reason for which the celebrated theorem of Abel and the achievements of Lagrange, Pfaff,\index{Pfaff, Johann Friedrich (1765--1825)} Jacobi for the general theory of differential equations, have so long remained unfruitful. Out of this general part of the science of extended magnitude in which nothing is assumed but what is contained in the notion of it, it will suffice for the present purpose to bring into prominence two points; the first of which relates to the construction of the notion of a multiply extended manifoldness, the second relates to the reduction of determinations of place in a given manifoldness to determinations of quantity, and will make clear the true character of an $n$-fold extent.
\end{quote}
He also describes the passage from one dimension to another:
\begin{quote}\small

If in the case of a notion whose specialisations form a continuous manifoldness, one passes from a certain specialisation in a definite way to another, the specialisations passed over form a simply extended manifoldness, whose true character is that in it a continuous progress from a point is possible only on two sides, forward or backwards. If one now supposes that this manifoldness in its turn passes over into another entirely different, and again in a definite way, namely so that each point passes over into a definite point of the other, then all the specialisations so obtained form a doubly extended manifoldness. In a similar manner one obtains a triply extended manifoldness, if one imagines a doubly extended one passing over in a definite way to another entirely different; and it is easy to see how this construction may be continued. If one regards the variable object instead of the determinable notion of it, this construction may be described as a composition of a variability of $n+1$  dimensions out of a variability of $n$  dimensions and a variability of one dimension. 
\end{quote}

  Riemann's ideas on topology are explained in some sections of Klein's booklet \cite{Klein-Riemann}. For instance, \S 8 carries the title \emph{Classification of closed surfaces according to the value of the integer $p$}. 
     Let us comment on a passage from Klein's booklet \cite{Klein-Riemann} concerning the classification of surfaces, as an example of his point of view on topology. We know that topology, which was an emerging subject, plays an important role in Riemann surface theory. We already mentioned that Riemann introduced several major notions on surface topology. Klein tried to make a more systematic exposition of these ideas. His book \cite{Klein-Riemann} contains a chapter in which the classification of closed surfaces according to genera is presented. On p. 32,  he writes: 
\begin{quote}\small
That it is impossible to represent surfaces having different $p$'s\footnote{We recall that $p$ denotes the genus.} upon one another, the correspondence being uniform, seems evident. It is not meant, however, that this kind of geometrical certainty needs no further investigation; cf. the explanations of G. Cantor (\emph{Crelle}, t. {\sc lxxxiv}. pp. 242 \emph{et seq.}). But these investigations are meanwhile excluded from consideration in the text, since the principle there insisted upon is to base all reasonings ultimately on intuitive relations.
\end{quote}  Klein then states the converse: between any two surfaces of the same genus it is possible to find a ``uniform correspondence." He declares that this statement is more difficult to prove, and he refers to the 1866 article by Jordan \cite{Jordan1866}.\footnote{Camille Jordan (1838-1922), who is mostly known for his results on the topology of surfaces and on group theory, also worked on function theory in the sense of Riemann. The title of the second part of his doctoral thesis is: ``On periods of inverse functions of integrals of algebraic differentials." The subject was proposed to him by Puiseux, whom we mention in this paper concerning uniformization. Jordan is among the first who tried to study the ideas of Galois, and he is also among the first who introduced group theory in the study of differential equations.} This paper is an important milestone in the history of topology, because it contains the first attempt to classify surfaces up to homeomorphism, although there was no precise definition of homeomorphism yet.\footnote{The word ``homeomorphism" was introduced by Poincar\'e\index{Poincar\'e, Henri (1854--1912)} in his article \cite{Poin1} but with a meaning that is different from the one it has today. There is a definition of homeomorphism in the 1909 article by Hadamard  \cite{Hadamard-Notions}, as being a one-to-one continuous map. This is not correct, unless Hadamard meant, by ``one-to-one continuous", ``one-to-one bi-continuous." We refer the reader to the paper \cite{Moore-H} on the rise and the development of the notion of homeomorphism. This paper contains several quotes, some of which are very intriguing.} Jordan's aim, in his paper, is to prove the following theorem, which he states in the introduction:
\begin{quote}\small
In order for two surfaces or pieces of flexible and inextensible surfaces to be be applied onto each other without tear and duplication, the following is necessary and sufficient:
\begin{enumerate}
\item That the number of separated contours that respectively bound  these portions of surfaces be the same. (If the surfaces considered are closed, this number is zero.)
\item That the maximal number of closed contours which do not self-intersect or intersect each other that we can draw on each of the two surfaces without cutting it into two separate regions be the same on both sides.
\end{enumerate}
\end{quote}

Jordan gives the following ``definition" of two surfaces being ``applicable onto each other." For a modern reader, this definition may seem fuzzy, but one has to remember that this paper was written in the heroic epoch of the foundations of modern topology, that the  notion of homeomorphism seems extremely natural for us today, but that it was not so in those times. Jordan writes:

\begin{quote} \small We shall rely on the following principle, which we can consider as evident, and take it if necessary as a definition: \emph{Two surfaces $S, S'$ are applicable onto each other if we can decompose them into  infinitely small elements such that to any contiguous elements of $S$ correspond contiguous elements of $S'$.}
\end{quote}

Besides Klein's booklet, several books and treatises explaining Riemann's ideas appeared in the decades that followed Riemann's work. We mention as examples Neumann's \emph{Vorlesungen \"uber Riemann's Theorie der Abel'schen Integrale} (Lectures on Riemann's theory of Abelian integrals) \cite{Neumann-Vorlesungen}, Picard's \emph{Trait\'e d'Analyse} \cite{Picard-cours}, Appell and Goursat's \emph{Th\'eorie des fonctions alg\'ebriques et de leurs int\'egrales} (Theory of Algebraic functions and their integrals) \cite{Appell-Goursat}, and  there are others. The last two treatises, together with several other French books on the theory of functions of a complex variable, are reviewed in Chapter 8 of the present volume \cite{Papa-Riemann3}.

 Among the other important topological notions that were introduced before Riemann and that were used by him, we must mention the notion of \emph{homotopy of paths} and its use in complex analysis (in particular, in the theory of line integrals), especially by Cauchy and Puiseux. This is discussed in detail in Chapter 7 of the present volume \cite{Papa-Puiseux}. Cauchy published his first work on the subject in 1825 \cite{Cauchy-1825}. This is also a topic on which Gauss was a forerunner, but he did not publish anything about it. This is attested in his letter to Bessel, December 18, 1811, published in Volume VIII of his \emph{Collected works} (p. 90--92), a letter in which Gauss makes the important remark that if one defines integrals along paths in the complex plane, then  the value obtained may depend on the path.

Regarding the history of Riemann's ideas on topology, we could have also commented on his predecessors regarding the notion of the discreteness and continuity of space, but this would have taken us too far. We make a few remarks on this matter in Chapter 6 of the present volume \cite{Papa-Riemann1} .

We end this section by quoting Alexander Grothendieck,\index{Grothendieck, Alexandre (1928--2014)} from his  \emph{R\'ecoltes et semailles} (Harvesting and Sowing),\footnote{The complete title is: \emph{R\'ecoltes et semailles : R\'eflexions et t\'emoignage sur un pass\'e de math\'ematicien} (Harvesting and Sowing: Reflections and testimony on the past of a mathematician). This is a long manuscript by Grothendieck in which he meditates on his  past as a mathematician and where he presents without any compliance his vision of the mathematical milieu in which he evolved, especially the decline in morals, for what concerns intellectual honesty.}  commenting on Riemann's reflections on this theme (\cite{RS} Chapter 2 \S 2.20):
\begin{quote} \small
It must be already fifteen or twenty years ago since, skimming the modest volume that constitutes Riemann's complete works, I was struck by a remark which he made ``incidentally." He observes there that it might be that the ultimate structure of space is ``discrete," and that the ``continuous" representations which we make of it may be an oversimplification (which may turn out to be excessive on the long run...) of a more complex reality; that for human thought, ``the continuous" was easier to grasp than ``the discontinuous," and that it serves us, subsequently, as an ``approximation" for apprehending the discontinuous. This is an amazingly penetrating remark expressed by a mathematician, at a time where the Euclidean model of physical space was not questioned. In a strictly logical sense, it was rather the discontinuous which, traditionally, served as a technical mode of approaching the continuous.

Moreover, the mathematical developments of the last decades showed an even more intimate symbiosis between continuous and discontinuous structures, which was not yet visualized in the first half of this century. The fact remains that finding a model which is ``satisfactory" (or, if need be, a collection of such models, linked in the most possible satisfying way...), whether the latter is ``continuous," ``discrete," or of a ``mixed" nature -- such a task will surely involve a great conceptual imagination, and a consummate intuition, in order to apprehend and disclose mathematical structures of a new type. This kind of imagination and ``intuition" seems to me a rare object, not only among physicists (where Einstein\index{Einstein, Albert (1879---1955)} and Schr\"odinger\index{Schr\"odinger, Erwin (1887--1961)} seem to be among the rare exceptions), but even among mathematicians (and I am talking in full knowledge of the facts).

To summarize, I foresee that the long-awaited  renewal (if ever it comes...) will rather come from someone who has the soul of a mathematician, who is well informed about the great problems of physics, rather than from a physicist. But above all, we need a man having the ``philosophical openness" required to grasp the crux of the problem. The latter is not at all of a technical nature, but it is a fundamental problem of ``natural philosophy."\footnote{Il doit y avoir d\'ej\`a quinze ou vingt ans, en feuilletant le modeste volume constituant l'\oe uvre compl\`ete de Riemann, j'avais \'et\'e frapp\'e par une remarque de lui ``en passant." Il y fait observer qu'il se pourrait bien que la structure ultime de l'espace soit ``discr\`ete", et que les repr\'esentations ``continues" que nous en faisons constituent peut-\^etre une simplification (excessive, peut-\^etre, \`a la longue...) d'une r\'ealit\'e plus complexe ; que pour l'esprit humain, ``le continu"  \'etait plus ais\'e \`a saisir que ``le discontinu", et qu'il nous sert, par suite, comme une ``approximation" pour appr\'ehender le discontinu. 
C'est l\`a une remarque d'une p\'en\'etration surprenante dans la bouche d'un math\'ematicien, \`a un moment o\`u le mod\`ele euclidien de l'espace physique n'avait jamais \'et\'e mis en cause ; au sens strictement logique, c'est plut\^ot le discontinu qui, traditionnellement, a servi comme mode d'approche technique vers le continu.

Les d\'eveloppements en math\'ematique des derni\`eres d\'ecennies ont d'ailleurs montr\'e une symbiose bien plus intime entre structures continues et discontinues, qu'on ne l'imaginait encore dans la premi\`ere moiti\'e de ce si\`ecle. Toujours est-il que de trouver un mod\`ele ``satisfaisant" (ou, au besoin, un ensemble de tels mod\`eles, se ``raccordant" de façon aussi satisfaisante que possible...), que celui-ci soit ``continu," ``discret" ou de nature ``mixte" -- un tel travail mettra en jeu s\^urement une grande imagination conceptuelle, et un flair consomm\'e pour appr\'ehender et mettre \`a jour des structures math\'ematiques de type nouveau. Ce genre d'imagination ou de ``flair" me semble chose rare, non seulement parmi les physiciens (o\`u Einstein et Schr\"odinger semblent avoir \'et\'e parmi les rares exceptions), mais m\^eme parmi les math\'ematiciens (et l\`a je parle en pleine connaissance de cause).

Pour r\'esumer, je pr\'evois que le renouvellement attendu (s'il doit encore venir \ldots) viendra plut\^ot d'un math\'ematicien dans l'\^ame, bien inform\'e des grands probl\`emes de la physique, que d'un physicien. Mais surtout, il y faudra un homme ayant ``l'ouverture philosophique" pour saisir le n\oe ud du probl\`eme. Celui-ci n'est nullement de nature technique, mais bien un probl\`eme fondamental de ``philosophie de la nature."}

\end{quote}

      \section{Differential geometry}\label{s:diff}

  In this section, 
we shall review some milestones in the history of differential geometry, concerning especially  the study of geodesics and of curvature, from its beginning until the work of Riemann.

               Differential geometry starts with the study of differentiable curves. The notion of curvature of planar curves already appears in works of Newton and of Johann I and Jakob Bernoulli. We mentioned,  in \S \ref{s:complex}, Johann Bernoulli's\index{Bernoulli, Johann (1667--1748)} paper \cite{Bernoulli1718} published in 1718 on the isoperimetry problem in the plane. 
               
            Differential geometry  is also the study of curvature.    In 1744, Euler published a book \cite{E65} in which he  sets the  bases of the calculus of variations. The title is \emph{Methodus inveniendi lineas curvas maximi minimive proprietate gaudentes, sive solutio problematis isoperimetrici lattissimo sensu accepti}
(Method for finding curved lines enjoying properties of maximum or minimum, or solution of isoperimetric problems in the broadest accepted sense). In that book, several applications of the new methods are presented, among them isoperimetry problems, the problem of  finding the shape of the brachystochrone, the study of the catenoid, and that of finding geodesics between two points on a surface. With this work of Euler, the methods of differential calculus, more precisely those of finding minima and maxima, were suddenly  generalized to the realm of a variable moving in an infinite dimensional space (even though the expression ``infinite dimensional" was not there yet), namely, in the question of looking for curves of minimal length or satisfying other geometric properties,   
among all curves joining two points.              
               
             Riemann's differential geometry is essentially about curvature, actually, the curvature of space, and we must talk now about the history of curvature, which starts with curvature of curves and surfaces. The history starts again with Euler.
             
              Volume II of Euler's\index{Euler, Leonhard (1707--1783)} \emph{Introductio in analysin infinitorum} \cite{Euler-Int-b}, published in 1748,  is concerned with the differential geometry of space curves  and  surfaces. Curves are given there by parametric equations of the form $x=x(t), y=y(t), z=z(t)$,  and surfaces by parametric equations of the form $x=x(u,v), y=y(u,v), z=z(u,v)$. It is possible that this is the first time where such a parametric representation of surfaces appears in print. About twenty years after the first edition of this treatise was published, Euler wrote a memoir entitled \emph{Recherches sur la courbure des surfaces} (Researches on the curvature of surfaces) \cite{E333} (1767), another work which transformed the subject. The aim of this memoir was to introduce and study the curvature\index{curvature} at a point on a surface. Euler's idea, which is very natural, was to introduce a notion of curvature at a point of a differentiable surface based on the curvature of curves that pass through that point. His intuition was that to understand curvature at a point of a surface, it suffices to study the curvature of curves that are intersections of that surface with Euclidean planes. Moreover, he showed that it is sufficient to consider the planes that are perpendicular to the surface, that is, the planes containing the normal vector to the surface at that point. Each such curve has an osculating circle, and the collection of radii of these circles contains all the information about the curvature of the surface at that point. Furthermore, Euler\index{Euler, Leonhard (1707--1783)}   proved that at any given point on the surface, the maximal curvature  and the minimal curvature associated to the normal planes determine all the other normal curvatures. To be more precise, given a point on the surface and a tangent vector $v$ at that point, let us call \emph{normal curvature} though $v$  the curvature of a curve obtained by intersecting the surface with a plane containing the vector $v$ and the normal vector at that point. The \emph{maximal} and \emph{minimal normal curvature} at the given point are the maximum and minimum of the normal curvatures taken over all the normal planes at that point. Likewise, the \emph{normal curvature radius} at the given point in the direction of the vector $v$ is the curvature radius of the associated curve. We have a similar notion of  \emph{maximal} and \emph{minimal normal curvature radii} at the given point. 
                  Euler showed that the directions of the planes that realize these extremal curvatures  (except in very special cases) are orthogonal to each other, and he proved that at a given point, if $\rho_1$ and $\rho_2$ are the maximal and minimal normal curvature radii respectively, then the normal curvature radius $\rho$ of the normal section through an arbitrary tangent vector $v$ is given by the equation
$$\rho={2\rho_1\rho_2\over (\rho_1+\rho_2)-(\rho_1-\rho_2)\cos(2\varphi)}\,,$$
where $\varphi$ is the angle between $v$ and the tangent vector to the normal plane with maximal curvature radius. 
 
It is usual to write Euler's equation in the following form:

$${1\over \rho}={\cos^2\varphi\over \rho_1}+{\sin^2\varphi\over \rho_2}\,.$$
We note that $\rho_1$ and $\rho_2$ may also take negative values and that Euler's equation has also a meaning when $\rho_1$ or $\rho_2$ is infinite; in the latter case, the curvature ${1\over \rho}$ is zero for all $\varphi$. There is a classical local classification of differentiable surfaces at a point  in terms of the signs of $\rho_1$ and $\rho_2$.

Euler writes (\cite{E333}, R\'eflexion VI, p.~143):
\begin{quote}\small
Thus, the judgement of the curvature of surfaces, however complicated it seems at the beginning, is reduced for each point to the knowledge of two osculating radii, one of which is the largest and the other the smallest at that element. These two objects determine entirely the nature of the curvature, displaying for us the curvature of all the possible sections that are perpendicular to the proposed element.\footnote{Ainsi le jugement sur la courbure des surface, quelque compliqu\'e qu'il ait paru au commencement, se r\'eduit pour chaque \'el\'ement \`a la conaissance de deux rayons osculateurs, dont l'un est le plus grand et l'autre le plus petit dans cet \'el\'ement ; ces deux choses d\'eterminent enti\`erement la nature de la courbure en nous d\'ecouvrant la courbure de toutes les sections possibles qui sont perpendiculaires sur l'\'el\'ement propos\'e.} 
\end{quote}

There are other memoirs by Euler on the differential geometry of surfaces. 
We mention his \emph{Solutio trium problematum difficiliorum ad methodum tangentium inversam pertinentium} (Solution of three rather difficult
problems pertaining to the inverse method of tangents) \cite{E771} published in 1826, that is, several years after Euler's death. In this memoir, Euler addresses ``inverse problems" in differential geometry, e.g., to reconstruct curves from information on their tangents.  We also mention Euler's \emph{De solidis quorum superficiem in planum explicare licet}
(On solids whose entire surface can be unfolded onto a plane) \cite{E419} in which for the first time the notion of a surface developable on the plane is introduced. This notion was thoroughly investigated in the later works of Monge and his students that we mention below, and much  later by Eugenio Beltrami. 
This paper \cite{E419} on developable surfaces also addresses a so-called ``inverse problem," namely, the question of giving a characterization of the surfaces that are developable.

Gauss's development of differential geometry attained a high degree of perfection in the 1820s, motivated by his works on geography, astronomy and geodesy. He was probably the first to formulate the question of finding the properties of surfaces which are independent of their embedding in 3-space. After Euler in \cite{E333} highlighted the role of the maximal and minimal curvature at a point of a surface, it was Gauss's idea to take the product of these quantities as a measure of the curvature at that point, and to show that the result   is an isometry invariant of the surface. This is (expressed in modern terms) the content of  Gauss's \emph{Theorema egregium}\index{Theorema egregium@\emph{Theorema egregium} (Gauss)}.\index{Gauss!\emph{Theorema egregium}} The result is contained in  
             Gauss's \emph{Disquisitiones generales circa superficies curvas}\index{Gauss!Disquisitiones generales circa superficies curvas} (General investigations on curved surfaces) \cite{Gauss-Pesic} (1828).  In the abstract he wrote for his \emph{Disquisitiones} \cite{Gauss-Pesic} (1828), Gauss writes (translation from \cite{Gauss-Pesic} p. 48):
 \begin{quote}\small
[...] These theorems lead to the consideration of the theory of curved surfaces from a new point of view, where a wide and still wholly uncultivated field is open to investigation. If we consider surfaces not as boundaries of bodies, but as bodies of which one dimension vanishes, and if at the same time we conceive them as flexible but not extensible, we see that two essentially different relations must be distinguished, namely, on the one hand, those that presuppose a definite form of the surface in space; on the other hand, those that are independent of the various forms that the surface may assume. This discussion is concerned with the latter. In accordance with what has been said, the measure of curvature belongs to this case. But it is easily seen that the consideration of figures constructed upon the surface, their angles, their areas and their integral curvatures, the joining of the points by means of shortest lines, and the like, also belong to this case. All such investigations must start from this, that the very nature of the curved surface is given by means of the expression of any linear element in the form $\sqrt{Edp^2+2Fdpdq+Gdq^2}$.
\end{quote}

             Gauss, in the same memoir, used the so-called Gauss map from a surface to the unit sphere, defined by sending the normal unit normal vector at a point to the corresponding point on the sphere and showing that one can recover the curvature of the surface, which he had defined as the product of the minimal and maximal curvatures. The curvature, using the Gauss map, is obtained by taking the ratio of the area of the image of the Gauss map by the area of the surface (the definition of the curvature at a point needs a passage to the infinitesimal level).
             
 Riemann's most important articles on differential geometry are his habilitation lecture \emph{\"Uber die Hypothesen welche der Geometrie zu Grunde liegen} (1854)
 which is mentioned several times in the present paper and in the other chapters of this book, 
 and the 
\emph{Commentatio mathematica,  qua  respondere tentatur quaestioni  ab Ill${}^{\mathrm{ma}}$    Academia Parisiensi propositae} (A mathematical note attempting to answer a question posed by the distinguished Paris Academy), a memoir which he wrote in 1861, at the occasion of his participation to a competition prize set by the Paris Academy of Sciences, and which we consider in some detail in Chapter 6 of the present volume \cite{Papa-Riemann1}. These two memoirs are unusual for opposite reasons: the first one lacks of formulae, and the second one is full of them. The second memoir contains the explicit form of the object which we call today the Riemann tensor.

         In his habilitation lecture, Riemann makes a clear reference to  Gauss's \emph{Disquisitiones} as one of his major sources of inspiration, a work which he  
 includes however in a broad philosophical  discussion on magnitude, measure, quantity and the possibility of geometric representation. It is always good to read Riemann:
          
\begin{quote}\small
 Having constructed the notion of a manifoldness of $n$ dimensions, and found that its true character consists in the property that the determination of position in it may be reduced to $n$ determinations of magnitude, we come to the second of the problems proposed above, viz. the study of the measure-relations of which such a manifoldness is capable, and of the conditions which suffice to determine them. These measure-relations can only be studied in abstract notions of quantity, and their dependence on one another can only be represented by formul\ae . On certain assumptions, however, they are decomposable into relations which, taken separately, are capable of geometric representation; and thus it becomes possible to express geometrically the calculated results. In this way, to come to solid ground, we cannot, it is true, avoid abstract considerations in our formul\ae , but at least the results of calculation may subsequently be presented in a geometric form. The foundations of these two parts of the question are established in the celebrated memoir of Gauss, \emph{Disqusitiones generales circa superficies curvas}.
\end{quote}

For the case of surfaces, he writes:
\begin{quote}\small
In the idea of surfaces, together with the intrinsic measure-relations in which only the length of lines on the surfaces is considered, there is always mixed up the position of points lying out of the surface. We may, however, abstract from external relations if we consider such deformations as leave unaltered the length of lines -- i.e., if we regard the surface as bent in any way without stretching, and treat all surfaces so related to each other as equivalent. Thus, for example, any cylindrical or conical surface counts as equivalent to a plane, since it may be made out of one by mere bending, in which the intrinsic measure-relations remain, and all theorems about a plane -- therefore the whole of planimetry -- retain their validity. On the other hand they count as essentially different from the sphere, which cannot be changed into a plane without stretching. According to our previous investigation the intrinsic measure-relations of a twofold extent in which the line-element may be expressed as the square root of a quadric differential, which is the case with surfaces, are characterized by the total curvature. Now this quantity in the case of surfaces is capable of a visible interpretation, viz., it is the product of the two curvatures of the surface, or multiplied by the area of a small geodesic triangle, it is equal to the spherical excess of the same. The first definition assumes the proposition that the product of the two radii of curvature is unaltered by mere bending; the second, that in the same place the area of a small triangle is proportional to its spherical excess. To give an intelligible meaning to the curvature of an $n$-fold extent at a given point and in a given surface-direction through it, we must start from the fact that a geodesic proceeding from a point is entirely determined when its initial direction is given. According to this we obtain a determinate surface if we prolong all the geodesics proceeding from the given point and lying initially in the given surface-direction; this surface has at the given point a definite curvature, which is also the curvature of the n-fold continuum at the given point in the given surface-direction. 
\end{quote}

In these passages, Riemann summarizes his ideas on the general metric on (what became known later on as) an $n$-dimensional differentiable manifold,  defined by a quadratic form on each tangent space, a broad generalization of Gauss's investigations on surfaces in which  the quadratic form determines the metric, permits to calculate distances, angles and the curvature at any point. The curvature is the product of two quantities  and is invariant by bending. The quadratic form represents the square of the line element. With these tools, one can study geodesic triangles on surfaces, prove that a geodesic is  determined by its initial vector, generalize these matters to immersed surfaces, etc. 
                       
One may also include in Riemann's list of works on differential geometry his two papers on minimal surfaces \cite{Riemann-minimal} and \cite{Riemann-Untersuchungen}. They are reviewed in Chapter 5 of the present volume, \cite{Yamada}. We also note that in his doctoral dissertation defended in 1880 in Paris and written under the supervision of  Bonnet\index{Bonnet, Pierre-Ossian (1819--1892)}   (cf. \cite{Niewenglowski}), Niewenglowski explains that Riemann, in his work on minimal surfaces, was inspired by Bonnet; cf. also the comments in Chapter 8 of the present volume \cite{Papa-Riemann3}. Again, minimal surfaces first appear in the work of Euler (cf. Chapter V, \S 47 of Euler's first treatise on the calculus of variations, \cite{E65}).

In a longer survey on Riemann's predecessors in the field of differential geometry, one would have analyzed the works of several French mathematicians who stand between Euler and Gauss, in particular Gaspard Monge (1764--1818),\index{Monge, Gaspard (1764--1818)}
 Jean-Baptiste Meusnier\index{Meusnier  de La Place, Jean-Baptiste Marie (1754--1793)}  (1754--1793), Sim\'eon-Denis Poisson\index{Poisson, Sim\'eon-Denis (1781--1840)} (1781--1840), 
 Charles Dupin\index{Dupin, Charles (1784--1873)} (1784--1873),  Olinde Rodrigues\index{Rodrigues, Olinde (1795--1851)} (1795--1851) and there are  others. We only mention some of these works.

  Monge, who was the founder of a famous school on projective and differential geometry, continued Euler's work on developable surfaces, cf. \cite{Monge-deve1} and \cite{Monge-deve2}. He worked in particular with two orthogonal line fields that are defined by Euler's minimal and maximal directions of curvature radii, and he coined the expression \emph{umbilical point} for points where the two curvature radii have the same value. (On the sphere, every point is  umbilical.) Monge expressed several times in his writings his debt to Euler. In \cite{Monge-deve2}, he writes:
 \begin{quote}\small
 Having resumed this matter, at the occasion of a memoir that Mr. Euler gave in the 1771 volume of the Petersburg Academy on developable surfaces, in which this famous geometer gives formulae to determine whether a given surface  may or may not be applied onto a plane, I reached results on the same subject which seem to me much simpler, and whose usage is much easier.\footnote{Ayant repris cette mati\`ere, \`a l'occasion d'un m\'emoire que M. Euler a donn\'e dans le volume 1771 de l'Acad\'emie de P\'etersbourg, sur les surfaces d\'eveloppables, et dans lequel cet illustre g\'eom\`etre donne des formules pour reconna\^itre si une surface courbe propos\'ee jouit ou non de la propri\'et\'e de pouvoir \^etre appliqu\'ee sur un plan, je suis parvenu \`a des r\'esultats qui me semblent beaucoup plus simples, et d'un usage bien plus facile pour le m\^eme sujet.}
 \end{quote}

 Poisson\index{Poisson,  Sim\'eon Denis (1781--1840)} was a student of Lagrange and Laplace. He wrote a memoir entitled \emph{M\'emoire sur la courbure des surfaces} (Memoir on the curvature of surfaces) \cite{Poisson} (1832) in which he studied singular points of the curvature. We mention by the way that there are several points in the work of Poisson that are related to Riemann's works, in particular, concerning definite integrals and Fourier series.

 Meusnier was a student of Monge. He gave a formula for the curvature of a curve obtained by intersecting a surface by a non-normal section, in terms of that of the normal sections at the given point that were considered by Euler. His paper on the subject carries almost the same title as Euler's paper \cite{E333},   \emph{M\'emoire sur la courbure des surfaces}  \cite{Meusnier1785} (1785). In this paper, Meusnier writes (p. 478):
\begin{quote}\small
Mr. Euler treated the same matter in a very beautiful memoir, printed in 1760 among those of the Berlin Academy. This famous geometer addresses the question in a manner which is different from the one which we just presented. He makes the curvature of an element of the surface dependent on that of the various sections that one can perform on it by cutting it with planes.\footnote{M. Euler a trait\'e la m\^eme mati\`ere dans un fort beau m\'emoire, imprim\'e en 1760 parmi ceux de l'Acad\'emie de Berlin. Cet illustre g\'eom\`etre envisage la question d'une mani\`ere diff\'erente de celle que nous venons d'exposer ; il fait d\'ependre la courbure d'un \'el\'ement de surface de celle des diff\'erentes sections qu'on peut y faire en le coupant par des plans.}
\end{quote}

 Meusnier's work is surveyed by Darboux in the paper \cite{Darboux-Meusnier}. 
 
 Dupin is mostly known for the so-called Dupin indicatrix, a geometric device which characterizes the shape of a surface at a given point and which turned out to be related to the Gaussian curvature at that point. His famous work on geometry bears the title \emph{D\'eveloppements de g\'eom\'etrie, avec des applications \`a la stabilit\'e des vaisseaux, aux d\'eblais et remblais, au d\'efilement, \`a l'optique, etc. pour faire suite \`a la g\'eom\'etrie descriptive et \`a la g\'eom\'etrie analytique de M. Monge : Th\'eorie.} (Developments in geometry, with applications to the stability of vessels, cuts and fills, scrolling, optics, etc. as a sequel to the descriptive and analytic geometry of Mr. Monge: Theory) \cite{Dupin}, 1813.

 Rodrigues introduced, before Gauss, the Gauss map, and showed that at a given point the ratio of the area of the image to the area on the surface is equal (at the infinitesimal level) to the product of the two principal curvatures (those defined by Euler), cf. \cite{Rodrigues}. This result pre-dates that of Gauss, but the fact that curvature is an isometry invariant is however absent from Rodrigues' work.

Finally, it is fair to mention that Riemann's  work in high dimensions was prepared by works of other mathematicians   done in higher dimensions, and we mention in particular Jacobi \cite{Jacobi-binis} on the reduction of pairs of quadratic forms, Grassmann on higher-dimensional linear algebra \cite{Grassmann1844}, and Cayley \cite{Cayley-higher} on higher-dimensional analytic geometry.

 \section{Trigonometric series} \label{s:trigo} 
    
 Riemann's habilitation dissertation (\emph{Habilitationsschrift})\index{Riemann! habilitation text}\index{habilitation text!Riemann} is concerned with trigonometric series. Its title is  \emph{\"Uber die Darstellbarkeit einer Function durch eine trigonometrische Reihe} (On the representability of a function by a trigonometric series)\index{trigonometric series} \cite{Riemann-Trigo}. The dissertation is divided into two parts. The first part is a survey of the history of the problem of representing by a trigonometric series a function which is arbitrary, but which -- according to Riemann's words at the beginning of his memoir -- is ``given graphically."  It is important to note the last fact, since, for instance, a function which is discontinuous at a dense set of points cannot be given graphically. Dealing with the most general functions is part of the subject of the second part of Riemann's memoir.

 We shall report on the historical part of Riemann's memoir. It involves several mathematicians, in particular Euler, although not directly his work on trigonometric series. We therefore note right away that Euler uses trigonometric series in his two memoirs on astronomy \emph{Recherches sur la question des in\'egalit\'es du mouvement de Saturne et de Jupiter} (Researches on the question of the inequalities in the motion of Saturn and Jupiter) 
 \cite{E120}, and \emph{De motu corporum coelestium a viribus quibuscunque perturbato}
(On the movement of celestial bodies perturbed by any number of forces) \cite{E232},  both presented  for competitions proposed by the Paris Acad\'emie des Sciences, in 1748 and 1756 respectively. The two memoirs are analyzed in the paper \cite{Golland} in which the authors show that Euler was much more concerned with convergence of series than what is claimed in several books and articles on his work.

The starting point of Riemann's historical survey is the controversy between Euler\index{Euler, Leonhard (1707--1783)} and d'Alembert\index{Alembert@d'Alembert, Jean le Rond (1717--1783)}  which originated in  the publication in 1747 of a memoir by the latter,\index{Alembert@d'Alembert, Jean le Rond (1717--1783)}  \emph{Recherches sur la courbe que forme une corde tendue mise en vibrations} (Researches on the curve formed by a taut string subject to vibrations) \cite{Al_recherches}. In this  memoir, d'Alembert gave the partial differential equation that represents the motion of a point on a vibrating string\index{vibrating string!equation}\index{equation!vibrating string} subject to small vibrations:\index{string vibration}\index{vibrating string} 
  \[  \frac{\partial^2y}{\partial t^2}=\alpha^2\frac{\partial^2y}{\partial x^2}.
\]
 In this equation, $t$ represents time, $\alpha$ is a constant and $y=y(t)$ is the oscillation of the string at a point whose coordinate is $x$ along the string. The main problem, after this discovery, was to characterize the functions that are solutions of that equation. 
 
 One first obvious (but wrong) guess for a necessary condition on the solution is that it should be order-two differentiable. However, it was soon realized that this condition is too restrictive. Understanding the exact nature of the solutions of the vibrating string equation led to a fierce controversy which involved some of the most brilliant mathematicians of the eighteenth century. Among them are Euler,\index{Euler, Leonhard (1707--1783)} d'Alembert, Lagrange\index{Lagrange, Joseph-Louis (1736--1813)} and Daniel Bernoulli. There is a large amount of primary literature  concerned with this debate, including several memoirs by each of these mathematicians and the correspondence between them. Let us recall a few points of the history of that controversy.

   In the memoir \cite{Al_recherches}, d'Alembert\index{Alembert@d'Alembert, Jean le Rond (1717--1783)} wrote that the general solution to the wave equation is a function of the form
\[y(x,t)=\frac{1}{2}\left(\phi(x+\alpha t)+\phi(x-\alpha t)\right),\] where $\phi$ is an ``arbitrary" periodic function whose period is the double of  the length of the string. The problem was to give a meaning to the adjective ``arbitrary."

At the beginning of his memoir \cite{Al_recherches}, d'Alembert\index{Alembert@d'Alembert, Jean le Rond (1717--1783)} declares that he will show that the problem admits infinitely many other solutions than the usual one represented by the sine curve (which he calls \emph{compagne de la cyclo\"\i de allong\'ee}). But from his point of view (like from Euler's one) the only acceptable functions were those given by a formula (functions which, as we recall, were termed ``analytic" by Euler).\index{analytic function (in the sense of Euler)}\index{function!analytic (in the sense of Euler)} The reason is that it was considered that the powerful methods of analysis can be applied only to such functions.  

 Euler\index{Euler, Leonhard (1707--1783)} published an article in the next volume of the Memoirs (1748) \cite{Euler_vib-corde} in which he 
gave an exposition of d'Alembert's\index{Alembert@d'Alembert, Jean le Rond (1717--1783)} results but where he expressed a different point of view on the nature of the solution of the wave equation.\footnote{Euler  published a Latin and a French version of his memoir, which appeared in the years 1749 and 1748 respectively. (The title of the French version, \emph{Sur la vibration des cordes, traduit du latin}, although it was published first, shows that it was written after the Latin one.)}  He claimed that a solution is not necessarily given by a formula, but that it might be  ``discontinuous"  in the sense that it could be a concatenation of functions defined on smaller intervals on which the restriction of the function is defined by formulae. We already mentioned this notion of ``discontinuity" in \S \ref{s:complex} of the present paper. His assertion was supported by physical evidence, more precisely, by the fact that the initial form of a string, in a musical instrument that is pinched in the usuel manner, is a concatenation of two segments with a corner at their intersection. More than that, Euler\index{Euler, Leonhard (1707--1783)} pointed out that one may give an arbitrary initial form to the string, and therefore the solution may be arbitrary.  Euler's paper introduced some doubts concerning the assertion made by d'Alembert\index{Alembert@d'Alembert, Jean le Rond (1717--1783)} that the solution must be twice differentiable and given by a formula.
 D'Alembert,\index{Alembert@d'Alembert, Jean le Rond (1717--1783)} who disagreed with Euler's claim, published the following year a memoir in which he confirmed his initial ideas. The rest of the controversy on the notion of function is very interesting and there are several articles on this subject. We recommend in particular the introduction, by Youschkevitch and Taton, of Volume V of Series IV A of Euler's \emph{Opera omnia}  containing Euler's correspondence with  Clairaut, d'Alembert and Lagrange \cite{YT}.

The difficulty of defining a general notion of function is never too much emphasized. We mention in this respect that in  1787, that is, four years after Euler's\index{Euler, Leonhard (1707--1783)} death, the Academy of Sciences of Saint-Petersburg proposed, as a competition question, to write an essay on the nature of an arbitrary function.\footnote{\emph{Histoire de l'Acad\'emie Imp\'eriale des Sciences, ann\'ee 1787}, p. 4.} The prize went to the Alsatian mathematician 
 Louis-Fran\c cois-Antoine Arbogast (1759--1803),\index{Arbogast, Louis-Fran\c cois-Antoine (1759--1803)}  who, in his  \emph{M\'emoire sur la nature des fonctions arbitraires qui entrent dans les int\'egrales des \'equations aux diff\'erentielles partielles} (Memoir on the nature of arbitrary functions that appear in the integrals of partial differential equations), \cite{Arbogast} (1791), adopted Euler's point of view: he accepted discontinuous functions in the sense Euler defined them, as solutions of partial differential equations. It is interesting to note that in the description of that problem, the Academy starts with the physical problem of vibrating strings: 
 \begin{quote}\small 
 The problem of the vibrating strings is without doubt one of the most famous problems of applied mathematics. The most celebrated geometers of our time, who solved it, have argued on the legitimacy of their solution, without ever being able to convince each other. It is not that it is difficult to reduce the problem itself to pure analysis, but as it has given the first occasion to treat three-variable differential equations which give, by integrating them, arbitrary and varying functions, the important question which divided the points of view of these great men is whether these functions are entirely arbitrary, whether they represent all the arbitrary curves and surfaces, formed by a voluntary motion of the hand, or whether they include only those that are comprised under an algebraic or transcendental equation. Besides the fact that on that decision depends the way of terminating the dispute on vibrating strings, the same question on the nature of arbitrary functions re-emerges each time an arbitrary problem leads to differential equations with three or more variables: this happens even very often, not only when we treat subjects of sublime mechanics, but most of all in the theory of fluid motion: in such a way that one cannot rigorously sustain that such a problem has been solved before setting precisely the nature of of arbitrary functions. The Academy invites then all the geometers to decide:
 
 \emph{Whether arbitrary functions, to which we are led by integrating equations with one or several variables, represent arbitrary curves or surfaces, either algebraic or transcendental, or mechanical, discontinuous or produced by a voluntary motion of the hand; or whether these functions only contain continuous curves represented by an algebraic or transcendental equation.}\footnote{Le probl\`eme des cordes vibrantes est sans contredit un des plus fameux probl\`emes de la math\'ematique appliqu\'ee. Les plus c\'el\`ebres g\'eom\`etres de notre temps, qui l'ont r\'esolu, se sont disput\'es sur la l\'egitimit\'e de leurs solutions, sans avoir jamais pu se convaincre l'un l'autre. Ce n'est pas que le probl\`eme en lui-m\^eme ne soit pas facilement r\'eduit \`a l'analyse pure ; mais comme il a \'et\'e le premier qui ait donn\'e occasion de traiter des \'equations diff\'erentielles \`a trois variables, par l'int\'egration desquelles on parvient \`a des fonctions arbitraires et variables, la question importante qui partagea les avis de ces grands hommes fut : si ces fonctions sont enti\`erement arbitraires ? si elles repr\'esentent toutes les courbes et surfaces quelconques, form\'ees par un mouvement volontaire de la main ? ou si elles ne renferment que celles qui sont comprises sous une \'equation soit alg\'ebrique soit transcendante ? Outre que c'est de cette d\'ecision que d\'epend le moyen de terminer cette dispute sur les cordes vibrantes, la m\^eme question sur la nature des fonctions arbitraires rena\^\i t toutes les fois qu'un probl\`eme quelconque conduit \`a des \'equations diff\'erentielles \`a trois et plusieurs variables : ce qui arrive m\^eme bien souvent, non seulement lorsqu'on traite des sujets de la m\'ecanique sublime, mais surtout des mouvements des fluides : de sorte qu'on ne saurait soutenir rigoureusement qu'un pareil probl\`eme ait \'et\'e r\'esolu, avant qu'on ait fix\'e exactement la nature des fonctions arbitraires. L'Acad\'emie invite donc tous les g\'eom\`etres de d\'ecider:
 
 \emph{Si les fonctions arbitraires, auxquelles on parvient par l'interm\'ediaire des \'equations \`a trois ou plusieurs variables, repr\'esentent des courbes ou surfaces quelconques, soit alg\'ebriques soit transcendantes, soit m\'ecaniques, discontinues, ou produites par un mouvement volontaire de la main ; ou si ces fonctions renferment seulement des courbes continues repr\'esent\'ees par une \'equation alg\'ebrique ou transcendante ?}}
 \end{quote}

  We now come to the problem of trigonometric series.\index{trigonometric series}
  
    Brook Taylor, in his 1713 memoir entitled \emph{De motu nervi tensi} (On the motion of a tense string) \cite{Taylor1713} (cf. also his \emph{Methodus incrementorum directa et inversa} (Direct and Indirect Methods of Incrementation), \cite{Taylor1715} (first edition 1715), showed that the vibration problem admits as a solution the sine and cosine functions.  For several reasons which we shall mention below, it was tempting to conjecture then that the general solution of the problem is obtained by taking an infinite sum of trigonometric functions. This was done by Daniel Bernoulli\index{Bernoulli, Daniel (1700--1782)}  (1700--1782).

In 1753, Bernoulli\index{Bernoulli, Daniel (1700--1782)} wrote a memoir on the vibration of strings.\index{string vibration}\index{vibrating string} Bernoulli had already thought about this question for several years.  In his approach to it, like in the other physical problems he considered, Bernoulli\index{Bernoulli, Daniel (1700--1782)} was an adept of Leibniz' calculus, rather than Euler's geometric methods (which were adopted by d'Alembert\index{Alembert@d'Alembert, Jean le Rond (1717--1783)}).  As a physicist, the mathematical notion of function was not a central theme in his research, and from his point of view, the function representing the solution of the question was simply identified with the shape of the vibrating string. While
   Taylor had considered each trigonometric solution individually, that is, he noticed that functions of the form
    \[y(x,t)=\sin \frac{n\pi x}{l}\cos \frac{n\pi \alpha t}{l}\]
    are solutions of  the wave equation\index{wave equation}\index{equation!wave}, Bernoulli stated that the general solution was an infinite sum of such functions.  Thus, he added to the debate the question of the convergence of trigonometric series.\index{trigonometric series}      
In the meanwhile, d'Alembert\index{Alembert@d'Alembert, Jean le Rond (1717--1783)} published a first supplement to his memoir \emph{Sur les vibrations des cordes sonores} (On the vibration of sonorous strings) \cite{dAL_Opuscules2} in which, referring to Daniel Bernoulli's work, he  writes: 
 \begin{quote}\small 
  The question is not to \emph{conjecture}, but to \emph{prove}, and it would be dangerous (although, to tell the truth, this  misfortune is unlikely to happen) that this kind of proof which is so odd enters into geometry. The only thing which seems surprising is that such reasonings are used in way of a proof by a famous mathematician.\footnote{Il ne s'agit pas de \emph{conjecturer}, mais de \emph{d\'emontrer}, et il serait dangereux (quoi qu'\`a la v\'erit\'e ce malheur soit peu \`a craindre) qu'un genre de d\'emonstration si singulier s'introduis\^it en G\'eom\'etrie. Ce qui pourra seulement para\^itre surprenant, c'est que de pareils raisonnements soient employ\'es comme d\'emonstratifs par un math\'ematicien c\'el\`ebre.}
    \end{quote}

After the publication of Bernoulli's memoir, Euler\index{Euler, Leonhard (1707--1783)}  wrote a new memoir in which he  generalizes Bernoulli's result, \emph{Remarques sur les m\'emoires pr\'ec\'edents de M. Bernoulli} (Remarks on the preceding memoirs by Mr. Bernoulli)  \cite{E213} (1755). He also confirms his own intuition that a solution of the vibrating string equation may be an arbitrary function. Today, we know that, in some sense, the solution he proposed is identical to that of Bernoulli, but the relation between trigonometric series\index{trigonometric series} and arbitrary functions was not yet discovered. 
 
                       In his memoir, Euler\index{Euler, Leonhard (1707--1783)} starts by declaring that, without any doubt, Bernoulli  developed the theory of formation of sound\index{sound!theory of}  infinitely better than any other scientist before him, that his predecessors stopped at the mechanical determination of the motion of a tight string without any thorough investigation of the nature of sound,\index{sound!theory of}  and that it was still not understood how a single string can emit several sounds at  the same time. He then expresses his doubts about the fact that Bernoulli's infinite series of sines could be the general solution of the problem. He writes that it is impossible for the curve made by  a vibrating string to be constituted by an infinite number of trocho\"\i ds (which is the name he used for the sine curves). He declares that there are infinitely many curves that are not included in that solution. 
                       
                       The solution of this problem was given by Fourier\index{Fourier, Joseph (1768--1830)} and completed by Dirichlet\index{Dirichlet, Johann Peter Gustav Lejeune (1805--1859)} and Riemann, in the following century, as we shall discuss below. In the same paper, Euler\index{Euler, Leonhard (1707--1783)} insists on the fact that the general solution of the equation of the vibrating string cannot be given by a formula. He  mentions his own conflict with d'Alembert,\index{Alembert@d'Alembert, Jean le Rond (1717--1783)} saying that he wishes very much that the latter explains why he is mistaken. Based on partial differential calculus, Euler gives a new explanation of the fact that by varying the initial shape of a string, any function becomes admissible as a solution of the problem.

Between November 1, 1759 and the end of the same year, Euler\index{Euler, Leonhard (1707--1783)} presented three memoirs,  \cite{LE_2}, \cite{LE_3} and \cite{LE_4}, on the propagation of sound.\index{sound!theory of}   In these memoirs,
he studies respectively the propagation in one, two and three dimensions. The differential equations that describe this propagation are the same as those which describe the vibration of strings.\index{string vibration}\index{vibrating string} Euler mentions the limitations of the works of ``Taylor, Bernoulli and some others." 
Despite the fact that the debate on  the vibrating string had already lasted many years, the relation between the scientists working on that subject was still very tense.
 
In a later memoir, titled \emph{M\'emoire   sur   les   vibrations   des   cordes   d'une   \'epaisseur   in\'egale} (Memoir of the vibration of strings of uneven width) \cite{Ber1767} (1767), Bernoulli gave an additional reason for the use of an infinite sum (p. 283): 
 \begin{quote}\small
 When a string makes several  vibrations at the same time, of whatever number, and in whatever order, the absolute curvature will always be expressed by the general equation 
  \[y=\alpha\sin\frac{x}{l}\pi + \beta\sin\frac{2x}{l}\pi+ \gamma\sin\frac{3x}{l}\pi+ \mathrm{etc.}\]
and since the number of arbitrary coefficients is infinite, one can make the curve pass by whatever number of points of positions that we wish, which indicates that all the curves belong to this case, provided we do not oppose the hypotheses. And it would be opposing them if we do not treat the quantities $y$, $dy$ and $ddy$ as infinitely smaller, at every point of the curve, than the quantities $x$, $dx$ and $\frac{dx^2}{l}$.
 \footnote{Lorsque la corde fait \`a la fois plusieurs esp\`eces de vibration, quel qu'en soit le nombre, et de quelque ordre qu'elles soient, la courbure absolue sera toujours exprim\'ee par cette \'equation g\'en\'erale
 \[y=\alpha\sin\frac{x}{l}\pi + \beta\sin\frac{2x}{l}\pi+ \gamma\sin\frac{3x}{l}\pi+ \mathrm{etc.}\]
et comme le nombre des coefficients arbitraires est infini, on peut faire passer la courbe par tant de points donn\'es de position qu'on voudra, ce qui marque que toutes les courbes se trouvent dans ce cas, pourvu qu'on ne fasse pas violence aux hypoth\`eses. Et ce serait leur faire violence, si on ne faisait pas les quantit\'es $y$, $dy$ et $ddy$ comme infiniment plus petites dans tous les points de la courbe, que les quantit\'es $x$, $dx$ et $\frac{dx^2}{l}$.}
 \end{quote}

 The interested reader may skim the volume of Euler's collected works containing the correspondence between Euler and Lagrange, \cite{EL}, and the volume containing the correspondence between Lagrange and d'Alembert, \cite{OL13}, not only in order to understand more deeply this multi-faced controversy,  but also in order to feel the cultural and scientific atmosphere in Europe during that period. Let us quote, as examples, two excerpts related to the discussion around the solution of the wave equation. 
In a letter to Lagrange, dated Octobre 1759, Euler writes: 
  \begin{quote}\small
  I was pleased to learn that you agree with my solution relative to the vibrating strings, which d'Alembert\index{Alembert@d'Alembert, Jean le Rond (1717--1783)} tried hard to refute using various sophisms, for the only reason that he did not propose it himself. He announced that he will publish an overwhelming proof of it; I don't know whether he did it. He thinks he will be able to impress people by his half-scholar eloquence. I doubt that he can seriously play such a role, unless he is profoundly blinded by self-esteem. He wanted to insert in our \emph{Memoirs}, not a proof, but a simple declaration according to which my solution was very deficient. On my side, I proposed a new proof which has all the required rigor.\footnote{J'ai appris avec plaisir que vous approuviez ma solution relative aux cordes vibrantes, que d'Alembert\index{Alembert@d'Alembert, Jean le Rond (1717--1783)} s'est efforc\'e de r\'efuter par divers sophismes, et ceci pour l'unique raison qu'il ne l'a pas propos\'ee lui-m\^eme. Il a annonc\'e qu'il en publierait une accablante r\'efutation ; j'ignore s'il l'a fait. Il croit qu'il pourra jeter de la poudre aux yeux avec son \'eloquence de demi-savant. Je doute qu'il joue ce r\^ole s\'erieusement, \`a moins qu'il ne soit profond\'ement aveugl\'e par l'amour-propre. Il a voulu ins\'erer dans nos M\'emoires non une d\'emonstration, mais une simple d\'eclaration suivant laquelle ma solution \'etait tr\`es d\'efectueuse ; pour ma part, j'ai propos\'e une nouvelle d\'emonstration poss\'edant toute la rigueur voulue.}
  \end{quote}

In a letter to  d'Alembert, dated March 20, 1765, Lagrange writes: 
\begin{quote}\small
Concerning [our discussion] on vibrating strings, it is now reduced to a point which, it seems to me, escapes any analysis. Moreover, I found, by a completely direct way, that if we admit in the intial figure the conditions that you ask, the solution reduces to the one of Mr. Bernoulli, namely,  $y=\alpha \sin \frac{\pi x}{a} + \beta \sin \frac{2 \pi x}{a} + ...$, and it is difficult for me to believe that this is the only one that can be found in nature. Besides, the phenomena of sound propagation can be explained only if we admit discontinuous functions, as I proved in my second dissertation.\footnote{\`a l'\'egard de [notre discussion] sur les cordes vibrantes, elle est maintenant r\'eduite \`a un point qui \'echappe, ce me semble, \`a l'Analyse. Au reste, j'ai trouv\'e par une voie tout \`a fait directe qu'en admettant dans la figure initiale les conditions que vous y exigez, la solution se r\'eduit \`a celle de M. Bernoulli, savoir : $y=\alpha \sin \frac{\pi x}{a} + \beta \sin \frac{2 \pi x}{a} + ...$, et j'ai peine \`a croire que celle-ci soit la seule qui puisse avoir lieu dans la nature. D'ailleurs, les ph\'enom\`enes de la propagation du son ne peuvent s'expliquer qu'en admettant les fonctions discontinues, comme je l'ai prouv\'e dans ma seconde dissertation.}
 \end{quote}

   Let us also quote Nicolaus Fuss,\index{Fuss, Nicolaus (1755--1826)} the famous biographer of Euler,\footnote{Nicolaus Fuss (1755--1826) was first hired as Euler's secretary, then he became successively his favorite student, his closest friend, his collaborator and colleague at the Saint Petersburg Academy of Sciences, and eventually his grandson-in-law (the husband of Euler's grand-daughter Albertine).} from his  \emph{\'Eloge} \cite{Fuss}: 
   \begin{quote}\small
   The controversy between  Messrs.  Euler, d'Alembert\index{Alembert@d'Alembert, Jean le Rond (1717--1783)} \& Bernoulli regarding the motion of the vibrating strings can be of interest only to professional geometers. Mr. D. Bernoulli, who was the first to develop the physical part which concerns the production of sound\index{sound!theory of}  generated by this motion, thought that Taylor's solution was sufficient to explain it.  Messrs. Euler\index{Euler, Leonhard (1707--1783)} and  d'Alembert,\index{Alembert@d'Alembert, Jean le Rond (1717--1783)} who had exhausted, in this difficult matter, everything exquisite and profound that an analytic mind may possess, showed that the solution of Mr. Bernoulli, extracted from Taylor's Trochoids, is not general, and that it is even deficient. This controversy, which lasted a long time, with all the consideration that such famous men owe to each other, gave rise to a   quantity of excellent memoirs; it really ended only at the death of Bernoulli.\footnote{La controverse entre MM. Euler, d'Alembert\index{Alembert@d'Alembert, Jean le Rond (1717--1783)} \& Bernoulli au sujet du mouvement des cordes vibrantes ne peut int\'eresser proprement que les G\'eom\`etres de profession. M. D. Bernoulli, qui fut le premier \`a en d\'evelopper la partie physique qui regarde la formation du son engendr\'e par ce mouvement, crut la solution de Taylor suffisante de l'expliquer. MM. Euler\index{Euler, Leonhard (1707--1783)} et d'Alembert, qui avaient \'epuis\'e, dans cette mati\`ere difficile, tout ce que l'esprit analytique a de sublime \& de profond, firent voir que la solution de M. Bernoulli, tir\'ee des Trocho\"\i des Tayloriennes, n'est pas g\'en\'erale, qu'elle est m\^eme insuffisante. Cette controverse qui a \'et\'e continu\'ee longtemps, avec tous les \'egards que des hommes aussi illustres se doivent mutuellement, a donn\'e naissance \`a quantit\'e d'excellents m\'emoires ; elle n'a fini proprement qu'\`a la mort de M. Bernoulli.}  \end{quote}     
   
D'Alembert\index{Alembert@d'Alembert, Jean le Rond (1717--1783)} eventually accepted functions that are discontinuous (in the sense of Euler) as solutions of partial differential equations; cf. his 1780 memoir entitled \emph{Sur les fonctions discontinues} (On discontinuous functions), published in \cite{Alembert-OP} (t. VIII, \S VI) in which he formulates a \emph{R\`egle sur les fonctions discontinues qui peuvent entrer dans l'int\'egration des \'equations aux d\'eriv\'ees partielles} (Rule on discontinuous functions that may enter into the integration of partial differential equations). We refer the interested reader to the papers \cite{YT}, \cite{Thiele}, \cite{Y} and \cite{Jehel} for more on the history of the subject.

The confirmation of Bernoulli's conjecture followed from Fourier's\index{Fourier, Joseph (1768--1830)} manuscript \emph{Th\'eorie de la propagation de la chaleur dans les solides} (Theory of heat propagation in solids),\footnote{This memoir, read to the Academy 1807, was never published, until it was edited with comments by Grattan-Guinness, see \cite{Grattan-F}.} \cite{Fourier1807} read to the Academy in 1807, that is, twenty-five years after Bernoulli's death. The manuscript carries the subtitle: ``M\'emoire sur la propagation de la chaleur avec notes s\'epar\'ees sur cette propagation -- sur la convergence des s\'eries $\displaystyle \sin x-\frac{1}{2}\sin 2x + \frac{1}{3}\sin 3x$ \&c." (Memoir on the propagation of heat,  with separate notes on that propagation -- on the convergence of the series  $\displaystyle \sin x-\frac{1}{2}\sin 2x + \frac{1}{3}\sin 3x$ etc.).  Let us quote an excerpt (\cite{Grattan-F} p. 183): 
   \begin{quote}\small
It follows from my researches on this object that the arbitrary functions, even discontinuous, can always be represented by the sine or cosine expansions of multiple arcs, and that the integrals  which contain these developments are precisely as general as those where arbitrary functions of multiple arcs enter. A conclusion that the celebrated Euler\index{Euler, Leonhard (1707--1783)} has always rejected.\footnote{Il r\'esulte de
mes recherches  sur cet objet que  les fonctions  arbitraires m\^eme discontinues peuvent toujours
\^etre   repr\'esent\'ees   par   les   d\'eveloppements   en   sinus   ou   cosinus   d'arcs   multiples,   et   que   les
int\'egrales   qui   contiennent   ces   d\'eveloppements   sont   pr\'ecis\'ement   aussi   g\'en\'erales   que   celles   o\`u
entrent   les   fonctions   arbitraires   d'arcs   multiples. Conclusion que le c\'el\`ebre Euler a toujours repouss\'ee.}
\end{quote}

In 1811, the Paris \emph{Acad\'emie des sciences} proposed a competition whose title was: \emph{Donner la th\'eorie math\'ematique des lois de la propagation de la chaleur et comparer les r\'esultats de cette th\'eorie \`a des exp\'eriences exactes} (To give the mathematical theory of the laws of propagation of heat and to compare the results of this theory with exact experiences).  Fourier\index{Fourier, Joseph (1768--1830)} submitted for the prize a very extensive work which included his 1807 manuscript. The jury of the competition consisted of Lagrange,\index{Lagrange, Joseph-Louis (1736--1813)} Laplace, Maus, Ha\"uy and Legendre. Darboux,\index{Darboux, Gaston (1842--1917)} in his review \cite{Darboux1888}, quotes part of the report on the work of Fourier:

      \begin{quote}\small
      
      This piece  contains genuine differential equations of heat transmission, either in the interior of bodies, or at their surface. And what is new in the subject, added to its importance, has led the Class to crown this treatise, while noting however that the manner with which the author arrives at his equation is not exempt of difficulties, and that his analysis, to integrate them, still leaves something to be desired, either relative to the generality, or even from the point of view of rigor.\footnote{Cette pi\`ece renferme de v\'eritables \'equations diff\'erentielles de la transmission de la chaleur, soit \`a l'int\'erieur des corps, soit \`a leur surface ; et la nouveaut\'e du sujet, jointe \`a son importance, a d\'etermin\'e la Classe \`a couronner cet Ouvrage, en observant cependant que la mani\`ere dont l'Auteur parvient \`a ses \'equations n'est  pas exempte de  difficult\'es,  et que  son  analyse,  pour les  int\'egrer, laisse  encore
quelque   chose   \`a   d\'esirer,   soit   relativement   \`a   la   g\'en\'eralit\'e,   soit   m\^eme   du   c\^ot\'e   de   la   rigueur.}
\end{quote}

    The sum of Fourier's\index{Fourier, Joseph (1768--1830)} work on the propagation of heat was collected in his masterpiece, \emph{Th\'eorie analytique de la chaleur}\index{Fourier!Th\'eorie analytique de la chaleur} (Analytic theory of heat) \cite{Fourier}, published in 1822.  The following  result is stated at the end of Chapter III of this memoir, as a summary of what has been done (art. 235, p. 258):
    \begin{quote}\small 
    It follows from all that was proved in this section, concerning the series expansion of trigonometric functions, that if we propose a function $fx$, whose value is represented on a given interval, from $x=0$ to $x=X$, by the ordinate of a curves line drawn arbitrarily, one can always expand this function as a series which will contain only the sines, or the cosines, or the sines and the cosines of multiple arcs, or only the cosines of odd multiples.\footnote{Il r\'esulte de tout ce qui a \'et\'e d\'emontr\'e dans cette section, concernant le d\'eveloppement des fonctions en s\'eries trigonom\'etriques, que si l'on propose une fonction $fx$, dont la valeur est repr\'esent\'ee dans un intervalle d\'etermin\'e, depuis $x=0$ jusqu'\`a $x=X$, par l'ordonn\'ee d'une ligne courbe trac\'ee arbitrairement on pourra toujours d\'evelopper cette fonction en une s\'erie qui ne contiendra que les sinus, ou les cosinus, ou les sinus et les cosinus des arcs multiples, ou les seuls cosinus des multiples impairs.} 
    \end{quote}
    Trigonometric functions were essential for the solution of the heat equation, as they were for the wave equation at the time of Bernoulli. What is important for our topic here is that trigonometric series became an essential tool in the field of complex analysis, independently of the heat flow.
    Fourier writes in the same section     (art. 235 p. 258):
   \begin{quote}\small We cannot entirely solve the fundamental questions of the theory of heat without reducing to this form the functions that represent the initial state of temperatures. 
   
   These trigonometric series, ordered according to the cosines or sines of the multiples of the arc, pertain to elementary analysis, like the series whose terms contain successive powers of the variable. The coefficients of the trigonometric series are definite areas, and those of power series are fractions given by differentiation, in which one also attributes to the variable a definite value.\footnote{On ne peut r\'esoudre enti\`erement les questions fondamentales de la th\'eorie de la chaleur, sans r\'eduire \`a cette forme les fonctions qui repr\'esentent l'\'etat initial des temp\'eratures. 

    Ces s\'eries trigonom\'etriques, ordonn\'ees selon les cosinus ou les sinus des multiples de l'arc, appartiennent \`a l'analyse \'el\'ementaire, comme les s\'eries dont les termes contiennent les puissances successives de la variable. Les coefficients des s\'eries trigonom\'etriques sont des aires d\'efinies, et ceux des s\'eries de puissances sont des fonctions donn\'ees par la diff\'erenciation, et dans lesquelles on attribue aussi \`a la variable une valeur d\'efinie.}    \end{quote}
    
Fourier\index{Fourier, Joseph (1768--1830)} then summarizes several properties of these series, including the integral formulae for the coefficients, and he also states the following:
    \begin{quote}\small
    The series, ordered according to the cosines or the sines of the multiple arcs, are always convergent, that is, when we give to the variable an arbitrary non-imaginary value, the sum of the terms converges more and more to a unique and fixed limit, which is the value of the expanded function.\footnote{Les s\'eries ordonn\'ees selon les cosinus ou les sinus des arcs multiples sont toujours convergentes, c'est-\`a-dire qu'en donnant \`a la variable une valeur quelconque non imaginaire, la somme des termes converge de plus en plus vers une seule limite fixe, qui est la valeur de la fonction d\'evelopp\'ee.}
    \end{quote} 
    
    It is interesting to recall that Fourier,\index{Fourier, Joseph (1768--1830)} in his  \emph{Th\'eorie analytique de la chaleur}\index{Fourier!Th\'eorie analytique de la chaleur} quotes  Archimedes,\index{Archimedes (c. 287 B.C.--c. 212 B.C.)} Galileo\index{Galilei, Galileo (1564--1642)} and Newton,\index{Newton, Isaac (1643--1727)} the three scientists mentioned by  Riemann in the last part of his Habilitation lecture, \emph{\"Uber die Hypothesen, welche der Geometrie zu Grunde liegen}\index{Riemann! habilitation lecture}\index{habilitation lecture!Riemann}. Fourier writes (\cite{Fourier}, p. i--ii):
\begin{quote}
The knowledge that the most ancient people could have acquired in rational mechanics did not reach us, and the history of that science, if we except the first theorems on harmony, does not go back further than the the discoveries of Archimedes.\index{Archimedes (c. 287 B.C.--c. 212 B.C.)} This great geometer explained the mathematical principles of the equilibrium of solids and of fluids. About eighteen centuries passed before Galileo,\index{Galilei, Galileo (1564--1642)} the first inventor of the dynamical theories, discovered the laws of motions of massive bodies. Newton\index{Descartes, Ren\'e (1596--1650)} encompassed in that new science all the system of the universe.\footnote{Les connaissances que les plus anciens peuples avaient pu acqu\'erir dans la m\'ecanique rationnelle ne nous sont pas parvenues, et l'histoire de cette science, si l'on excepte les premiers th\'eor\`emes sur l'harmonie, ne remonte point au-del\`a des d\'ecouvertes d'Archim\`ede. Ce grand g\'eom\`etre expliqua les principes math\'ematiques de l'\'equilibre des solides et des fluides. Il s'\'ecoula environ dix-huit si\`ecles avant que Galil\'ee, premier inventeur des th\'eories dynamiques, d\'ecouvrit les lois du mouvement des corps graves. Newton\index{Newton, Isaac (1643--1727)} embrassa dans cette science nouvelle tout le syst\`eme de l'univers.}
\end{quote}

Riemann claims in his memoir on trigonometric series that the question of finding conditions under which a function can be represented by a trigonometric series was completely settled in the work of Dirichlet,\index{Dirichlet, Peter Gustav Lejeune (1805--1859)} ``for all cases that one encounters in nature. [...] The questions to which Dirichlet's researches do not apply do not occur in nature." We quote  Dirichlet, from his 1829 memoir \emph{Sur la convergence des s\'eries trigonom\'etriques qui servent \`a repr\'esenter une fonction arbitraire entre des limites donn\'ees} (On the convergence of trigonometric series used to represent an arbitrary function between two given bounds) \cite{Dirichlet1829} (1829), in which he gives a solution to the convergence problem, and where he also mentions the work of Cauchy on the same problem. His paper starts as follows:

\begin{quote}\small
The series of sines and cosines, by means of which one can represent an arbitrary function on a given interval, enjoy among other remarkable properties the one of being convergent. This property  has not escaped the attention of the famous geometer who opened a new field for the applications of analysis, introducing there the way of expressing the arbitrary functions which are our  subject here. 
It is contained in the memoir containing his first researches on heat. But, to my knowledge, nobody up to now gave a general proof. I know only on this subject a work due to Mr. Cauchy and which is part of his \emph{M\'emoires de l'Acad\'emie des sciences de Paris pour l'ann\'ee 1823}. The author of this work confesses himself that his proof fails for certain functions for which convergence is nevertheless indisputable.  A careful examination of that memoir led me to the belief that the proof which is presented there is insufficient even for the cases to which the authors thinks it applies.\footnote{Les s\'eries de sinus et de cosinus, au moyen desquelles on peut repr\'esenter une fonction arbitraire dans un intervalle donn\'e, jouissent entre autres propri\'et\'es remarquables aussi de celles d'\^etre convergentes. Cette propri\'et\'e n'avait pas \'echapp\'e au g\'eom\`etre illustre qui a ouvert une nouvelle carri\`ere aux applications de l'analyse, en y introduisant la mani\`ere d'exprimer les fonctions arbitraires dont il est question ; elle se trouve \'enonc\'ee dans le M\'emoire qui contient ses premi\`eres recherches sur la chaleur. Mais personne, que je sache, n'en a donn\'e jusqu'\`a pr\'esent une d\'emonstration g\'en\'erale. 
Je ne connais sur cet objet qu'un travail d\^u \`a M. Cauchy et qui fait partie des M\'emoires de l'Acad\'emie des sciences de Paris pour l'ann\'ee 1823. L'auteur de ce travail avoue lui-m\^eme que sa d\'emonstration tombe en d\'efaut pour certaines fonctions pour lesquelles la convergence est pourtant incontestable. Un examen attentif du m\'emoire cit\'e m'a port\'e \`a croire que la d\'emonstration qui y est expos\'ee n'est pas m\^eme suffisante pour les cas auxquels l'auteur la croit applicable.}

\end{quote}

Riemann, after quoting the work of Dirichlet,  gives two reasons for investigating the cases to which Dirichlet's methods do not apply directly. The first reason is that sorting out these questions will bring more clarity and precision to the principles of infinitesimal calculus. The second reason is that Fourier series\index{Fourier series} will be useful not only in physics, but also in number theory.  Riemann says that in this field, it is precisely the functions which Dirichlet\index{Dirichlet, Johann Peter Gustav Lejeune (1805--1859)}  did not consider that seem to be the most important.

The so-called Dirichlet conditions for a function defined on the interval $[0,2\pi]$ to  have a Fourier trigonometric expansion is now classical.  Picard, in \cite{Picard-developpement} (p. 8) writes that Dirichlet's memoir on Fourier series remained a model of rigor. We conclude this section by quoting Riemann. 
Talking about Dirichlet's\index{Dirichlet, Johann Peter Gustav Lejeune (1805--1859)} work on trigonometric series, he writes: 
 \begin{quote}\small
This work of Dirichlet\index{Dirichlet, Johann Peter Gustav Lejeune (1805--1859)} gave a firm foundation to a great number of important analytic researches. Highlighting a point on which Euler\index{Euler, Leonhard (1707--1783)} was mistaken, he succeeded in clearing out a question that had occupied so many eminent geometers for more than seventy years. 
\end{quote}

 Riemann, in the same memoir, developed his integration theory  in order to build a general theory for Fourier series,\index{Fourier series} in particular for functions which have an infinite number of discontinuity  points. This is the subject of the next section.

\section{Integration}\label{s:integration}

The second part of Riemann's memoir on trigonometric functions \cite{Riemann-Trigo} carries the title ``On the notion of definite integral and on the scope of its applicability."\index{Riemann integral}\index{integral!Riemann} The relation between integration and trigonometric series is based on Fourier's formulae which give the coefficients of a trigonometric series in the form of integrals.
Riemann starts the second part of the memoir by formulating a question: ``First of all, what do we mean by 
\[\int_a^b f(x)dx?"\]
The rest of the memoir is the answer to this question.

We explained at length that one of the fundamental questions in eighteenth century mathematics was ``\emph{What is a function}," and how this question led to a celebrated controversy.  Riemann had to deal with the same question in his memoir on trigonometric series, more than a century after the controversy started, and he gave it the definitive answer. The problem in Riemann's memoir was addressed in a new context, namely, his integration theory, which was developed in a few pages at the end of that memoir. More particularly, the  question became in that context:  ``\emph{What are the functions that can be integrated}?" and in particular, whether the known functions were sufficient for the theory that became known as the Riemann integral\index{Riemann integral} or whether a new class of functions was needed.

We recall that since Newton and Leibniz, and passing by Euler, integration was defined as an anti-derivarive.
Cauchy started an approach to integrals as limits of sums associated to partitions of the interval of definition, that is, sums of the form
\[\sum_{k=1}^\infty f(x_{k-1})(x_k-x_{k-1}),\]
cf. e.g. Cauchy's\index{Cauchy, Augustin-Louis (1789--1857)} \emph{R\'esum\'e des leçons donn\'ees \`a l'\'Ecole Royale Polytechnique sur le calcul infinit\'esimal} (Summary of lectures  on infinitesimal calculus given at the \emph{\'Ecole Royale Polytechnique}) (1823)  \cite{Cauchy-Resume}. In Cauchy's setting, the limit always exists because he\index{Cauchy integral}\index{integral!Cauchy} considered only continuous functions. It was soon realized that the definition may apply to more general functions.  Dirichlet, in his work on trigonometric functions, used Cauchy's theory applied to discontinuous functions. 
 Riemann states in his memoir that Cauchy's integration theory involves some random definitions which cannot make it a universal theory.

Riemann answered the question of how far discontinuity is allowed. He was led to the most general functions, which he termed ``integrable."   In  \S VII, VIII and IX of his memoir, he applies his new  integration theory to the problem of representing arbitrary functions by trigonometric series.\index{trigonometric series} The main results are stated in three theorems in \S VIII, and the propositions concerning the representation of   functions by trigonometric series\index{trigonometric series} are contained in \S IX. \S X and XI contain results on the behavior of the coefficients of a trigonometric series. The last sections (\S XII and XIII) concern particular cases, more precisely, cases where the Fourier series\index{Fourier series} is not convergent. 
  
   It is curious that Riemann mentions Cauchy\index{Cauchy, Augustin-Louis (1789--1857)}  several times in this memoir on trigonometric series, but he never refers to him in his dissertation on the theory of functions of a complex variable.

 There is a section on the history of integration in Lebesgue's book \cite{Lebesgue1904}. In particular, Lebesgue summarizes Cauchy's theory,\index{Cauchy integral}\index{integral!Cauchy} as well as an unpublished work of 
Dirichlet on the subject, which reached us through a description of Lipschitz (\cite{Lebesgue1904}  p. 9). Dirichlet's\index{Dirichlet integral}\index{integral!Dirichlet} work applies to functions with an infinite number of discontinuity points, but forming a non-dense subset.  Riemann, using series, constructed functions to which the preceding techniques do not apply and which may still be integrated.  These functions of Riemann do not have a graphical representation.   We are far from Euler's ``arbitrary drawable function" which, indeed, he thought exceeded the power of the calculus (by not being differentiable).

Chapter 2 of Lebesgue's treatise is a survey of the Riemann integral.\index{Riemann integral}\index{integral!Riemann} This theory allows one to prove theorems such as the fact that a uniformly convergent integrable sequence of functions is an integrable function (Lebesgue p. 30), and that a uniformly convergent series of integrable functions may be integrated term by term. 
Lebesgue also mentions the work of Darboux, involving the notions of  upper and lower limits. He then presents   his own geometric theory (as opposed to the analytic theory of Riemann), based on set theory and measure theory. There are more comments on Lebesgues' integration theory in Chapter 8 of the present volume \cite{Papa-Riemann3}.

  To conclude this section, let us mention that Riemann's ideas about the general notion of function in relation with integration theory underwent several developments in the twentieth century (one may think about the difficulties in the introduction of general measurable functions).

    Riemann's memoir on trigonometric series was published 13 years after it was written.
It was translated into French by Darboux and Hou\"el.

It is interesting to note that trigonometric series are used in the proof of the so-called Poincar\'e lemma, a lemma which plays an essential role in the proof of the modern version of the Riemann--Roch theorem which is presented in Chapter 13 of the present volume.

\section{Conclusion} \label{s:conclusion}

In the preceding sections, we reviewed part of the historical origins of Riemann's mathematical works. One should write another article about the roots of his ideas in  physics and philosophy. The intermingling between the old and new ideas of physics and philosophy is yet another subject. In this respect, and since the present book is also about relativity, we quote Kurt G\"odel from his article \emph{A remark about the relationship between relativity theory and idealistic philosophy} \cite{Godel}. Talking about the insight that this theory brings into the nature of time, he writes:
\begin{quote}\small
In short, it seems that one obtains an unequivocal proof for the view of those philosophers who, like Parmenides, Kant and the modern idealists, deny the objectivity of change and consider change as an illusion or an appearance due to our special mode of perception. 
\end{quote}

It would be stating the obvious to say that mathematicians should read the works of mathematicians from the past, not only the recent past, but most of all the founders of the theories they are working on. Yet, very few do it. I would like to conclude the present paper by quoting some pre-eminent mathematicians who expressed themselves on this question. I start with Chebyshev.

We learn from his biographer in \cite{bio} that Chebyshev's\index{Chebyshev, Pafnouti Lvovitch (1821--1894)}  thoroughly studied the works of Euler, Lagrange, Gauss, Abel,\index{Abel, Nils Henrik (1802--1829)} and other pre-eminent mathematicians. The biographer also writes that, in general, Chebyshev was not interested in reading the mathematical works of his contemporaries, considering that spending time on that would prevent him of having original ideas.

On the importance of reading the old masters, we quote again Andr\'e Weil,\index{Weil, Andr\'e (1906--1998)} from his 1978 ICM talk (\cite{Weil-ICM-H} p. 235):
\begin{quote}\small
From my own experience I can testify about the value of suggestions found in Gauss and in Eisenstein.\index{Eisenstein, Ferdinand Gotthold Max (1823--1852)} Kummer's congruences for Bernoulli numbers, after being regarded as little more than a curiosity for many years, have found a new life in the theory of $p$-adic and $L$-functions, and Fermat's ideas on the use of the infinite descent in the study of Diophantine equations of genus one have proved their worth in contemporary work on the same subject.
\end{quote}

Among the more recent mathematicians, I would like to quote again Grothendieck.\index{Grothendieck, Alexandre (1928--2014)} During his apprenticeship, like most of us, Grothendieck was not encouraged to read ancient authors. He writes, in \emph{R\'ecoltes et semailles}\index{Grothendieck!R\'ecoltes et semailles} (Chapter 2, \S\,2.10):
\begin{quote}\small
 In the teaching I received from my elders,  historical references were extremely rare, and I was nurtured, not by reading authors which were slightly ancient, nor even contemporary, but only through communication, face to face or through correspondence with others mathematicians, and starting with those who were older than me.\footnote{Dans l'enseignement que j'ai re\c cu de mes a\^\i n\'es, les r\'ef\'erences historiques \'etaient rarissimes, et j'ai \'et\'e nourri, non par la lecture d'auteurs tant soit peu anciens ni m\^eme contemporains, mais surtout par la communication, de vive voix ou par lettres interpos\'ees, avec d'autres math\'ematiciens, \`a commencer par mes a\^\i n\'es.}
\end{quote}
 
 In the same work, we read (Chapter 2, \S  2.5):
 
 \begin{quote}\small
I personally feel that I belong to a lineage of mathematicians whose spontaneous mission and joy is to constantly construct new houses. [...] I am not strong in history, and if I were asked to give names of mathematicians in that lineage, I can think spontaneously of Galois and Riemann (in the past century) and Hilbert (at the beginning of the present century).\footnote{Je me sens faire partie, quant \`a moi, de la lign\'ee des math\'ematiciens dont la vocation spontan\'ee et la joie est de construire sans cesse des maisons nouvelles. [...] Moi qui ne suis pas fort en histoire, si je devais donner des noms de math\'ematiciens dans cette lign\'ee-l\`a, il me vient spontan\'ement ceux de Galois et de Riemann (au si\`ecle dernier) et celui de Hilbert (au d\'ebut du si\`ecle pr\'esent).}

\end{quote}

Grothendieck's\index{Grothendieck, Alexandre (1928--2014)} attitude towards mathematics and mathematicians changed drastically at the time he decided to quit the mathematical milieu, in 1970, twenty years after he obtained his first job, putting an end to an extraordinarily productive working life and to his relation with his contemporary mathematicians. One thing which is not usually mentioned about him is that his writings, during the period that followed, contain many references to mathematicians of the past, to whom Grothendieck expresses his debt, and among them stands Riemann. 
  In his \emph{R\'ecoles et semailles}\index{Grothendieck!R\'ecoltes et semailles} \cite{RS}, Grothendieck\index{Grothendieck, Alexandre (1928--2014)} writes (Chap. 2, \S 2.5):
 \begin{quote}\small
 Most mathematicians are led to confine themselves in a conceptual framework, in a ``\emph{Universe}," which is fixed once and for all -- essentially, the one they found ``ready-made" at the time they were students. They are like the heirs of a big and completely furnished beautiful house, with its living rooms, kitchens and workshops, its kitchen set and large equipment, with which there is, well, something to cook and to tinker. How this house was progressively constructed,  over the generations, and why and how such and such tool (and not another) was conceived and shaped, why the rooms are fit out in such a manner here, and in another manner there -- these are as many questions as these heirs will never think to ask.  That is the ``Universe," the ``given" in which we must live, that's it!  Something which will seem great (and, most often, we are far from having discovered all the rooms), {\bf familiar}\footnote{The emphasis is Grothendieck's.} at the same time, and, most of all, {\bf unchanging}.\footnote{La plupart des math\'ematiciens sont port\'es \`a se cantonner dans un cadre conceptuel, dans un ``{\it Univers}" fix\'e une fois pour toutes -- celui, essentiellement, qu'ils ont trouv\'e ``tout fait" au moment o\`u ils ont fait leurs \'etudes. Ils sont comme les h\'eritiers d'une grande et belle maison toute install\'ee, avec ses salles de s\'ejour et ses cuisines et ses ateliers, et sa batterie de cuisine et un outillage \`a tout venant, avec lequel il y a, ma foi, de quoi cuisiner et bricoler. Comment cette maison s'est construite progressivement, au cours des g\'en\'erations, et pourquoi et comment ont \'et\'e con\c cus et fa\c conn\'es tels outils (et pas d'autres...), pourquoi les pi\`eces sont am\'enag\'ees de telle fa\c con ici, et de telle autre l\`a -- voil\`a autant de questions que ces h\'eritiers ne songeraient pas \`a se demander jamais. C'est ça ``l'Univers", le ``donn\'e" dans lequel il faut vivre, un point c'est tout ! Quelque chose qui para\^\i t grand (et on est loin, le plus souvent, d'avoir fait le tour de toutes ses pi\`eces), mais {\bf familier} en m\^eme temps, et surtout: {\bf immuable}.}
\end{quote}

We conclude with Grothendieck's reference to Riemann. In his  \emph{Sketch of a program}, \cite{Gro-esquisse}, he writes (p. 240 of the English translation):\footnote{The English translation is by Lochak and Schneps.}
\begin{quote}\small
Whereas in my research before 1970, my attention was systematically directed towards objects of maximal generality, in order to uncover a general language adequate for the world of algebraic geometry [...] here I was brought back, via objects so simple that a child learns them while playing, to the beginnings and origins of algebraic geometry, familiar to Riemann and his followers!
\end{quote}

\noindent \emph{Acknowledgements.--- } I would like to thank Vincent Alberge, Jeremy Gray and Marie-Pascale Hautefeuille who read a preliminary version of this paper and suggested corrections.
  \newpage
\smaller\smaller\smaller\smaller
  \begin{center}

\rowcolors{1}{gray!50!gray!50}{gray!50}
 \begin{tabular}{ | l | l | l   | }
 \hline

 {\bf Topic} & {\bf Euler} & {\bf Riemann}\\

 \hline
   \hline

      Functions of a  & $\bullet$  Introductio in  analysin   & $\bullet$  Grundlagen f\"ur eine  \\

 complex variable &  infinitorum (1748) &  allgemeine Theorie der  \\

  &  $\bullet$ De repraesentatione  &   Functionen einer\ ver\"anderlichen \\

  &   superficiei sphaericae   & complexen Gr\"osse (1851)  \\

  & super plano (1777) &  $\bullet$ Theorie der Abel'schen \\

  &  & Functionen (1857)  \\

    \hline
    Elliptic and & $\bullet$ Specimen de constructione  & $\bullet$  Grundlagen f\"ur eine \\
   Abelian integrals & aequationum differentialium & allgemeine Theorie der  \\
  & sine   indeterminatarum  &  Functionen einer\ ver\"anderlichen \\
  &  separatione (1738)  & complexen Gr\"osse (1851)  \\
  
    & $\bullet$ Observationes de comparatione  &   $\bullet$ Theorie der Abel'schen  \\
   & arcuum curvarum &  Functionen (1857)  \\  
  &  irrectificibilium (1761) &$\bullet$  \"Uber das Verschwinden der  \\

    & $\bullet$ De integratione  &   $\vartheta$-Functionen (1857)  \\
   & aequationis differentialis &  \\
  &  $ \frac{mdx}{\sqrt{1-x^4}}=\frac{ndy}{\sqrt{1-y^4}}$  (1761)&   \\

    \hline
Hypergeometric  &   $\bullet$  De summatione innumerabilium  & $\bullet$  Beitr\"age zur Theorie der  \\
 series &   progressionum (1738)  & durch  die Gauss'sche Reihe  \\
 & $\bullet$  Institutionum calculi integralis &  $F(\alpha,\beta,\gamma,x)$ darstellbaren   \\
  
  & volumen secundum (1769) & Functionen (1857)   \\
 &  $\bullet$ Specimen transformationis &   \\
 &  singularis serierum (1778) &    \\

 \hline

    The zeta function & $\bullet$ Variae observationes circa  & $\bullet$ \"Uber die Anzahl der   \\

  & series infinitas (1744)  &  Primzahlen unter einer   \\
    &  $\bullet$ Remarques sur un beau rapport & gegebenen Gr\"osse (1859)  \\
  & entre les series des puissances tant   &   \\
   & directes que r\'eciproques (1749)  &   \\

 \hline
   Integration  & $\bullet$   Institutionum calculi integralis  & $\bullet$ \"Uber die Darstellbarkeit \\
         
  &  (1768--1770) & einer Function durch eine   \\
      
    &    & trigonometrische Reihe (1854)  \\

     \hline
    Space and & $\bullet$   Anleitung zur Naturlehre   & $\bullet$ Grundlagen f\"ur eine    \\
 philosophy   of  &   (1745) &   allgemeine Theorie der  \\
    nature  &  $\bullet$ Reflexions sur l'espace   &  Functionen einer\ ver\"anderlichen   \\
         &    et le temps (1748)  &   complexen Gr\"osse (1851)  \\

     \hline
     Topology  & $\bullet$ Solutio problematis ad    & $\bullet$ Grundlagen f\"ur eine    \\
   &  geometriam situs pertinentis  &   allgemeine Theorie der  \\
      &   (1741) &  Functionen einer\ ver\"anderlichen   \\
         &   $\bullet$  Elementa doctrinae   &   complexen Gr\"osse (1851)  \\
            &  solidorum (1758)   &   $\bullet$ Theorie der Abel'schen \\
               &   $\bullet$ Demonstratio nonnullarum   &   Functionen (1857)  \\
                  &    insignium proprietatum,   &  $\bullet$ \"Uber die Hypothesen,    \\
                     & quibus  solida  hedris  planis    &  welche der Geometrie   \\
                        &    inclusasunt praedita   (1758) & zu Grunde liegen (1854)   \\

    \hline

 Differential  &   $\bullet$  Introductio in  analysin  & $\bullet$  \"Uber die Hypothesen,   \\
  geometry &   infinitorum (1748) & welche der Geometrie  \\
 &   & zu Grunde liegen (1854)   \\
  
  &  $\bullet$ Recherches sur la courbure & $\bullet$  Commentatio mathematica,    \\
 &   des surfaces (1767) & qua  respondere tentatur \\
 &   &  quaestioni  ab Ill${}^{\mathrm{ma}}$   Academia \\
 &   &  Parisiensi propositae (1861)  \\

   &   &  $\bullet$ Ein 
   beitrag zu den \\
            
               &   & Untersuchungen \"uber die    \\
                 &   &   fl\"ussigen Bewegung eines  \\
               
                  &   &gleichartigen Ellipsoides (1861)  \\

   &   &  $\bullet$  \"Uber die Fl\"ache vom  \\
   
      &   &  kleinsten Inhalt bei gegebener \\
      
         &   &  Begrenzung (1867) \\

    \hline

    Trigonometric  & $\bullet$  Recherches sur la question   & $\bullet$ \"Uber die Darstellbarkeit \\
         
      series  & des  inegalit\'es  du mouvement  & einer Function durch eine   \\
      
    &   de Saturne   et de Jupiter (1748)  & trigonometrische Reihe (1854)  \\

   \hline

       Acoustics & $\bullet$ Dissertatio physica de sono &  $\bullet$ \"Uber die Fortpflanzung  \\
  &  (1727) &  ebener Luftwellen von endlicher  \\
  & $\bullet$ Sur la vibration des cordes  &  Schwingungsweite  (1860) \\
  & (1748)   & \\
  
    \hline

 \end{tabular}

\end{center}

%\newpage 
%            \medskip
%            \begin{center}
%              \begin{figure}
%\centering
%\includegraphics[width=1.05\linewidth]{gravure-Euler.pdf}    \caption{\small {Leonhard Euler}}   \label{Euler}  
%\end{figure}
%
%\end{center}
% 

\printindex
  

\begin{thebibliography}{99}
\bibitem{Wallis2} J. Wallis, A treatise on algebra, both historical and practical, shewing, the original, progress, and advancement thereof, from time to time, and by what steps it hath attained to the heighth at which now it is. With some additional treatises, I. Of the cono-cuneus; being a body representing in part a conus, in part a cuneus. II. Of angular sections; and other things relating there-unto, and to trigonometry. III. Of the angle of contact; with other things appertaining to the composition of magnitudes, the inceptives of magnitudes, and the composition of motions, with the Results thereof. IV. Of combinations, alternations, and aliquot parts.London, John Playford for Richard Davis, 1685.
\end{thebibliography}
 \end{document}